\documentclass[11pt]{amsart}
\usepackage{amsmath,amsthm,amssymb,latexsym,color}
\usepackage{hyperref}
\usepackage{lipsum}

\date{}

\usepackage{graphicx}
\newtheorem{theorem}{Theorem}[section]
\newtheorem*{theorem*}{Theorem}
\newtheorem*{theoremA}{Theorem~A}
\newtheorem*{theoremB}{Theorem~B}

\newtheorem{lemma}[theorem]{Lemma}
\newtheorem{cor}[theorem]{Corollary}
\newtheorem*{corC}{Corollary~C}

\newtheorem{prop}[theorem]{Proposition}

\newtheorem{claim}[theorem]{Claim}
\theoremstyle{definition}
\newtheorem{Remark}[theorem]{Remark}

\theoremstyle{plain}

\newcommand{\N}{\mathbb{N}}
\newcommand{\Z}{\mathbb{Z}}
\newcommand{\R}{\mathbb{R}}

\newcommand{\E}{\mathbb{E}}
\newcommand{\Var}{{\rm Var}}
\newcommand{\Event}{\mathcal{E}}
\def\Prob{{\mathbb P}}

\newcommand{\net}{\mathcal N}
\newcommand{\matrixset}{\mathcal M}

\def\Row{{\rm row}}

\newcommand{\degree}{{\rm deg}}
\newcommand{\edg}{\leftrightarrow}
\newcommand{\inc}{{\rm \mathbf{Inc}}}
\newcommand{\cycle}{\mathcal C}

\newcommand{\diam}{\mathbf{Diam}}
\newcommand{\tree}{\mathbb{T}}

\newcommand{\dmax}{d_{\max}}

\def\path{{\mathcal P}}
\def\disc{{\bf Discov}}
\def\uncov{{\bf Uncov}}
\def\nbr{{\mathcal N}}
\def\V{{\mathcal V}}
\newcommand{\w}{\widetilde}
\def\up{\uparrow}
\def\down{\downarrow}

\def\ind{{{\mathcal S}\text{-}{\bf Ind}}}
\def\cind{{\cycle\text{-}{\bf Ind}}}
\def\dind{{HV\text{-}{\bf Ind}}}

\newcommand{\Hset}{\mathcal{H}}

\newcommand{\hdset}{{\mathcal B}}
\newcommand{\height}{H}
\newcommand{\weight}{{\mathcal W}}
\newcommand{\cset}{{\mathcal A}}

\def\rearr{{\mathcal R}}
\def\classif{{\bf\zeta}}
\def\pathw{{\Psi}}
\def\stnd{{\mathcal Y}}
\def\badptclass{{\mathcal B\mathcal C}}

\def\cycleone{{\mathcal C_1}}
\def\cycleedges{{\mathcal C}}
\def\maxnorm{{\mathcal M_1}}
\def\multone{{\bf m_1}}
\def\multon{{\bf p_1}}

\title{Outliers in spectrum of sparse Wigner matrices}
\author{Konstantin Tikhomirov and Pierre Youssef}

\begin{document}

\begin{abstract}
In this paper, we study the effect of sparsity on the appearance of outliers in the semi-circular law.  
Let $(W_n)_{n=1}^\infty$ be a sequence of random symmetric matrices such that
each $W_n$ is $n\times n$ with i.i.d entries above and on the main diagonal equidistributed with the product
$b_n\xi$, where $\xi$ is a real centered uniformly bounded random variable of unit variance and
$b_n$ is an independent Bernoulli random variable with a probability of success $p_n$.
Assuming that $\lim\limits_{n\to\infty}n p_n=\infty$, we show that for the random sequence $(\rho_n)_{n=1}^\infty$
given by
$$\rho_n:=\theta_n+\frac{n p_n}{\theta_n},\quad \theta_n:=\sqrt{\max\big(\max\limits_{i\leq n}\|\Row_i(W_n)\|_2^2-np_n,n p_n\big)},$$
the ratio
$\frac{\|W_n\|}{\rho_n}$ converges to one in probability. 
A non-centered counterpart of the theorem allows to obtain
asymptotic expressions for eigenvalues of the Erd\H os--Renyi graphs,
which were unknown in the regime $n p_n=\Theta(\log n)$.
In particular,
denoting by
$A_n$ the adjacency matrix of the Erd\H os--Renyi graph
$\mathcal{G}(n,p_n)$ and by $\lambda_{|k|}(A_n)$ its $k$-th largest (by the absolute value) eigenvalue, 
under the assumptions $\lim\limits_{n\to\infty }n p_n=\infty$ and $\lim\limits_{n\to\infty}p_n=0$ we have
\begin{itemize}
\item (No non-trivial outliers) If $\liminf\frac{n p_n}{\log n}\geq\frac{1}{\log (4/e)}$
then for any fixed $k\geq2$, $\frac{|\lambda_{|k|}(A_n)|}{2\sqrt{n p_n}}$ converges to $1$ in probability;

\item (Outliers) If $\limsup\frac{n p_n}{\log n}<\frac{1}{\log (4/e)}$ then there is $\varepsilon>0$ such that
for any $k\in\N$, we have

$\lim\limits_{n\to\infty}\Prob\Big\{\frac{|\lambda_{|k|}(A_n)|}{2\sqrt{n p_n}}>1+\varepsilon\Big\}=1$.

\end{itemize}
On a conceptual level, our result reveals similarities in
appearance of outliers in spectrum of sparse matrices
and the so-called BBP phase transition phenomenon in deformed Wigner matrices.
\end{abstract}

\maketitle

\tableofcontents

\section{Introduction}

Spectral analysis of large random matrices is a very active area of research motivated by questions
in statistics, mathematical physics, computer science.
A quantity of particular interest is the empirical spectral distribution. 
Given an 
$n\times n$ symmetric matrix $A$, its {\it empirical spectral distribution} is a measure on $\R$ defined by $$\mu_A:=\frac{1}{n}\sum_{j=1}^n \delta_{\lambda_j(A)},$$ where $\lambda_j(A)$ denote the eigenvalues of $A$.

One of the classical results in the random matrix theory asserts
that whenever $(\Xi_n)_{n\geq 1}$ is a sequence of $n\times n$ symmetric matrices whose entries on and above the diagonal are independent
and equidistributed with a given random variable $\xi$ of zero mean and unit variance,
the sequence of (random) measures $\mu_{\frac{1}{\sqrt{n}} \Xi_n} $ converges almost surely to
the {\it Wigner semi-circular} distribution $\mu_{sc}$ with the density $\frac{1}{2\pi} \sqrt{4-x^2}\,\mathbf{1}_{[-2,2]}(x)$
\cite{wigner}; thus, the distribution of $\xi$ does not affect the limiting measure. 
Moreover, by considering the support of $\mu_{sc}$, it follows 
that almost surely $$\Big\Vert \frac{1}{\sqrt{n}} \Xi_n\Big\Vert\geq 2-o(1),$$ where $\Vert \cdot\Vert$ stands for the spectral norm and $o(1)$ denotes a quantity vanishing to $0$ as $n\to\infty$.
Whenever the entries of the matrix have a finite fourth moment,
the extreme eigenvalues converge to the edges of the support of the limiting measure \cite{GH, G, bai-yin}:
one has almost surely $\big\Vert \frac{1}{\sqrt{n}} \Xi_n\big\Vert\leq 2+o(1)$.
These relations determine the location of the spectrum on the macroscopic scale and
show, in particular, that under the fourth moment assumption
there are no spectral outliers (i.e.\ eigenvalues asymptotically detached from the support of the limiting measure).

In this paper, we study the effect of sparsity on the existence of spectral outliers.
We start with an $n\times n$ symmetric random matrix $\Xi_n=(\xi_{ij})_{1\leq i,j\leq n}$ as above and suppose that its entries are uniformly bounded. 
Next, 
we randomly
zero out some of the matrix entries. 
To implement this, let $B_n= (b_{ij})_{1\leq i,j\leq n}$ be an $n\times n$ symmetric matrix whose entries on and above the diagonal are i.i.d Bernoulli variables with probability of success $p_n$ 
and suppose that $B_n$ and $\Xi_n$ are independent. We consider the random matrix $W_n$ obtained as the entry-wise product of $B_n$ and $\Xi_n$. 
It is known that the Wigner semi-circular law is stable under the
sparsification as long as the average number of the non-zero entries in each row is infinitely
large. More precisely, as long as $np_n\to \infty$, we have 
$$
\mu_{\frac{1}{\sqrt{np_n}} W_n} \underset{n\to\infty}{\overset{a.s.}{\longrightarrow}} \mu_{sc}.
$$
The situation with spectral outliers is more complicated. 
Whenever $\frac{np_n}{\log n}\to \infty$, it is known that $\|W_n\|=(2+o(1))\sqrt{n p_n}$ with probability tending to one with $n$. This result
was verified in a series of works where the assumptions on the matrix sparsity
and the matrix entries were sequentially relaxed (see \cite{Furedi-komlos, Khorunzhy, Vu, BGBK, LVY}). 
On the other hand, when $p_n\to 0$ with $n$ relatively fast, one can easily verify that the extreme eigenvalues get asymptotically detached from the bulk of the spectrum. 
For example, taking $\xi_{ij}$ to be standard Rademacher variables and taking $p_n$ sufficiently small (say $p_n= \log \log n/n$), standard estimates on the tails of binomial random variables show 
that with probability going to $1$ with $n$
$$
\max_{1\leq i\leq n} \Vert \Row_i(W_n)\Vert_2\geq \sqrt{\frac{\log n }{\log \log n}},
$$
where we denoted by $\Row_i(W_n)$ the $i$-th row of $W_n$ and by $\Vert\cdot\Vert_2$ the Euclidean norm in $\R^n$. Since deterministically $\Vert W_n\Vert \geq \max\limits_{1\leq i\leq n} \Vert \Row_i(W_n)\Vert_2$, this indicates that when $p_n= \log \log n/n$, the extreme eigenvalue(s) do not converge to the edges of the support of the limiting measure. 
More generally, when $\frac{np_n}{\log n}\to 0$, this phenomenon was observed in the case of Rademacher variables \cite{Khorunzhy} and in the case of the Erd\H{o}s--Renyi random graphs \cite{BGBK sparse} which will be discussed later on. In the window around $\log n$
it is known that, up to constant multiples, the matrix norm is of order $\sqrt{\log n}$ \cite{seginer, BGBK sparse, BGBK, LVY},
however, to the best of our knowledge, there have been no results on its exact asymptotic behavior.
In this connection, we can ask the following questions: 
\vskip 0.3cm
{\it 
\begin{enumerate}
\item Is there a sharp phase transition (in terms of sparsity) in the appearance/disappearance of outliers in the semi-circular law? \\
\item For concrete distributions, say, sparse Bernoulli matrices,
what is the explicit formula for the sparsity threshold (if it exists)? \\
\item What is a conceptual explanation of why the outliers appear at a particular level of sparsity? \\
\item What is the exact asymptotic value of an outlier? 
\end{enumerate}
}
\vskip 0.3 cm

In this paper, we partially answer the above questions
by characterizing the norm (and, more generally, $k$--th largest eigenvalue)
of a sparse matrix. The first main result of this paper is the following theorem.

\begin{figure}
\begin{center} 
\includegraphics[width=0.5\textwidth]{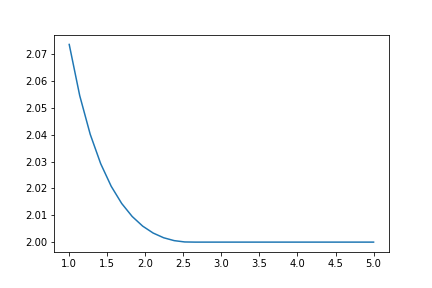}
\caption{\small \sl The value of $\lim_n(\rho_n/ \sqrt{np_n})$ (vertical axis) for the adjacency matrix of the
Erd\H os--Renyi graph,
viewed as a function of $\lim_n(np_n/\log n)$ (horizontal axis). 
The phase transition happens at 
$\frac{1}{\log(4/e)}\approx 2.59$. \label{figure}}
\end{center} 
\end{figure}

\begin{theoremA}
Let $\xi$ be a real centered uniformly bounded random variable of unit variance.
For each $n$, let $W_n$ be an $n\times n$ symmetric random matrix with i.i.d.\ entries above and on the main diagonal,
with each entry equidistributed with the product $b_n\xi$, where $b_n$ is a $0/1$ (Bernoulli) random variable independent of $\xi$,
with probability of success equal to $p_n$.
Assume further that $n p_n\to\infty$ with $n$. For each $n$, define the random quantities
$$\rho_n:=\theta_n+\frac{n p_n}{\theta_n},\quad \theta_n:=\sqrt{\max\big(\max\limits_{i\leq n}\|\Row_i(W_n)\|_2^2-np_n,n p_n\big)}.$$
Then the sequence $\big(\frac{\|W_n\|}{\rho_n}\big)_{n\geq 1}$ converges to one in probability.
More generally, denoting by $\lambda_{|k|}(W_n)$ the $k$-th largest (by the absolute value) eigenvalue
of $W_n$, for any fixed $k$ the sequence $\Big(\frac{|\lambda_{|k|}(W_n)|}{\rho_n}\Big)_{n\geq 1}$ converges to one in probability.
\end{theoremA}
The theorem is obtained as a combination of Theorems~\ref{th: centered upper} and~\ref{th: lower top} of this paper.
Let us make a few remarks. 
The quantity $\rho_n/\sqrt{np_n}$ is equal to $2$ 
iff $\max_{i} \Vert \Row_i(W_n)\Vert_2^2 \leq 2np_n$.
Combined with standard concentration inequalities and simple continuity properties of $\rho_n$, this implies
that if $(p_n)$ is a sequence satisfying
$$\limsup\limits_n\frac{\E\max_{i} \Vert \Row_i(W_n)\Vert_2^2}{n p_n}\leq 2$$
then there are no asymptotic spectral outliers for the sequence of matrices $(W_n)$, i.e.\
$\|W_n\|\leq (2+o(1))\sqrt{n p_n}$ with probability tending to one with $n$.
On the other hand, if
$$\liminf\limits_n\frac{\E\max_{i} \Vert \Row_i(W_n)\Vert_2^2}{n p_n}> 2$$
then there is $\varepsilon>0$ such that for any fixed $k$, $|\lambda_{|k|}(W_n)|\geq (2+\varepsilon)\sqrt{n p_n}$
with probability going to one.
We will revisit this statement in context of the Erd\H{o}s--Renyi graphs (see Corollary~C below and Figure~\ref{figure}).

Futher, let us discuss the result at a more conceptual level. 
For a fixed $k$, large enough $n$ and under the assumption $\max\limits_{i\leq n} \Vert \Row_i(W_n)\Vert_2^2 \geq 2np_n$, the $k$--th largest eigenvalue
$\lambda_{|k|}\big(\frac{1}{\sqrt{np_n}} W_n\big)$ is of order
$$
\frac{\sqrt{\max_{i}\|\Row_i(W_n)\|_2^2-np_n}}{\sqrt{np_n}}+ \frac{\sqrt{np_n}}{\sqrt{\max_{i}\|\Row_i(W_n)\|_2^2-np_n}}.
$$ 
A very similar formula has appeared multiple times in a different context --- in the study of perturbed random matrices.
The spectrum of random matrices perturbed by fixed matrices of a given structure 
has been subject of very active research. More specifically,
consider the spectrum of $M_n+H_n$, where $M_n$ is an $n\times n$ Wigner matrix and $H_n$ is a fixed deterministic symmetric perturbation. When $H_n$ is of a finite (or relatively small) rank, the limiting spectral distribution of $M_n+H_n$
is not affected by the perturbation (remains semi-circular) 
due to the interlacing property of the eigenvalues.
However, the perturbation can affect largest eigenvalues forcing some of them
to get asymptotically detached from the rest of the spectrum.
This phenomenon was first considered in \cite{Furedi-komlos} where the authors, motivated by estimating the largest eigenvalue of the adjacency matrix of an Erd\H{o}s--Renyi graph, studied rank one deformations of a Wigner matrix.
Later on, considerable interest in deformed random matrices was also connected with the work \cite{BBP},
where the famous BBP phase transition phenomenon was put forward.
A large number of articles was devoted to investigating the phase transition in a variety of models as well as to studying fluctuations of the largest eigenvalues detached or not detached from the bulk \cite{BBP, BC, BGM, BR,BS,C,CDF,CDF2,CDF3,FP,KY,KY1,OM,peche,PRS,Ra,RS,tao}. We refer, among others, to survey \cite{peche-icm} for a review of the subject. 
\begin{theorem*}[{BBP Phase Transition, see \cite[Theorem~2.1]{peche-icm}}]
For each $n$, let $M_n$ be an $n\times n$ Wigner matrix whose entries on and above the main diagonal are
independent copies of $\frac{1}{\sqrt{n}}\xi$, where $\xi$ is centered variable of unit variance.
Suppose further that $\xi$ has a finite fourth moment. 
Fix $r\in \N$ and $\theta_1\geq \ldots\geq\theta_r>0$,
and for each $n$, let $H_n$ be an $n\times n$ deterministic symmetric matrix of rank $r$ with non-zero eigenvalues $\theta_1,\dots,\theta_r$.
Then for any $1\leq i\leq r$, 
\begin{itemize}
\item If $\theta_i\leq 1$, then $\lambda_i\underset{n\to\infty}{\overset{a.s.}{\longrightarrow}} 2$;
\item If $\theta_i>1$, then $\lambda_i\underset{n\to\infty}{\overset{a.s.}{\longrightarrow}} \theta_i+\frac{1}{\theta_i}$,
\end{itemize}
where $\lambda_1\geq \ldots\geq \lambda_n$ denote the eigenvalues of $M_n+H_n$. 
\end{theorem*}

Note that if $\theta_i$ in the above theorem was replaced with
$\sqrt{(n p_n)^{-1}\max_{i}\|\Row_i(W_n)\|_2^2-1}$ then the theorem would
describe exactly the same asymptotic behavior
as revealed in Theorem~A. 
We can give the following {\it non-rigorous} justification for this similarity.
Let us reconsider the matrix $W_n$ from Theorem~A, and note that due to concentration inequalities,
most of the rows have their norms squared concentrated around $np_n$. 
Only a small fraction of these rows can have their norm far from $\sqrt{n p_n}$.
For simplicity, suppose that only one row of $W_n$ (say, the first one after an appropriate permutation)
has the Euclidean norm significantly larger than $\sqrt{n p_n}$.
Then we may decompose our matrix as 
$$
\frac{1}{\sqrt{np_n}} W_n\approx M_n+ H_n,
$$
where $M_n$ is obtained from $\frac{1}{\sqrt{np_n}}W_n$ by a {\it regularization} procedure of reducing the entries of the first row
and column
(the ones with the largest norms) in such a way
that the Euclidean norm of the transformed row and column is equal $1$,
and $H_n$ is the remainder: the $n\times n$ symmetric zero diagonal matrix whose first row/column's 
Euclidean norm is equal to $u_n:=\sqrt{(n p_n)^{-1}\|\Row_1(W_n)\|_2^2-1}
\approx\sqrt{(n p_n)^{-1}\max_{i}\|\Row_i(W_n)\|_2^2-1}$,
and the matrix entries not in the first row or column of $H_n$ are zeros
(here, we use ``$\approx$'' instead of the equality sign to emphasize that
our model describes the actual distribution of $W_n$ only approximately).
In a sense, we treat the extra mass in the rows and columns of $W_n$ of large Euclidean norms
as a deformation of the regularized matrix $M_n$.
We take this extra mass from the first row and column of $\frac{1}{\sqrt{np_n}}W_n$ and 
transfer it to the matrix $H_n$ which is perceived as a perturbation of $M_n$.
Clearly, $H_n$ is of rank $2$ with eigenvalues $\pm u_n$. On the other hand, by our assumption,
all rows of $M_n$ have their norms concentrated around $\sqrt{n p_n}$
which could suggest that the spectrum of $M_n$ has no outliers. 
Now, Theorem~A  states that if $u_n\leq 1$, then $\frac{1}{\sqrt{np_n}} W_n\approx M_n+H_n$
has all its eigenvalues asymptotically bounded by $2$; whereas, if $\liminf_n u_n>1$, then $M_n+H_n$
has an outlier and its value is given by $u_n+\frac{1}{u_n}$.
This parallels the BBP phase transition phenomenon. 
Let us emphasize once more that the above discussion
is meant only to suggest similarities between the two models and is not developed rigirously.
It seems interesting to understand if such a connection could be elaborated.
We remark that there has been several works recently
concerned with {\it regularizations} of random graphs/matrices
i.e procedures designed to reduce the norm of random matrices by changing a few of its entries (see \cite{FO,LLV,RV,Re}).

\begin{figure}
\begin{center} 
\includegraphics[width=0.5\textwidth]{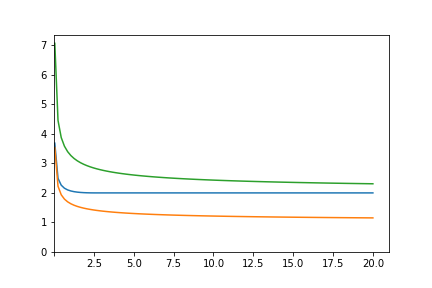}
\caption{\small \sl This figure illustrates the relation between the maximum row norm
and the quantity $\rho_n$ in the setting of the Erd\H os--Renyi graphs $\mathcal{G}(n,p_n)$ with adjacency matrices $A_n$,
with $\lim\limits_{n\to\infty }n p_n=\infty$ and $\lim\limits_{n\to\infty}p_n=0$.
The horizontal axis is the value of the limit $\lim_n \frac{n p_n}{\log n}$.
The blue curve is the corresponding values of $\lim_n\,(\rho_n/\sqrt{np_n})$ (convergence in probability).
The orange curve is the values of $\lim_n\,\sqrt{(n p_n)^{-1}\max_{i\leq n} \degree(i)}$.
The green curve --- the values of $\lim_n \,2\sqrt{(n p_n)^{-1}\max_{i\leq n} \degree(i)}$. 
When $\frac{np_n}{\log n}\to 0$, the left side of inequality \eqref{eq: norm-maxnorm} (with $W_n$
replaced with $A_n-\E A_n$)
is asymptotically sharp while for
$\frac{np_n}{\log n}\to \infty$, the right side of \eqref{eq: norm-maxnorm} is sharp. 
\label{fig}
}
\end{center} 
\end{figure}

Another consequence of Theorem~A is that the operator norm of the matrix
$W_n$ has the same order of magnitude as $\max_i \Vert \Row_i(W_n)\Vert_2$.
As previously stated, one has deterministically $\max_i \Vert \Row_i(W_n)\Vert_2\leq \Vert W_n\Vert$,
and Theorem~A implies that the reverse inequality is true up to a universal constant. Indeed, an analysis of the parameter $\rho_n$ shows that for any $\varepsilon>0$, $$\rho_n\leq (2+\varepsilon)\max\limits_{i\leq n}\|\Row_i(W_n)\|_2,$$ 
with probability going to $1$ with $n$. To view this, note that when $\max_{i}\|\Row_i(W_n)\|_2^2\leq 2np_n$, we have $\rho_n=2\sqrt{np_n}$ while otherwise $\rho_n\leq \sqrt{2}  \max_{i}\|\Row_i(W_n)\|_2$.
Moreover, in view of standard concentration inequalities, for any $\varepsilon>0$ one has 
$\max\limits_{i\leq n}\|\Row_i(W_n)\|_2\geq (1-\varepsilon)\sqrt{np_n}$ with probability going to one with $n$. Therefore, Theorem~A implies that for any $\varepsilon>0$, we have 
\begin{equation}\label{eq: norm-maxnorm}
\max\limits_{i\leq n} \Vert \Row_i(W_n)\Vert_2\leq \Vert W_n\Vert\leq (2+\varepsilon) \max\limits_{i\leq n} \Vert \Row_i(W_n)\Vert_2,
\end{equation}
with probability going to one with $n$ (see Corolllary~C and Figure~\ref{fig} for the case of the Erd\H{o}s--Renyi graphs).
This phenomenon was first observed by Seginer \cite{seginer} who showed that for matrices with i.i.d entries, the operator norm is comparable, up to a constant, to the maximum Euclidean norm of its rows/columns. We note that Seginer's result applies to a much wider
class of distributions of the entries (as long as they are i.i.d) and while it is stated in \cite{seginer}
for expectations and for non-symmetric matrices, it is not difficult to obtain its extension to tail estimates
for norms of symmetric matrices.
Theorem~A recovers Seginer's observation in the setting of uniformly bounded entries
and gives asymptotically optimal
relation between the spectral norm and the maximum Euclidean norm of the rows.
While one might be tempted to think that this comparison is valid for any random matrix with independent entries of different
variances, Seginer \cite{seginer} provided an example showing that it is not the case even 
for the class of inhomogeneous matrices with subgaussian entries.
Nevertheless, it was shown in \cite{LVY} that for inhomogeneous Gaussian matrices with independent Gaussian entries
having arbitrary variances, the spectral norm is equivalent, up to constant multiples,
to the maximum Euclidean norm of rows. We refer to \cite{LVY} for further discussion and references concerning this phenomenon. 

Theorem~A deals with sparse matrices with centered entries and does not directly provide information on the magnitude of the largest
eigenvalues in the non-centered setting.
Assume $W_n$ is an $n\times n$ symmetric matrix whose entries on and above the main diagonal are i.i.d copies of $b_n\xi$,
where $b_n$ is a $0/1$ Bernoulli with probability of success $p_n$ and $\xi$ is a uniformly bounded real random variable of unit second moment (not necessarily centered) independent of $b_n$.
Naturally, one could recenter the matrix $W_n$ and consider the matrix $W_n-\E\, W_n$
in order to estimate $\lambda_{|2|}(W_n)$.
However, the centered matrix is no longer sparse and Theorem~A cannot be applied.
Moreover, standard symmetrization technique replacing $W_n-\E\, W_n$
with a difference of two independent copies of $W_n$, would result in extra multiplicative constants. Indeed, one can write 
$$
\E\Vert W_n-\E\, W_n\Vert \leq \E\, \Vert W_n -\w W_n\Vert,
$$
where the matrix $W_n-\w W_n$ has entries of the form $b_{ij}(\xi_{ij}-\xi_{ij}')$, 
with $\xi_{ij}'$ being independent copies of $\xi_{ij}$. This matrix is sparse and has centered entries so that Theorem~A can be applied. However, the rows of $W_n-\w W_n$ have on average by $\sqrt{2}$ larger Euclidean norms than rows of $W_n$,
resulting in an extra constant factor in the upper bound for $\lambda_{|k|}(W_n)$ obtained by this procedure.
In order to capture the true asymptotic behavior of $\lambda_{|k|}(W_n)$, we develop a special 
procedure relating the spectrum of the non-centered matrix to a specially chosen centered model.
This reduction will be discussed in more detail in the next section. Let us state the second main result of this paper. 

\begin{theoremB}
Let $\xi$ be a uniformly bounded real random variable with $\E\xi^2=1$.
For each $n$, let $W_n$ be an $n\times n$ symmetric random matrix with i.i.d.\ entries above and on the main diagonal,
with each entry equidistributed with the product $b_n\xi$, where $b_n$ is $0/1$ (Bernoulli) random variable independent of $\xi$,
with probability of success equal to $p_n$.
Assume further that $n p_n\to\infty$ with $n$ and $\lim\limits_{n\to\infty}p_n=0$. Then,
defining $\rho_n$ as in Theorem~A, i.e
$$\rho_n:=\theta_n+\frac{n p_n}{\theta_n},\quad \theta_n:=\sqrt{\max\big(\max\limits_{i\leq n}\|\Row_i(W_n)\|_2^2-np_n,n p_n\big)},$$
for each fixed $k\geq 2$, the sequence $\Big(\frac{|\lambda_{|k|}(W_n)|}{\rho_n}\Big)_{n\geq 1}$ converges to one in probability,
where $\lambda_{|k|}(W_n)$ denotes the $k$-th largest (by the absolute value)
eigenvalue of $W_n$.
\end{theoremB}
Theorem~B is obtained as a combination of Theorems~\ref{th: non-centered upper} and~\ref{th: lower top} of the paper.
The main application of Theorem~B concerns the random Erd\H os--Renyi graphs $G_n:=\mathcal{G}(n,p_n)$, by taking
$\xi$ to be constant $1$.
In \cite{KS}, it is shown that the largest eigenvalue of $G_n$  
almost surely satisfies 
$$
\lambda_1(G_n)=\big(1+o(1)\big) \max\Big( \sqrt{\max\nolimits_{i\leq n}\degree(i)}, np_n\Big),
$$
where $o(1)$ tends to $0$ as $ \max\big( \sqrt{\max_{i\leq n}\degree(i)}, np_n\big)$ tends to infinity,
and $\degree(i)$ is the degree of the $i$-th vertex of $G_n$.
Of particular interest is the second eigenvalue $\lambda_{2}(G_n)$ as the difference
$\lambda_1(G_n)-\lambda_{2}(G_n)$ may be viewed as a measure of the graph expansion properties.
As stated previously, when $\frac{np_n}{\log n}\to \infty$,
it is known that $\vert \lambda_{\vert 2\vert}(G_n) \vert \leq (2+o(1))\sqrt{np_n}$ \cite{Furedi-komlos, Vu, BGBK, LVY}. 
With further constraints on $n p_n$,
more precise information, including fluctuation intervals and limiting distribution of the extreme eigenvalues
is available in the literature \cite{EKYY1,EKYY,LS, HKM, HLY}. 
In contrast, when $\frac{np_n}{\log n}\to 0$, it is shown in \cite{BGBK sparse}
that $\vert \lambda_{\vert 2\vert}(G_n)\vert$ concentrates around the square root of the maximum degree in the graph
(in the same paper, the authors study the $k$-th largest eigenvalue, for arbitrary $k\leq n^{1-\varepsilon}$ and $\varepsilon >0$). 

In the window $n p_n\approx\log n$, 
no asymptotically sharp results for $\vert \lambda_{\vert 2\vert}(G_n)\vert$
were previously available. Moreover, it was not even known if there is a {\it sparsity threshold} (a multiplicative factor of $\log n$)
where the phase transition between existence and absence of non-trivial outliers in the spertrum can be observed.

It follows from Theorem~B that $\Big(\frac{|\lambda_{|2|}(G_n)|}{\rho_n}\Big)_{n\geq 1}$ converges to $1$ in probability, where 
$$
\rho_n=\sqrt{\max\big(\max\nolimits_{i\leq n}\degree(i) -np_n,n p_n\big)}+ \frac{np_n}{\sqrt{\max\big(\max\nolimits_{i\leq n}\degree(i)-np_n,n p_n\big)}}.
$$
The distribution of the maximum degree of the
Erd\H{o}s--Renyi graph is very well understood
(see, for example, \cite[Theorem~3.1]{book-bollobas}).
This leads to an explicit formula for $\rho_n$
and thus an asymptotic formula for $\vert\lambda_{\vert 2\vert}\vert$. 
The phase transition for the Erd\H os--Renyi graphs is considered in the following statement.

\begin{corC}[Outliers in the spectrum of the Erd\H os--Renyi graphs]
For each $n\geq 1$, let $G_n$ be the Erd\H os--Renyi random graph on $n$ vertices,
with parameter $p_n$. Assume that $n p_n\to \infty$ and $p_n\to 0$. Let $\lambda_{|k|}(G_n)$
be the $k$-th largest by absolute value eigenvalue of the adjacency matrix of $G_n$.
Then, denoting
$$
\rho_n^G:=\theta_n^G+\frac{n p_n}{\theta_n^G},\quad \theta_n^G:=
\sqrt{\max\bigg(enp_n \exp\Big[\mathcal{W}_0\Big( \frac{\log n-np_n}{enp_n}\Big)\Big]-np_n,np_n\bigg)},
$$
for any $k\geq 2$ the ratio $\frac{|\lambda_{|k|}(G_n)|}{\rho_n^G}$
converges to one in probability.
In particular,
\begin{itemize}

\item (No non-trivial outliers) If $\liminf\frac{n p_n}{\log n}\geq\frac{1}{\log(4/e)}$
then for any $k\geq2$, $\frac{|\lambda_{|k|}(G_n)|}{2\sqrt{n p_n}}$ converges to
$1$ in probability.

\item (Outliers) If $\limsup\frac{n p_n}{\log n}<\frac{1}{\log(4/e)}$ then there is $\varepsilon>0$ such that
for any $k\in\N$, we have

$\lim\limits_{n\to\infty}\Prob\Big\{\frac{|\lambda_{|k|}(G_n)|}{2\sqrt{n p_n}}>1+\varepsilon\Big\}=1$.

\end{itemize}
Here, $\mathcal{W}_0$ denotes the main branch of the Lambert function defined
by $z=\mathcal{W}_0(z)e^{\mathcal{W}_0(z)}$.
\end{corC}

We will provide a proof of the corollary in Section~\ref{s: outliers-erdos-renyi}. 
The corollary is illustrated in Figure~\ref{figure}. 

\medskip

Shortly after this manuscript was posted on arXiv, related results 
appeared in the work \cite{ADK}. In particular, the phase transition above happening at $\log n/ \log (4/e)$ for the appearing of outliers in the spectrum of the Erd\H os--Renyi graphs 
was also captured in \cite{ADK} and the results of \cite{ADK} extend as well to the Wigner matrix model studied in this paper. 
The authors of \cite{ADK} apply a completely different technique which cleverly exploits a tridiagonal representation of a hermitian matrix and a relation with the spectrum of the associated non-backtracking matrix. 
We refer to \cite{ADK} for more details.

\subsection*{Acknowledgements} A part of this work was done
while the second named author was visiting Georgia Tech in July 2018.
He would like to thank the institution for the great working conditions. 
P.Y was supported by grant ANR-16-CE40-0024-01.

\section{Overview of the proof}\label{s: overview}

This section is intended to give a fairly detailed overview of the proofs of our main results,
giving an emphasis to those parts of the argument which, in our opinion,
may turn out useful in future works on the subject.
As a starting point, we consider a simplified model that shows how (and why) the quantity $\rho_n$
defined in the main theorems, appears in the proof.

As it was mentioned in the introduction, it seems instructive to think of our model as of a standard (dense)
Wigner matrix being perturbed by a small number of rows/columns of relatively large norms.
These rows and columns distort the matrix spectrum and (if the magnitude of the norms exceeds a certain threshold)
shift the largest eigenvalue to a non-classical location.

We will consider a simpler {\it deterministic} model as an illustration.
Assume that the entries of our symmetric $n\times n$ matrix $A_n$ take values $\{0,1\}$
and that the locations of non-zero elements are fixed (non-random), and that $G_n$ is the corresponding simple deterministic graph
which does not contain any cycles.
Assume further that the support length (i.e.\ the vertex degree)
of every row/column, except for the first one, is at most $d=d(n)$,
while the support of the first row/first column
has length $\widetilde d=\w d(n)\geq d$.
We will estimate from above the norm of $A_n$ using the trace method. Fix $k\geq 1$. A standard formula gives
$$
\|A_n\|^{2k}\leq \sum\limits_{\path}1,
$$
where the summation is taken over all closed paths $\path$ of length $2k$ on $G_n$.
For each vertex $v$ of the graph $G_n$, and its neigborhood $\nbr_v$, we fix a bijective mapping
$ind_v$ from $\nbr_v$ into the integer interval $[1,\deg(v)]$ --- a ``local indexation'' of neighbors of $v$.

First, we consider the paths starting at vertex $1$ (of degree $\w d$).
Each such path corresponds to a diagram $H_\path$ on $[0,2k]$ i.e.\ a mapping $H_\path:[0,2k]\to\Z$
with $H_\path(0)=0$ and, for every $t\in[2k]$, $H_\path(t)-H_\path(t-1)=1$ whenever $\path(t)$ is farther from vertex
$1$ than $\path(t-1)$,
and $H_\path(t)-H_\path(t-1)=-1$ otherwise. Note that the diagram
is a Dyck path i.e.\ it is non-negative everywhere and is equal to zero at $2k$.
To each moment of time $t\geq 1$ with $H_\path(t)-H_\path(t-1)=1$
we can put in correspondence the local index $ind_{\path(t-1)}(\path(t))$.
Then the data structure consisting of the diagram $H$ and the $k$ indices corresponding to times with $H_\path(t)-H_\path(t-1)=1$,
will uniquely identify the path, i.e.\ in order to estimate the number of paths it is sufficient to estimate
the number of such data structures.

Note that whenever $H_\path(t)-H_\path(t-1)=1$ and $H_\path(t-1)>0$, the corresponding index can only take values in $[1,d]$,
while in the case $H_\path(t-1)=0$ the index takes values in $[1,\w d]$.
For each $1\leq u\leq k$, let $N_u$ be the total number of the Dyck paths of length $2k$ with $u$ returns to zero (counting
the point $2k$). 
By a standard formula, $N_u=\frac{u}{2k-u}{2k-u\choose k}$. Thus, the total number of the data structures can be estimated
by
$$
\sum\limits_{u=1}^k N_u d^{k-u}{\w d}^u
= \sum\limits_{u=1}^k \frac{u}{2k-u}{2k-u\choose k} d^{k-u}{\w d}^u
\leq d^k\,\frac{\max(\w d/d,2)^{2k-1}}{\big(\max(\w d/d,2)-1\big)^{k-1}}
$$
(see Lemma~\ref{lem: calculation-binomial} of this paper for a proof of the last relation).

We omit computations related to the setting when the starting vertex of a path is not $1$;
the upper bound is essentially the same as above. Overall, assuming an appropriate growth condition for $k=k(n)$
(in particular, $\lim\limits_{n\to\infty}k/\log n=\infty$), we can show that
\begin{align*}
\|A_n\|^{2k}\leq \sum\limits_{\path}1&\leq (1+o(k))^k d^k\,\frac{\max(\w d/d,2)^{2k-1}}{\big(\max(\w d/d,2)-1\big)^{k-1}}\\
&\leq (1+o(k))^k \,\frac{\max(\w d,2d)^{2k}}{\big(\max(\w d,2d)-d\big)^{k}}\\
&=(1+o(k))^k \Bigg(\sqrt{\max(\w d-d,d)}+\frac{d}{\sqrt{\max(\w d-d,d)}}\Bigg)^{2k}.
\end{align*}
Taking into account that $\w d=\max\limits_{i\leq n}\|\Row_i(A_n)\|_2^2$,
the last expression in the brackets perfectly agrees with the definition of $\rho_n$ in the
theorems from the introduction.

\medskip

Proving that the above model accurately describes the situation in case of sparse Wigner matrices
is the main technical problem within the proof.
More specifically, we need to show that the norm of a typical realization of a sparse Wigner matrix
with uniformly bounded entries is essentially determined by local tree-like structures similar
to the one in the above example. In order to implement this strategy, we need to resolve a number of issues;
among them:
\begin{itemize}
\item Show that the contribution of paths with cycles is not much larger than the contribution of paths on trees.
\item Deal with the fact that there are multiple vertices of large degrees within the graph,
and rather than taking two distinct values the degrees are ``continuously'' distributed within some integer interval.
\item Develop a procedure to condition on a ``good'' realization of the matrix $A_n$ and the underlying graph $G_n$.
Clearly, in the sparsity regime we study, taking the unconditional expectation of $\|A_n\|^{2k}$
would result in highly suboptimal bound on the norm. While taking a conditional expectation
given a ``good'' realization of the graph $G_n$ and of absolute values of the matrix entries
may seem a reasonable strategy when the entries are symmetrically
distributed, in the case of non-symmetric distributions a different approach has to be used.
\item Transfer the results obtained for centered matrices to the non-centered setting, in particular,
the adjacency matrices of the Erd\H os--Renyi graphs.
Since the centered adjacency matrices are no longer sparse, this problem requires a special symmetrization
procedure.
\item Show that the upper bound obtained using this strategy is optimal i.e.\ prove a matching lower bound.
\end{itemize}

\bigskip

As a common starting point, given an $n\times n$ symmetric matrix $A$, we write
$$
\lambda_1(A)^{2k}\leq \sum\limits_{\path}\prod\limits_{e\sim\path}a_e,
$$
where the sum is taken over all closed paths $\path$ of length $2k$ on the complete graph $K_{[n]}$,
the product --- over all edges $e$ in $\path$ (counting multiplicities), and where $a_e$ is the
matrix entry corresponding to the edge $e$.
While an averaging argument (taking the expectation of the expressions on the left and right of the above relation)
is a usual step in classical applications of the trace method,
in this paper we rely on computing the expectations only when considering paths
with many edges of multiplicity one, whereas for other paths we estimate the products $\prod\limits_{e\sim\path}a_e$
for {\it every} realization of the matrix $A$ within some special event of probability close to one.

\medskip

{\bf Data structure.}
Classical applications of the trace method often involve defining
an auxiliary structure associated with a path, which simplifies counting;
for example, diagrams and auxiliary sets marking cyclic elements within the path.
In our proof, the data structure associated with a path plays a fundamental role,
and, in addition to ``usual'' information (times of discovering new vertices/edges, traveling directions along a previously
discovered edge) also contains data about the magnitude of the Euclidean norm and the distribution of mass 
across the rows/columns corresponding to the graph vertices.

At an abstract level, our approach can be described as follows:
we define an {\it injective} mapping ${\bf Data}$ from the set of paths $\path$ into a ``data space'' ${\bf S}$,
and for each element $s\in{\bf S}$ we define a weight $w(s)$ in such a way that
$w({\bf Data}(\path))\geq \prod\limits_{e\sim\path}a_e$ for all paths. Then, in view of the injectivity,
$$
\sum\limits_{\path}\prod\limits_{e\sim\path}a_e
\leq \sum\limits_{s\in{\bf S}}w(s).
$$
This way, analysis of the paths can be completely replaced by counting on the data space.
A crucial part of this approach is to define the mapping and the data space in such a way that,
on the one hand, ${\bf S}$ is sufficiently ``rich'' and both ${\bf S}$ and $w(\cdot)$ are
simply structured so that injectivity can be easily established
and the sum on the right hand side --- (relatively) easily computed;
on the other hand, ${\bf S}$ is not too large so that the sum $\sum\limits_{s\in{\bf S}}w(s)$
can be efficiently controlled from above.
A formal definition of our mapping ${\bf Data}$ and the proof of injectivity is given in Section~\ref{s: data-structure};
some structural properties of the data space ${\bf S}$ are discussed in Section~\ref{s: prop-data}. 
A satisfactory definition of the weight $w(\cdot)$ presents an issue on its own.
As we attempt to make counting over the data structures {\it simpler} than counting
over the paths, we inevitably lose information about the matrix when transfering the problem to the space ${\bf S}$;
in particular, we do not have in possession the precise information about the value of the product of the matrix entries
corresponding to a given data structure.
In order to deal with this issue, we introduce {\it vector majorizers}.

\medskip

{\bf Vector majorizers.}
A {\it majorizer} of a vector $x$ in $\R^n$ is any vector $y$ such that $y\geq x^*$
coordinate-wise, i.e.\ $y_i\geq x_i^*$ for all $i\leq n$, where $x^*$ denotes the non-increasing rearrangement of
the vector of absolute values of components of $x$.
The crucial observation, which we make in a rather general deterministic setting in Section~\ref{subs: majorizer1}
and in the more specific probabilistic setting in Section~\ref{subs: major},
is that there exists a {\it small} collection of majorizers $\net_{mjr}$ such that for any typical
realization of our matrix $A$, every row
can be majorized by a vector from this collection having only slightly larger norm.
Now, in order to implement the weight function $w(s)$ on our data space ${\bf S}$,
it is sufficient to record in each structure $s$ which majorizers from $\net_{mjr}$
have to be used, and then define $w(s)$ by analogy with the weight of a path, replacing entries $a_e$
with corresponding components of majorizers.

Let us provide a simple example to illustrate the idea.
Assume that $n=4$, and that $\net_{mjr}$ contains two vectors --- $(1,1,0,0)$ and $(1,0.7,0.5,0)$.
Assume that a typical realization of our matrix $A$ is

\begin{center}
\begin{tabular}{| c | c | c | c |}
\hline
  0 & 1 & -0.8 & 0 \\
\hline
  1 & 0 & 0.6 & 0.5 \\
\hline
  -0.8 & 0.6 & 0 & 0.3 \\
\hline
  0 & 0.5 & 0.3 & 0 \\
\hline
\end{tabular} 
\end{center}
Then we can assign to rows $1$ and $4$ the majorizer $(1,1,0,0)$, and to rows $2$ and $3$ --- majorizer $(1,0.7,0.5,0)$.
Consider a path $1\to 3 \to 2 \to 3 \to 1$ of length $4$.
Obviously, the absolute value of the path weight is $0.8*0.6*0.6*0.8$.
We can interpret this product as ``the second largest component of first row squared times
the second largest component of the second row squared''.
Now, replace the rows with corresponding majorizers.
We get ``the second largest component of majorizer $(1,1,0,0)$ squared times
the second largest component of majorizer $(1,0.7,0.5,0)$ squared'',
i.e.\ $1*1*0.7*0.7$. This is the weight of the data structure ${\bf Data}(1\to 3 \to 2 \to 3 \to 1)$.

Let us note that in the actual proof, it will be more convenient for us to define majorizers for vectors of
squares of the matrix entries, i.e.\ vectors of the form $(a_{ij}^2)_{j=1}^n$;
otherwise, our approach is very similar to the above example.
Since the number of majorizers is much smaller than the space of possible realizations of the matrix rows
, adding the information about the majorizers does not increase
the data space ${\bf S}$ by too much, and a satisfactory upper bound for $\sum\limits_{s\in{\bf S}}w(s)$
is possible. In fact, standard concentration inequalities imply that with very large probability a vast majority
of the matrix rows can be efficiently majorized by a single vector which we call a {\it standard majorizer}.
For those rows, no additional information should be added to the data structure which further controls the complexity of the space ${\bf S}$.

\medskip

{\bf Paths with many edges of multiplicity one.}
In classical applications of the trace method, paths having edges of multiplicity one do not participate
in the counting process since the expectation of the corresponding path weight is zero.
As we already mentioned above, in our setting taking the unconditional expectation of the trace cannot give
a satisfactory upper bound whereas conditioning on a ``good'' event (say, matrix realizations with predefined
statistics of norms of their rows and columns) produces complex dependencies within the matrix,
and computation of the conditional expectation becomes challenging.
For paths with relatively few
edges of multiplicity one, instead of the averaging,
we compute an upper bound for the sum of path weights valid everywhere inside the ``good'' event,
by bounding (the absolute value of) every path weight by the data structure weight, followed by some computations
which rely on the structure of the data space (see discussion above).
However, this approach is not applicable for paths having many edges of multiplicity one.
For a typical realization of our matrix, the path weights of such paths are split into approximately equal parts
according to their sign, and multiple cancellations occur.
Bounding each path weight individually by its absolute value destroys these cancellations and
cannot produce a satisfactory estimate.
The approach we take is to track these cancellations while conditioning on a ``good'' realization of the matrix,
even though this conditioning produces dependencies across the matrix.
Our method provides satisfactory estimates only for paths having relatively many multiplicity one edges,
and thus complements the argument based on the data structures.
The structure of the event we are conditioning on plays a crucial role.
The formal description of the method is given in Section~\ref{s: many mult1}.
Here, we would like to give a geometric viewpoint to it.
Let $A$ be our $n\times n$ symmetric random matrix, and let $G$ be the corresponding random graph on $[n]$
(whose edges mark non-zero entries of $A$).
Assume that we need to compute the conditional expectation of a product $\prod\limits_{e\sim \path}a_e$ given the event
$$\Event=\Big\{\sum\limits_{j=1}^n a_{ij}^2\leq T\mbox{ for all $i\in[n]$}\Big\},
$$
where $T>0$ is an appropriately chosen parameter. Geometrically, the event $\Event$
defines a bounded set $\mathcal S$ within the linear space of symmetric $n\times n$ matrices $M=(\mu_{ij})$,
and the conditional expectation of $\prod\limits_{e\sim \path}a_e$
can be viewed as integration of $\prod\limits_{e\sim \path}\mu_e$ over $M\in\mathcal S$ with respect to an appropriate measure.
It turns out that neighborhoods of points within $\mathcal S$ located far from the boundary $\partial \mathcal S$,
contain approximately equal mass of positive and negative realizations of $\prod\limits_{e\sim \path}\mu_e$,
and a cancellation within the integral can be verified. More precisely, we can show that if the number of surfaces
$Q_i:=\{M=(\mu_{ij}):\;\sum\limits_{j=1}^n \mu_{ij}^2= T\}$, $i\leq n$ which are close to our chosen point, is not very large
then necessarily a neighborhood of this point is ``well balanced'' in terms of the values of $\prod\limits_{e\sim \path}\mu_e$.
Thus, the parts of $\mathcal S$ in which cancellation does not happen are located close to the ``corners'' of $\mathcal S$,
and their measure is very small if the number of multiplicity one edges in $\path$ is large.
Therefore, the conditional expectation of $\prod\limits_{e\sim \path}a_e$ given $\Event$ is very close to zero.

\medskip

{\bf Non-centered matrices.} Let us discuss here how to get  an asymptotically sharp upper estimate for $|\lambda_{|k|}(W_n)|$ ($k\geq 2$)
for a sequence of non-centered matrices $(W_n)$, with the help of Theorem~A for sparse centered matrices.
As it was already mentioned in the introduction, standard centralization procedures ---
replacing $W_n$ with $W_n-\E W_n$ or with the difference of two independent copies of $W_n$ ---
do not allow to reduce the problem of bounding the eigenvalues of a non-centered sparse matrix to the setting of Theorem~A.
Indeed, in the first case we would obtain a non-sparse random matrix whereas the second approach indroduces
a suboptimal constant factor in the final estimate.
To deal with this issue, we carry out a different symmetrization procedure designed specially for sparse matrices.
To avoid technical details in this overview, we describe a model different at certain points from the one we actually use
but easier to discuss informally.

Assume that $W_n$ is a sparse Bernoulli matrix (i.e.\ the adjacency matrix of the Erd\H os--Renyi graph;
we can assume for simplicity that loops are allowed)
and let $p_n$ be the success probability of each entry.
Fix a small constant $\varepsilon>0$ and consider another Bernoulli random matrix $B_n$ with the entries having
probability of success $p_n/\varepsilon$. Then the difference
$W_n- \varepsilon B_n$ is a centered (and still sparse) random matrix, and we could apply Theorem~A
to get an upper bound on its largest eigenvalue. Note that the sparsity parameter
of the matrix $W_n-\varepsilon B_n$ is $1-(1-p_n)(1-p_n/\varepsilon)\approx p_n/\varepsilon$,
whereas the variance of each entry is $p_n + \varepsilon p_n\approx p_n$. Therefore, to apply Theorem~A
we should rescale the matrix by the factor $\varepsilon^{-1/2}$ (approximately) .
We get that for large enough $n$, with large probability
\begin{align*}
\varepsilon^{-1/2}\,\|W_n-\varepsilon B_n\|&\lesssim
\sqrt{\max\big(\max\nolimits_{i}\|\varepsilon^{-1/2}\,\Row_i(W_n-\varepsilon B_n)\|_2^2-np_n/\varepsilon,np_n/\varepsilon\big)}\\
&\hspace{2cm}+ \frac{np_n/\varepsilon}{\sqrt{\max\big(\max_{i}\|\varepsilon^{-1/2}\,\Row_i(W_n-\varepsilon B_n)\|_2^2-np_n/\varepsilon,np_n/\varepsilon\big)}}.
\end{align*}
Our next observation is that the maximal norm of the rows of $\varepsilon B_n$ is typically much smaller
than the maximal norm of the rows of $W_n$. As an informal justification, we can point to the trivial relation
$\E\|\Row_i(\varepsilon B_n)\|_2^2=\varepsilon n p_n\ll n p_n=\E\|\Row_i(W_n)\|_2^2$, $i=1,2,\dots,n$.
This allows us to replace $\|\Row_i(W_n-\varepsilon B_n)\|_2^2$ with $\|\Row_i(W_n)\|_2^2$
in the last display formula, so that, after cancelling $\varepsilon$, we obtain
\begin{align*}
\|W_n-\varepsilon B_n\|&\lesssim
\sqrt{\max\big(\max\nolimits_{i}\|\Row_i(W_n)\|_2^2-np_n,np_n\big)}\\
&\hspace{3cm}+ \frac{np_n}{\sqrt{\max\big(\max_{i}\|\Row_i(W_n)\|_2^2-np_n,np_n\big)}}.
\end{align*}
As the last step, we observe that $\|W_n-\varepsilon B_n\|=\|W_n-\E W_n+\varepsilon\E B_n- \varepsilon B_n\|$,
where $\varepsilon\E B_n- \varepsilon B_n$ is a {\it relatively small} perturbation of the matrix $W_n-\E W_n$,
and cannot significantly decrease the spectral norm.
The way it is actually done in our proof is to consider a random unit vector $X$ measurable with respect to $W_n$
and such that $\|(W_n-\E W_n)X\|_2=\|W_n-\E W_n\|$ everywhere on the probability space.
Then, using our assumption that $\varepsilon$ is small and that $X$ and $B_n$ are independent,
it is possible to show that $\|(B_n-\E B_n)X\|_2$ is small with probability tending to one, implying that
$$\|W_n-\varepsilon B_n\|\gtrsim \|W_n-\E W_n\|.$$
Combining this relation with the previous formula, and letting $\varepsilon\to 0$, we obtain an asymptotically
sharp upper bound for $\lambda_{|2|}(W_n)\leq\|W_n-\E W_n\|$.

\medskip

{\bf Lower bound for the $k$--th largest eigenvalue.}
As the last element of the proof, we discuss the lower bound for the largest eigenvalues.
Unlike the upper bound, this part of the argument does not use the trace representation,
and instead is based on a mixture of combinatorial arguments (which give us necessary structural
information on the underlying graph) and ``geometric'' methods which use an explicit construction of a random
vector capturing the values of the leading eigenvalues.
For the sake of simplicity, we will only discuss the bound for the operator norm; an estimate of $\lambda_{|k|}(W_n)$
for $k\geq 2$ is obtained by decomposing the matrix $W_n$ into $k$ blocks and carrying out the argument sketched
below on each of the blocks.

We recall that the basic test case in the study of the {\it upper bound} is a tree of finite (but large) depth
rooted at the vertex corresponding to the row/column with the largest Euclidean norm.
In a sense, the whole proof of the upper bound can be viewed as a justification of the fact that this
test case presents a significant contribution into the sum-of-paths representation of the trace.
For the lower bound, we take the same tree in the underlying random graph but this time we construct a special
random vector $Y$ (modelled in accordance with the tree structure) such that $\|W_n Y\|_2/\|Y\|_2$ is close
to the spectral norm of $W_n$.

Let $W_n=(w_{ij})$ be a matrix with the entries taking values in $\{0,\pm 1\}$, and let $G$ be the graph on $[n]$
with the edge set corresponding to non-zero entries of $W_n$. We will further assume that the vertex of $G$ having the largest
degree is $1$, and that the $q$--neighborhood of $1$ in $G$ (for a large constant $q$) is a tree whose nodes,
except for the root and the leaves, have a degree $d$.
Our assumption that the neighborhood is a tree, is reasonable
in view of the sparsity of our model (we avoid going into technical details here). 
For any integer $r\in[0, q]$, we let $V_r$ be the set of all vertices of the tree
having depth $r$ (so that, in particular, $V_0=\{1\}$). 
We then define our vector $Y$ as
$$Y:=\sum\limits_{r\in \Z_2\cap[0,q]}\sum\limits_{u\in V_r}Y_u,$$
where for each $u\in V_r$, we set
$$Y_u:=\delta_r\, \sum\limits_{\tiny\mbox{$z$: $z$ is a child of $u$}} w_{uz} e_z;$$
$\delta_r>0$, $r\in\Z_2\cap[0,q]$, are some parameters
and $(e_z)_{z=1}^n$ is the canonical basis in $\R^n$.
The vectors $Y_u$ have disjoint supports, and thus $Y$ can be viewed as a weighted combination
of tree nodes of odd depth; with the weight $\delta_r$ shared by all nodes (basis vectors)
of depth $r+1$. The latter condition is reasonable as the nodes of the tree having the same depth are
indistinguishable when unlabelled. The condition that we take only even values of $r$ is not crucial;
at the same time any given matrix row is supported on vertices which either all have odd depth or all have even depth, and
this produces a recursive relation between $\delta_r$ and $\delta_{r-2}$ in the computations
making the separation of odd and even tree layers somewhat natural.

The formula for the vector $Y$ was obtained by trial and error although
the above remarks suggest that the choice of such structure is quite natural.
The values of the parameters $\delta_r$ which produce maximal (or close to maximal)
value of $\|W_n Y\|_2/\|Y\|_2$ can be easily deduced; in fact, we take $(\delta_r)$ as a geometric
series determined by $d$ and the degree of vertex $1$ (we refer to the proof of Lemma~\ref{l: main lower} for details).

\section{Preliminaries}

\vskip 0.2cm


Given two integers $a\leq b$, we denote by $[a,b]$ the corresponding integer interval.
For the interval on the real line with the same boundary points, we will use notation $[a,b]_\R$.

Let $G=(V,E)$ be an undirected (simple) graph.
Given vertices $i$ and $j$ in $V$, the edge connecting $i$ and $j$ is denoted by $i\edg j$.
A {\it path} $\path$ of length $s$ on $G$ is a sequence of $s+1$ vertices of $G$
where each pair of successive vertices forms an edge in $G$.
We allow the vertices in the path to repeat.
A path is {\it closed} if the first and last vertex of the sequence coincide.
It will be convenient for us to view a path on the graph as a mapping from an integer interval to the set of vertices $V$.
If $\path:[0,a]\to V$ is a path on $G$ then $\path(0),\path(1),\dots,\path(a)$ are vertices travelled by $\path$.
For any $0\leq b\leq a$, by $\path[0,b]$ we denote the subpath $\path(0),\path(1),\dots,\path(b)$.

For any connected graph $G=(V,E)$ there is a natural metric $d(\cdot,\cdot)$ on its set of vertices $V$
induced by the graph distance i.e.\ for every $u,v\in V$, $d(u,v)$ is the length of the shortest path in
$G$ starting at $u$ and ending in $v$.
For any $r\geq 0$, the $r$-neighborhood of a vertex $v\in V$
is the subgraph of $G$ obtained by removing vertices with distance more than $r$ to $v$.
The {\it diameter} of $G$ is the largest distance between any pair of its vertices.
Further, a subset $S\subseteq V$ is {\it $r$--separated} if $d(u,v)>r$ for any two distinct elements of $S$.
We say that $S$ is a maximal $r$--separated subset of $V$ if it is $r$--separated and there is no vertex
$v\in V$ whose distance to $S$ is strictly greater than $r$.
Everywhere in the text below,
the terms ``distance'', ``diameter'',  ``$r$-neighborhood'' and ``$r$-separated set'' are used in the above sense. 

Given an edge $e\in E$ of $G=(V,E)$, we denote by $\inc(e)$ the set of two vertices incident to $e$
i.e.\ if $e=i\edg j$ then $\inc(e)=\{i,j\}$.
Further, for any $S\subseteq E$ we define $\inc(S):=\bigcup_{e\in S} \inc(e)$.
Conversely, given a subset $V'\subseteq V$, by $\inc(V')$ we denote the collection of all edges of $G$
which are incident to some vertex in $V'$. 
We say that $e\in E$ is a {\it cycle edge} of $G$
if it belongs to a cycle within $G$; otherwise we say that $e$ is a {\it non-cycle edge}.
We say that a graph $G=(V,E)$ is {\it $\ell$-tangle free} if every neighborhood of radius $\ell$ in $G$ has at most one cycle.

The next lemma provides a standard upper bound on the number of $r$-separated points in a connected graph. 
\begin{lemma}\label{lem: card-net}
Let $G=(V,E)$ be a connected graph containing $s$ vertices which are $r$-separated. Then $s\leq 2\vert E\vert/r$. 
\end{lemma}

It is a standard fact that any connected graph admits a subgraph with the same vertex set and without cycles
(a spanning tree). We will need a relative of that property which tells us that we can delete half of cycle edges
from a certain prescribed collection without destroying connectivity.
\begin{lemma}\label{lem: edg-remov-connect}
Let $G$ be a connected graph and $\w G=(\w V,\w E)$ be a connected subgraph of $G$. Further, let $S$ be a collection of cycle edges of $G$ not contained in $\w E$ but incident to $\w V$. 
Then there exists $S'\subseteq S$ of cardinality at least $\vert S\vert/2$
such that the graph obtained from $G$ by removing edges from $S'$ is still connected.
\end{lemma}
\begin{proof}
Let $S'$ be a subset of $S$ of maximal cardinality such that removing the edges in $S'$ keeps the graph connected, and assume that $\vert S'\vert <\vert S\vert/2$. 
For every edge $e\in S\setminus S'$, let $v(e)$ be the endpoint of $e$ not contained in $\w V$ (note that such vertex exists by the connectivity of $\w G$ and the maximal property of $S'$). 
Take any edge $e\in S\setminus S'$. We will show that there exists a path $\path_e$ on $G$ starting at $v(e)$, ending at a vertex in $\inc(S')\setminus \w V$  and not passing through edges in $\w E\cup S$.

Indeed, if $v(e)\in \inc(S')\setminus \w V$ then there is nothing to show.
Otherwise, since $e$ is a cycle edge, there is a cycle $\cycle_e$ in $G$ containing $e$. Let $\path_e'$ be the trail on $\cycle_e$ starting at $v(e)$, ending at $\inc(e)\setminus v(e)$, and not passing through $e$. Denote by $t_e\geq 1$ the first time when $\path_e'(t_e)\in \w V\cup \inc(S')$ (such time always exists since the path $\path_e'$ ends in $\w G$).
Observe that if $\path_e'(t_e)$ belonged to $\w V$ then we would obtain a path on $G$ from $v(e)$ to $\w V$ not passing through $S'$, which contradicts the maximality of $S'$ (since the removal of $e$ would keep the graph connected). Hence, $\path_e'(t_e)\in \inc(S')\setminus \w
V$, and $\path_e'([0,t_e])$ cannot contain any edges from $S$ since otherwise one of its vertices would belong to $\w V$.
It remains to set $\path_e:=\path_e'[0,t_e]$ to obtain a path starting at $v(e)$, ending at a vertex in $\inc(S')\setminus \w V$  and not passing through edges in $\w E\cup S$. 

Since $\vert S'\vert <\vert S\vert/2$, we have
$\vert S\setminus S'\vert > \vert \inc(S')\setminus \w V\vert$ and by the pigeonhole principle, there are $2$ edges $e,e'\in S\setminus S'$ such that 
the corresponding paths; $\path_e$ and $\path_{e'}$, end at the same vertex $v\in \inc(S')\setminus \w V$. Concatenating these two paths,
get that $v(e)$ and $v(e')$ are connected through a path on the edges in $(\w E\cup S)^c$. Therefore, one of the edges $e$ or $e'$ could be removed without destroying connectivity of the graph.
Since both of these edges do not belong to $S'$, this contradicts the maximality of $S'$
and we deduce that $\vert S'\vert \geq \vert S\vert/2$. 
\end{proof}

As a consequence of the last lemma, we can bound from above the number of cycle edges of an $\ell$-tangle free graph $G$
which are incident to a connected subgraph of small size.

\begin{lemma}\label{lem: cycle-edge-neighbor}
Let $G=(V,E)$ be a connected $\ell$-tangle free graph, $\w G=(\w V, \w E)$ be a connected subgraph of $G$ with $\diam(\w G)\leq (\ell-2)/3$, and $S$ be the collection of all cycle edges of $G$ incident to $\w V$ but not contained in $\w E$. Then $\vert S\vert \leq 6+\frac{16\vert E\vert}{\ell}$.
\end{lemma}
\begin{proof}
Let $S'\subseteq S$ be the set obtained from Lemma~\ref{lem: edg-remov-connect} and let $G'$ be the (connected) subgraph of $G$
with the same vertex set $V$ and with the edge set $E\setminus S'$. 
Let $\net$ be a maximal $\ell/2$--separated net in $G'$ and note that, by Lemma~\ref{lem: card-net},
we have $\vert \net\vert \leq 4\vert E\vert/\ell$.  
For each vertex $v\in \inc(S')\setminus \w V$, let $u_v\in \net$ be such that $d(v,u_v)\leq \ell/2$. 

Assume first that $\vert \inc(S')\setminus \w V\vert \geq \vert \net\vert +2$. Then, by the pigeonhole principle, there are at least $2$ couples of vertices $\{i_1,j_1\}, \{i_2,j_2\}\subseteq \inc(S')\setminus \w V$ with $i_1\neq j_1, i_2\neq j_2$, $u_{i_1}=u_{j_1}$ and $u_{i_2}=u_{j_2}$. This implies that $i_1$ and $j_1$ can be connected by a path in $G'$ of length at most $\ell$, and similarly for $i_2$ and $j_2$. At the same time, $i_1$ and $j_1$ can be connected by a path on $S'\cup \w E$, and similarly for $i_2$ and $j_2$. Therefore, there are at least two distinct cycles in $G$ of length at most $\ell+2+\diam(\w G)$, and with distance between the two cycles at most $\diam(\w G)$. Hence, there is a vertex $v$ in $G$ such that 
its $\big(\frac{\ell+2+\diam(\w G)}{2}+\diam(\w G)\big)$-neighborhood contains at least $2$ cycles, which contradicts
the $\ell$-tangle free property of $G$ by our restriction on $\diam(\w G)$. 

Thus, we deduce that $\vert \inc(S')\setminus \w V\vert\leq \vert \net \vert +1$. Further, for any $v\in \inc(S')\setminus \w V$, there are at most two edges in $S'$ incident to it, as otherwise we would get at least two cycles in $G$ of length at most $2+\diam(\w G)$ in the neighborhood of $v$ contradicting the $\ell$-tangle free property of $G$. Finally, observe that there is at most one edge in $S'$ with both end points in $\w V$, as otherwise we would get $2$ cycles of lengths at most  $1+\diam(\w G)$ at distance at most $\diam(\w G)$ from one another, leading to a contradiction as before. 

Hence, we get $\vert S'\vert\leq 1+2\vert \inc(S')\setminus \w V\vert \leq 3+2\vert \net\vert$. Combining this with the bound $\vert S\vert\leq 2\vert S'\vert$ from Lemma~\ref{lem: edg-remov-connect} and the one on the cardinality of $\net$, we complete the proof. 
\end{proof}

\vskip 0.2cm

\subsection{Edge and vertex discovery, and statistics of special vertex types}\label{subs: stat of cycle vertices} \ \vskip 0.2cm

Let $G=(V,E)$ be an undirected graph. 
As we mentioned above, it will be convenient to see a path on $G$ of length $\kappa$ as a function $\path: [0,\kappa]\to V$
where for every $t\in[0,\kappa]$, $\path(t)$ indicates the vertex reached by the path at the step (``time'') $t$.
For each path $\path:[0,\kappa]\to V$ on $G$, we denote by $G_\path= (V_\path, E_\path)$
the undirected subgraph of $G$ obtained by deleting the edges and vertices not visited by $\path$.  

Let $\cycleedges$ be the set of times $t\in [1,\kappa]$ such that $\path(t-1)\edg \path(t)$ is a cycle edge of $G_\path$ travelled either once or at least three times by $\path$. 
Now, for any  cycle edge $e$ in $G_\path$, {\it the discovery time $\disc(e)$} of $e$ is the smallest $t\in[1,\kappa]$ 
such that 
$$
\text{ Either }\quad t\in \cycleedges\text{ and }e=\path(t-1)\edg \path(t) \quad \text{ or } \quad G_{\path[0,t]} \text{ has a cycle containing $e$.}
$$
Note that the {\it discovery time for a cycle edge of the graph may be different from the first time the edge is traveled by the path.}

\medskip

In this subsection, we discuss three types of special vertices of the graph $G_\path=(V_\path,E_\path)$:
cycle meeting points, splitting points and cycle completion points.
The main purpose is to show that under the assumption that $G_\path$ is $\ell$-tangle free
(for a large enough $\ell$), the total 
number of the special vertices is very small.

Let $\path'$ be some path on $G=(V,E)$ (we change notation here as $\path'$ may serve both as the path $\path$ or as
its sub-path).
For every $e=i\edg j\in E_{\path'}$, we define the {\it direction of uncovering}
of $e$ in $\path'$ as the direction in which it first appeared in the $\path'$ i.e. 
$$
\uncov_{\path'}(e)=\left\{
	\begin{array}{ll}
		i\to j  & \text{ if there is $t$ such that $\path'(t-1) =i$, $\path'(t)=j$ and $e\not\in E_{\path'[t-1]}$} \\
		j\to i & \mbox{otherwise}.
	\end{array}
\right.
$$
The {\it time of uncovering} for $e$ is the smallest $t$ such that $e$ is an edge of $G_{\path'[t]}$.
We once again would like to turn the Reader's attention to the distinction between times of discovery and uncovering
of a cycle edge in the graph.

Let $v$ be a vertex of $G_{\path'}$.
\begin{itemize}

\item $v$ is {\it cycle meeting point} of $G_{\path'}$
if there are at least three distinct vertices $i,j,u$ of $G_{\path'}$ such that $v\edg i$, $v\edg j$ and $v\edg u$ are cycle edges of $G_{\path'}$.

\item $v$ is a {\it splitting point} of $G_{\path'}$
if there exist vertices $i\neq j$ such that
$v\edg i$, $v\edg j$ are cycle edges in $G_{\path'}$ and 
$\uncov_{\path'}(v\edg i)=v\to i$ and $\uncov_{\path'}(v \edg j)= v\to j$.

\item $v$ is a {\it cycle completion point} if there is time $t$ such that
$\path'(t)=v$ and the number of cycles in $G_{\path'[t-1]}$ is strictly less than the number of cycles in $G_{\path'[t]}$.
\end{itemize}

All vertices of $G_{\path'}$ which belong to one of the above types, will go under the name {\it special cycle vertices} of $G_{\path'}$.
Note that a given vertex $v$ can be simultaneously of more than one of the above types.
Next, we define the concept of {\it discovery times} for the special cycle vertices.
{\it{}The discovery time $\disc(v)$} of a special cycle vertex $v$ of $G_{\path'}$
is the smallest $t$ such that $v$ is a special cycle vertex of $G_{\path'[t]}$.
Notice that the discovery time of a special cycle vertex {\it is not} equal to the first time the path visits this vertex.

\bigskip

As a simple corollary of Lemma~\ref{lem: cycle-edge-neighbor},
we get the following statement which bounds the number of cycle edges incident to a cycle meeting point.
\begin{cor}\label{cor: nb-cycles meeting}
Assume that $G_\path$ is $\ell$-tangle free with $\ell\geq 5$,
and let $v$ be a cycle meeting point in $G_\path$.
Then the total number of cycle edges of $G_\path$ incident to $v$ does not exceed $6+\frac{16\vert E_\path\vert}{\ell}$. 
\end{cor}

In the next lemma, we bound the total number of the cycle meeting points. 

\begin{lemma}\label{lem: nb-cycle meeting points}
Assume that the graph $G_\path$ is $\ell$-tangle free, with $\ell \geq 2$.
Then the total number of cycle meeting points in $G_\path$ is at most $C\vert E_\path\vert^2/\ell^2$ for some universal constant $C$. 
\end{lemma}
\begin{proof}
Let $N$ be the number of cycle meeting points in $G_\path$.
We will assume that $N\geq 32 \vert E_\path\vert/\ell$ (otherwise, there is nothing to prove).
Denote by $\net$ a maximal $\ell/8$-separated net in $G_\path$
and note that by Lemma~\ref{lem: card-net},
we have $\vert \net\vert \leq 16\vert E_\path\vert/\ell$. Clearly, there exists a vertex $v\in \net$ such that its $\ell/8$-neighborhood
$B$ contains at least $\frac{\ell}{16\vert E_\path\vert} N\geq 2$ cycle meeting points. 

Let $\tree$ be a spanning tree for $B$, and denote by $\tree'$ a subtree of $\tree$ obtained by
successively removing edges incident to degree one vertices which are not cycle meeting points of $G_\path$
(and then throwing away the obtained isolated vertices).
Note that this operation will necessarily keep all cycle meeting points which fall into $B$, inside $\tree'$
as any two such vertices are connected by a path within $B$, whose edges cannot be removed.
In particular, $\tree'$ is non-empty.
Note also that all leaves of $\tree'$ are by construction the cycle meeting points of $G_\path$. 

Let $S$ be the set of all cycle edges of $G_\path$ which are incident to vertices in
$\tree'$ but not contained in $\tree'$. Note that, by the definition of a cycle meeting point,
for any leaf of $\tree'$, there are at least $2$ edges from $S$ incident to it.
Further, for any cycle meeting point which is a node of degree $2$ in $\tree'$, there is at least one edge from $S$ incident to it. 
Finally, it is easy to check that the number of nodes of degree at least $3$ in any tree is less than the number of its leaves. 
Since any edge in $S$ is incident to at most $2$ cycle meeting points in $\tree'$, the above observations imply that 
\begin{align*}
\frac{\ell}{16\vert E_\path\vert} N &\leq \vert\{\mbox{cycle meeting points of $G_\path$ which are leaves of $\tree'$}\}\vert\\
&\hspace{0.5cm}+\vert\{\mbox{cycle meeting points of $G_\path$ which are nodes of degree $2$ in $\tree'$}\}\vert\\
&\hspace{0.5cm}+\vert\{\mbox{cycle meeting points of $G_\path$ which are nodes of degree at least $3$ in $\tree'$}\}\vert\\
&\leq 2\vert\{\mbox{cycle meeting points of $G_\path$ which are leaves of $\tree'$}\}\vert\\
&\hspace{0.5cm}+\vert\{\mbox{cycle meeting points of $G_\path$ which are nodes of degree $2$ in $\tree'$}\}\vert\\
&\leq 2\vert S\vert+2\vert S\vert=4\vert S\vert. 
\end{align*}
Applying Lemma~\ref{lem: cycle-edge-neighbor} with $\w G= \tree'$
(note that $\diam(\tree')\leq \lfloor \ell/4\rfloor\leq (\ell-2)/3$), we get that 
$$
\frac{\ell}{16\vert E_\path\vert} N\leq 4\vert S\vert\leq 4\Big(6+\frac{16\vert E_\path\vert}{\ell}\Big), 
$$
whence $N\leq C\vert E_\path\vert^2/\ell^2$ for some appropriate constant $C$. 
\end{proof}

In order to bound the number of splitting and cycle completion points, it will be convenient to introduce
a notion of a {\it cycle interval}.
By a {\it cycle interval} in $G_\path$
we understand an (unordered)
collection of distinct cycle edges in $G_\path$ of the form $i_1\edg i_2, i_2\edg i_3,\ldots, i_{q-1}\edg i_q$, where 
none of the vertices $i_2,i_3,\dots,i_{q-1}$ are cycle meeting points of $G_\path$ and,
if $i_1\neq i_q$, then $i_1$ and $i_q$ are cycle meeting points of $G_\path$. 
Note that with this definition, every cycle edge of $G_\path$ belongs to a unique cycle interval,
and no two distinct cycle intervals share more than $2$ common vertices.
Furthermore, any cycle interval $I$ belongs to one of the following two types:
\begin{itemize}
\item Either $I$ is a full cycle of $G_\path$ which has at most one common vertex with other cycles of $G_\path$,
\item Or $I$ is a collection of consecutive edges of some cycle in $G_\path$ which
is bounded from both sides by two distinct cycle meeting points.
\end{itemize}
In the next lemma, we give an upper bound on the number of cycle intervals. 

\begin{lemma}\label{lem: cycle-interval}
Assume that the graph $G_\path$ is $\ell$-tangle free, with $\ell\geq 5$.
Then the number of cycle intervals in $G_\path$ is at most $C\vert E_\path\vert^3/\ell^3$ for some universal constant $C>0$.
\end{lemma}
\begin{proof}
The definition of a cycle interval implies that each interval $I$ in $G_\path$
is either a full cycle of $G_\path$ or has cycle meeting points as its boundary vertices.
In view of Lemma~\ref{lem: nb-cycle meeting points} and Corollary~\ref{cor: nb-cycles meeting},
the number of intervals of the second type is bounded above by $C' \vert E_\path\vert^3/\ell^3$ for some constant $C$. 
To count the number of intervals in $G_\path$ which are full cycles,
note that since the graph is $\ell$-tangle free, the number of cycles of $G_\path$ of length at most $\ell$ is at most $4\vert E_\path\vert/\ell$.
Indeed, this can be verified by constructing an $\ell/2$-separated set in $G_\path$ consisting of ``representative''
vertices which belong to distincs short cycles, and then applying Lemma~\ref{lem: card-net}.
Finally, it remains to note that the number of disjoint cycles of length greater than $\ell$ is at most $\vert E_\path\vert/\ell$. 
\end{proof}

\begin{lemma}\label{l: number of splitting}
Let $I=\{i_1\edg i_2,i_2\edg i_3,\dots,i_{q-1}\edg i_q\}$ be a cycle interval in $G_\path$. 
Then,
among the vertices $\{i_2,i_3,\dots,i_{q-1}\}$, there is at most one splitting point.

Moreover, if $G_\path$ is $\ell$-tangle free, then the total number of splitting points in $G_{\path}$
is at most $C\vert E_\path\vert^3/\ell^3$, for some universal constant $C$.  
\end{lemma}
\begin{proof}
Assume that the cycle interval $I$ contains $2$ distinct splitting points in its interior,
say, $i_u, i_{u'}$ for some $u,u'\in [2,q-1]$. Assume further that the point $i_u$ is visited for the first time by $\path$ earlier
than $i_{u'}$. The vertex $i_{u'}$ cannot be reached for the first time by $\path$ through any of the edges $i_{u'-1}\edg i_{u'}$
or $i_{u'+1}\edg i_{u'}$ since this would violate their uncovering directions.
But then $i_{u'}$ must be accessed through another edge which necessarily would become a cycle edge of $G_\path$.
But then $i_{u'}$ is a cycle meeting point -- contradiction. Thus, each cycle interval contains at most one splitting
point in its interior.

To prove the second part of the lemma,
note that the number of splitting points is at most three times the number of cycle intervals in $G_\path$
(if we count $i_1$, $i_q$ as possible splitting points in $I$),
and Lemma~\ref{lem: cycle-interval} implies the estimate. 
\end{proof}

\begin{lemma}\label{l: n of cycle compl}
Let $G_\path$ be $\ell$-tangle free.
Then the total number of cycle completion points of $G_\path$ is bounded above by
$C\vert E_\path\vert^3/\ell^3$ for some universal constant $C>0$.
\end{lemma}
\begin{proof}
Note that each cycle interval may contain at most $1$ cycle completion point in its interior.
The result follows by applying Lemma~\ref{lem: cycle-interval}.
\end{proof}
\medskip

\vskip 0.2cm
\subsection{Majorizers}\label{subs: majorizer1}\ \vskip 0.2cm

Given a vector $x$ in $\R^n$, let $x^*$ be the non-increasing rearrangement of the
absolute values of coordinates of $x$. Further, given a non-increasing vector $y\in\R^n_+$,
we say that $y$ is {\it a majorizer} for $x$, and write $y\geq x^*$, if $y_i\geq x^*_i$ for all $i\leq n$.

We will consider a collection of majorizers for a special class of vectors. 
Given $h\geq 2$ and $\gamma,s\in \R_+$, 
define the set $\rearr(h,\gamma,s)$ of all $n$--dimensional vectors $x=(x_1,x_2,\dots,x_n)$,
such that $h\geq x_1\geq x_2\geq\dots\geq x_n\geq 0$, $\|x\|_1\leq \gamma$, and $x$ has at most $s$ non-zero coordinates.
We have the following lemma:
\begin{lemma}\label{l: net of majorizers}
For $h\geq 2$, $\gamma,s\in \R_+$ and $ \varepsilon\in(0,1/2]$ satisfying $\frac{\varepsilon\gamma}{h}\geq
C_{\text{\tiny \ref{l: net of majorizers}}}$
and $\frac{hs}{\varepsilon \gamma}\geq 2$,
there exists a subset $\net=\net(h,\gamma, s,\varepsilon)\subset\R^n_+$ of cardinality at most
$\Big(\frac{C_{\text{\tiny \ref{l: net of majorizers}}}\log_2
(\frac{hs}{\varepsilon \gamma})}{\varepsilon}\Big)^{C_{\text{\tiny \ref{l: net of majorizers}}} \varepsilon^{-2}\log_2\frac{h}{ \varepsilon}}$
such that for every $x\in \rearr(h, \gamma,s)$ there exists $y\in\net$ which is a majorizer for $x$,
and, moreover, $\|y\|_1\leq (1+\varepsilon)\gamma$.
Here, $C_{\text{\tiny \ref{l: net of majorizers}}}>0$ is a universal constant.
\end{lemma}
\begin{proof}
For any vector $x\in\rearr(h,\gamma,s)$, we write 
$$
x=\bar{x}+\w x,
$$
where $\bar{x}$ has at most $s$ non-zero coordinates all of which are smaller than $\frac{\varepsilon \gamma}{s}$ while $\w x$ has all its non-zero coordinates lying between 
$\frac{\varepsilon \gamma}{s}$ and $h$. Note that the rearranged vector $\bar{y}$ having $s$ non-zero coordinates all of which are equal to $\frac{\varepsilon \gamma}{s}$ is a majorizer for $\bar{x}$ for any $x\in \rearr(h,\gamma,s)$. Moreover $\Vert \bar{y}\Vert_1\leq \varepsilon\gamma$. Therefore, our task is to construct a majorizer $\w y$ for any given $\w x$ such that $\Vert \w y\Vert_1\leq (1+\varepsilon)\gamma$, then take $y=\w y\oplus \bar{y}$ which is a majorizer of $x$ satisfying $\Vert y\Vert_1\leq (1+2\varepsilon)\gamma$. 

Having in mind
the above reduction procedure,
we may therefore consider the set $\w \rearr(h,\gamma,s)$ of all $n$--dimensional vectors $x=(x_1,x_2,\dots,x_n)$,
such that $h\geq x_1\geq x_2\geq\dots\geq x_n\geq 0$,
$\|x\|_1\leq \gamma$, $x$ has at most $s$ non-zero coordinates, and each non-zero coordinate is at least $\frac{\varepsilon \gamma}{s}$.
We now construct a majorizer for any given vector $x\in\w\rearr(h,\gamma,s)$ and verify that it satisfies the required properties.
Define $p_0:=\lfloor\frac{\varepsilon \gamma}{4h}\rfloor$,
$$
p_i:=\bigg\lfloor 2^{c_1\varepsilon^{\,2} i}\,\frac{c_1\varepsilon^{\,2} \gamma}{h}\bigg\rfloor,\quad i=1,2,\dots;
\quad r_i:=\sum\limits_{j=0}^i p_j,\quad i=0,1,2,\dots,
$$
where $c_1>0$ is a universal constant which will be determined later.
Further, for every $i\geq 1$ we let
$$
w_i:=\inf\big\{h 2^{-m\varepsilon/4}:\;m\geq 0,\;h 2^{-m\varepsilon/4}\geq x_{1+r_{i-1}}\big\},
$$
and, finally,
define the $n$--dimensional vector $y=(y_1,\dots,y_n)$
by
$$
y_j:= h\mbox{ for $j\leq r_0$,}\;\;\mbox{ and }\;\;
y_j:=w_i,\quad 1+r_{i-1}\leq j\leq r_{i},\quad i=1,2,\dots
$$
It is not difficult to see that $y\geq x$, just by our construction.
Next, we estimate the $\|\cdot\|_1$--norm of $y$. Observe that the $n$--dimensional vector
$x':=(x_{\lceil i/(1+\varepsilon/4)\rceil})_{i=1}^n$ satisfies $\|x'\|_1\leq (1+\varepsilon/4)\|x\|_1$
(this can be easily checked, say, by embedding $\ell_1^n$ into $L_1[0,n]$).
Further, by our construction, assuming that the constant $c_1$ is sufficiently small (and $C$ is sufficiently large),
we have
$$p_i\leq \frac{\varepsilon}{4}r_{i-1},\quad i=1,2,\dots,$$
implying $\lceil r_i/(1+\varepsilon/4)\rceil\leq 1+r_{i-1}$, so that
$$
x'_{j}\geq x_{1+r_{i-1}}\quad\mbox{ for all admissible }1+r_{i-1}\leq j\leq r_{i},\quad i=1,2,\dots
$$
Hence, the vector $2^{\varepsilon/4}x'+h\sum\limits_{i\leq \varepsilon \gamma/4h}e_i$
is a majorizer for $y$. This gives
$$\|y\|_1\leq 2^{\varepsilon/4}(1+\varepsilon/4)\gamma
+\varepsilon \gamma/4\leq (1+\varepsilon)\gamma.$$

Let $\net$ be the set of all vectors $v\in \R^n$ satisfying the following conditions:
\begin{itemize}
\item $h\geq v_1\geq v_2\geq \ldots\geq v_n\geq 0$;
\item $\Vert v\Vert_1\leq (1+\varepsilon)\gamma$;
\item $\{v_i\}_{i=1}^n\subset \big\{h 2^{-m\varepsilon/4}:\;m\geq 0,\;h 2^{-m\varepsilon/4}\geq \frac{\varepsilon \gamma}{s}\big\}\cup\{0\}$;
\item $v$ is constant on $[p_0]$ and on $[r_{i-1}+1,r_i]\cap[n]$, $i=1,2,\dots$.
\end{itemize}
It is immediate that the vector $y$ constructed above belongs to $\net$.
Thus, to finish the proof, it remains to estimate the cardinality of $\net$.
Observe that coordinates of any vector $v$ from $\net$ may take at most
$\frac{4\log_2 (hs/\varepsilon \gamma)}{\varepsilon}+2$ different values, and
for any $i\geq 1$ such that $r_{i-1}\geq 2h\gamma$, we necessarily have that $v$ is zero on $[r_{i-1}+1,r_i]\cap[n]$
(as otherwise the assumption on its $\|\cdot\|_1$--norm will be violated).
Thus, there are at most $C' \varepsilon^{-2}\log_2\frac{h}{ \varepsilon}$
``non-trivial'' levels of $v$, whence
$$
|\net|\leq \bigg(\frac{4\log_2 \big(\frac{hs}{\varepsilon \gamma}\big)}{\varepsilon}+2\bigg)^{C' \varepsilon^{-2}\log_2\frac{h}{ \varepsilon}},
$$
for an appropriate constant $C'>0$.
The result follows.
\end{proof}

\section{Mapping to a data structure}\label{s: data-structure}

Given an $n\times n$ symmetric matrix $M=(\mu_{ij})_{1\leq i,j\leq n}$ with zero diagonal, we  denote by $G_M=([n],E_M)$ the graph
with the edge set $E_M:=\{i\edg j:\,\mu_{ij}\neq 0\}$. With some abuse of terminology, we will say that
a vertex $v$ of $G_M$ is {\it majorized} by a vector $y\in\R^n_+$ if
the non-increasing rearrangement of
the sequence $(\mu_{vi}^2)_{i=1}^n$ is majorized (coordinate-wise) by $y$. Given a vector $\stnd\in \R_+^n$, we say that a vertex $v$ of the graph $G_M$ is {\it $\stnd$--heavy}
if it is not majorized by $\stnd$.

Let $n$ be a large natural number, $\ell, d,\dmax\in \N$ and $h\geq 2$, $\maxnorm\in \R_+$, $\stnd\in \R_+^n$ satisfying 
\begin{equation}\label{eq: condition-d-dmax}
\begin{split}
d\leq \dmax\leq \min(d^{\frac43},\maxnorm);\quad
\frac{h\dmax\,\log\log\log n}{\maxnorm/(1+(\log\log\log n)^{-1})}\geq 2;\\
\maxnorm\geq \|\stnd\|_1;\quad
\frac{\maxnorm/(1+(\log\log\log n)^{-1})}{C_{\text{\tiny \ref{l: net of majorizers}}}\,h\log\log\log n}\geq 1.
\end{split}
\end{equation}
Define
\begin{equation}\label{eq: definition-matrixset}
\matrixset:= \matrixset(n,\ell,d, \dmax, h, \maxnorm, \stnd)
\end{equation} as the set of 
all $n\times n$ symmetric matrices $M=(\mu_{ij})_{1\leq i,j\leq n}$ with zero diagonal satisfying the following conditions: 
\begin{itemize}
\item For any $v\in [n]$, $\degree_{G_M}(v)\leq \dmax$. 
\item $G_M$ is $\ell$-tangle free. 
\item Each non zero entry $\mu_{ij}$ of $M$ satisfies $\mu_{ij}^2\leq h$. 
\item All non-zero entries of $M$ are distinct (up to the symmetry constraint). 
\item For any $i\in [n]$, $\sum_{j=1}^n \mu_{ij}^2\leq \frac{\maxnorm}{1+(\log\log\log n)^{-1}} $. 
\item For any vertex $v\in [n]$, the number of its $\stnd$--heavy neighbors is at most $d^{\frac89}$. 
\end{itemize}

\bigskip

Following the previous section, define the discrete set
$$
\net_{mjr}:=\net\big(h,\maxnorm/(1+(\log\log\log n)^{-1}), \dmax, (\log\log\log n)^{-1}\big),
$$
where $\net(\cdot,\cdot,\cdot,\cdot)$ is taken from Lemma~\ref{l: net of majorizers}.
Note that, by the lemma and since $\dmax\leq \maxnorm$, we have 
\begin{equation}\label{eq: aux 5239565039853-50}
|\net_{mjr}|\leq e^{C_{mjr} \log^2 h\,(\log\log\log n)^3}
\end{equation}
for a universal constant $C_{mjr}>0$,
and for any vector $x\in\rearr(h,\maxnorm/(1+(\log\log\log n)^{-1}),\dmax)$ there is $y\in\net_{mjr}$ with $y\geq x$
and $\|y\|_1\leq \maxnorm$.

Define a mapping $\classif$ which will assign to every element of $\rearr(h,\maxnorm/(1+(\log\log\log n)^{-1}),\dmax)$
a majorizer from $\net_{mgr}$.
We will call $\classif$ {\it the standard classifier}.

\bigskip

Let $M\in \matrixset$, and denote by $G$ its associate graph $G_M=([n], E_M)$. 
We are interested in estimating the quantity 
\begin{equation}\label{eq: main quantity}
\sum_{\path} \prod_{i\edg j\in E_{\path}} \mu_{ij},
\end{equation}
where the sum is over all closed paths of length $2k$ on $K_{[n]}$, the complete graph on $n$ vertices. Note that any $\path$ having an edge not in $E_M$ do not contribute to the above sum, and we are therefore left with paths on $G$. To this aim, we will associate to each path a data structure which would help us encode the path and its contribution to the above quantity. 
Everywhere below, $\path:[0,2k]\to [n]$ is a path on $G$. 
We start by associating a diagram to this path, where by a {\it diagram} on an integer interval $[a,b]$ we understand a function $D:[a,b]\to\Z$
such that $D(t)-D(t-1)\in\{\pm 1\}$ for all $t\in[a+1,b]$.
The graph of $D$ consists of a sequence of ``diagonal up arrows'' and ``diagonal down arrows'' connecting 
the neighboring points $(t-1,D(t-1))$ and $(t,D(t))$ (``up arrow'' if $D(t)-D(t-1)=1$ or ``down arrow'' if $D(t)-D(t-1)=-1$). 
In the sequel, we will write $D(t-1)\up D(t)$ (resp. $D(t-1)\down D(t)$) to indicate a diagonal up (resp. down) arrow between the points $(t-1,D(t-1))$ and $(t,D(t))$. 

\medskip 
With each closed path $\path$ on $G$ of length $2k$, we associate a diagram $H_{\path}:[0,2k]\to\Z$ 
which can be iteratively constructed as follows.
First we set $H_\path(0)=0$ and 
for every $t\in[1,2k]$,
\begin{itemize}

\item If $\path(t-1)\edg \path(t)$ was not traveled before time $t-1$ then we set $H_{\path}(t):=H_{\path}(t-1)+1$;

\item If $\path(t-1)\edg \path(t)$ is a non-cycle edge of $G_{\path}$ which was traveled before time $t-1$, and the first time it was traveled in direction
$\path(t-1)\to \path(t)$, then we set $H_{\path}(t):=H_{\path}(t-1)+1$;

\item If $\path(t-1)\edg \path(t)$ is a non-cycle edge of $G_{\path}$ which was traveled before time $t-1$, and the first time it was traveled in direction
$\path(t)\to \path(t-1)$, then we set $H_{\path}(t):=H_{\path}(t-1)-1$;

\item If $\path(t-1)\edg \path(t)$ is a cycle edge of $G_{\path}$ which was traveled before time $t-1$ (in any direction) then we set
$H_{\path}(t):=H_{\path}(t-1)-1$;

\end{itemize}

Note that the value of $H_\path(t)$ at $2k$ is not necessarily zero: in general, the number of up-arrows and down-arrows do not agree.
However, we have
\begin{lemma}\label{lem: up-down-total}
Let $\path$ be a closed path on $G$ of length $2k$. 
Then 
for every non-cycle edge $e$ of $G_\path$,
the number of up arrows in $H_\path$ corresponding to $e$
is equal to the number of down arrows for $e$.
\end{lemma}
\begin{proof}
Let $e$ be a non-cycle edge of $G_\path$. Since $\path$ is closed, this condition
implies, in particular, that $e$ is traveled an even number of times, with each of the two possible directions
traveled the same number of times. Then the algorithm of constructing $H_\path$ implies
the desired conclusion for $e$.
\end{proof}

\bigskip
The diagram above will give us partial information on
whether we are discovering a new edge, or re-traversing it, as well as it's uncovering direction. 
To further encode where the path is heading at each step, we will be introducing {\it local indexation} of vertices. We will use three types of indexation.

Given a vertex $u\in [n]$ and its neighbor $v$ in $G$,
the {\it simple local index $\ind_u(v)$} of $v$ with respect to $u$
is the position of $v$ in the sequence of all neighbors of $u$
ordered according to the 
magnitude of $\mu_{uv}$;
formally
$$\ind_u(v):=\big|\big\{y\in[n]:\;|\mu_{uy}|\geq |\mu_{uv}|\mbox{ and $y$ is a neighbor of $u$ in $G$}\big\}\big|.$$

Further,
for every vertex $u\in [n]$,
define $\nbr_{HVY}(u)$ as the set of all $\stnd$--{\it heavy} neighbors $v$ of $u$ (as defined at the beginning of the section).
Then the {\it heavy-vertex local index} of $v\in\nbr_{HVY}(u)$
is the position of $v$ in the sequence of elements of $\nbr_{HVY}(u)$ arranged in 
the order determined by the magnitude of $\mu_{uv}$:
$$\dind_u(v):=\big|\big\{y\in\nbr_{HVY}(u):\;|\mu_{uy}|\geq |\mu_{uv}|\big\}\big|.$$

Unlike the first two, the third type of local indexation depends on $\path$.
Let $t\in[2k]$ and let
\begin{equation}\label{eq: aux 072-50962035350}
\begin{split}
\nbr_{\cycle}(\path,t):=
\big\{&\mbox{The set of all neighbors $v$ of $\path(t-1)$ in $G_\path$
such that}\\
&\mbox{either $v\edg \path(t-1)$ is a cycle edge 
with discovery time at most $t-1$}\\
&\mbox{or $v$ is a special cycle vertex in $G_{\path}$
with discovery time at most $t-1$}
\big\}.
\end{split}
\end{equation}
The {\it cycle local index} of a vertex $v\in \nbr_{\cycle}(\path,t)$
is defined by analogy with the above types of indexation and utilizes the same ordering:
$$\cind_\path(v,t):=\big|\big\{y\in \nbr_{\cycle}(\path,t):\;|\mu_{\path(t-1),y}|\geq |\mu_{\path(t-1),v}|\big\}\big|.$$
Note that the cycle local index of a vertex depends on time and can vary for different $t$.
We record the following simple lemma (which is a consequence of the estimates from Section~\ref{subs: stat of cycle vertices})
for future reference:
\begin{lemma}\label{l: neighborhoods}
Assuming that $M\in\matrixset(n,\ell,d, \dmax, h, \maxnorm, \stnd)$ (in particular, $G_M$
is $\ell$-tangle free for some $\ell\geq 5$) we have for any closed path $\path$ on $G_M$ of length $2k$ and
any time $t$:
$$
|\nbr_{\cycle}(\path,t)|\leq C_{\text{\tiny\ref{l: neighborhoods}}} k^3/\ell^3;
\quad |\nbr_{HVY}(\path(t))|\leq d^{\frac89},
$$
where $C_{\text{\tiny\ref{l: neighborhoods}}}>0$ is a universal constant.
\end{lemma}

\bigskip

For each closed path $\path: [0,2k]\to [n]$ on $G$, {\it we define
a data structure consisting of}
\begin{itemize}
\item The diagram $H_\path$;
\item The initial vertex $v_\path\in[n]$;
\item Sets of times $\cset_\path, \cycleedges_\path\subset [2k]$ and $\hdset_\path\subset[0,2k]$;
\item A {\it weight function} $\weight_\path: [2k]\to \Z$;
\item A vector-valued mapping $\badptclass_\path$.
\end{itemize}

The data structure is rather complex, and before giving a formal definition of its components let us briefly
describe their purpose.
The diagram $H_\path$ gives {\it partial} information on times when a new edge is discovered
as well as direction in which an edge is traveled at a given time.
The set $\cset_\path$ will be used to store information about special cycle vertices associated with the path.
Further, $\cycleedges_\path$ will encode information about the cycle edges of $G_\path$ having multiplicity either one
or at least three: there are several (technical) reasons for treating the cycle edges of multiplicity two differently.
The set $\badptclass_\path$ stores data about $\stnd$--heavy vertices; roughly speaking, about the rows
and columns of $M$ which have a big support, or a big $\ell_2$--norm, or an ``unusual'' profile.
The weight function $\weight_\path$ will store an index of a currently traveled edge. The type of indexation used
will depend, in particular, on the type of the vertex/edge.

Now, we turn to the formal description.
Construction of the diagram $H_\path$ was discussed above.
The initial vertex $v_\path$ is simply the vertex $\path(0)$ in $G$.

\medskip

{\bf Definition of $\cycleedges_\path$.} We set
\begin{align*}
\cycleedges_\path:=\big\{
t\in[2k]:\;&\path(t-1)\edg \path(t)\mbox{ is a cycle edge of $G_\path$ which is travelled by $\path$}\\
&\mbox{either a single time or at least three times}
\big\}.
\end{align*}

\medskip

{\bf Definition of $\cset_\path$.} For each $t\geq 1$, we add $t$ to the set $\cset_\path$ if one of the two conditions is satisfied:
\begin{itemize}
\item If $\path(t-1)$ is a special cycle vertex of $G_\path$ with discovery time at most $t-1$
and there is down arrow from $t-1$ to $t$ in the diagram $H_\path$, or
\item $\path(t)$ is a special cycle vertex of $G_\path$ with discovery time at most $t-1$
and there is an up arrow from $t-1$ to $t$ in the diagram.
\end{itemize}

\medskip

{\bf Definition of $\hdset_\path$.} For any $t\geq 0$, we add $t$ to the collection $\hdset_\path$
if the vertex $\path(t)$ is $\stnd$--heavy. 

\medskip

{\bf Definition of the weight function.} The weight function is constructed as follows:
\begin{itemize}
\item If ``$t\in \cset_\path\cup\cycleedges_\path$ and $\path(t-1)\edg\path(t)$ is a cycle edge of $G_\path$
with discovery time at most $t-1$''

or

``$\path(t)$ is a special cycle vertex of $G_\path$ with discovery time at most $t-1$'',

then $\weight_\path(t)$ is equal to $-\big(\cind_{\path}\big(\path(t),t\big)\big)$.
\item Otherwise, if $t\in \cset_\path$ and $\path(t-1)\edg\path(t)$ is either a non-cycle edge
or a cycle edge of $G_\path$ with discovery time at least $t$, then $\weight_\path(t)$ is equal to $0$.
\item Otherwise, if $t\in \hdset_\path$ and there is an up arrow from $t-1$ to $t$, then $\weight_\path(t)$ is equal
to $\dind_{\path(t-1)}\big(\path(t)\big)$, where the heavy-vertex local indexation is taken in $G$. 
\item Otherwise, if there is an up arrow (resp. down arrow) from $t-1$ to $t$,
then $\weight_\path(t)= \ind_{\path(t-1)}\big(\path(t)\big)$ (resp.\ $\weight_\path(t)=1$), where the simple local indexation is taken in $G$.
\end{itemize}
Note that the weight function is negative when the cycle indexation is used.
The purpose of this convention is to make sure that we can see that the local cycle indexation
is applied simply by looking at the value of $\weight_\path(t)$, regardless of the structure of $\path$.
This will be important below when discussing injectivity of our mapping.

\medskip

It will be convenient to define
$$
\hdset_\path^\down:=\big\{t\in\hdset_\path\cap[1,2k]:\; H_\path(t-1)\down H_\path(t)\big\},
\quad \hdset_\path^\up:=\hdset_\path\setminus\hdset_\path^\down.
$$
Note that with this definition
$$
\hdset^\up_\path=\big\{t\in\hdset_\path\cap[1,2k]:\; H_\path(t-1)\up H_\path(t)\big\}\cup(\hdset_\path\cap\{0\}),
$$
i.e.\ we interpret the initial time, if the starting vertex is $\stnd$--heavy, as an ``up-time''.
We will use similar notations $\cset^\up_\path$, $\cset^\down_\path$ for subsets of $\cset_\path$, although in that case the
definition is simpler since the zero time can never be included in $\cset_\path$:
$$
\cset^\up_\path:=\big\{t\in\cset_\path:\; H_\path(t-1)\up H_\path(t)\big\},\quad
\cset^\down_\path:=\big\{t\in\cset_\path:\; H_\path(t-1)\down H_\path(t)\big\}.
$$
We define by analogy the sets $\cycleedges_\path^\up,\cycleedges_\path^\down$.

\medskip

{\bf Definition of $\badptclass_\path$.}
First, define an auxiliary set
\begin{align}\label{eq: definition-V}
\V_\path:=\big\{&t\in \hdset^{\down}_\path\setminus (\cset_\path\cup \cycleedges_\path):\,
``H_{\path}(s_t)\neq H_{\path}(t)"\nonumber\\ 
& \text{ or } ``[s_t, t]\cap (\cset_\path\cup \hdset^\up_\path\cup \cycleedges_\path)\neq\emptyset"
\text{ where } s_t=\max\{t'\in \hdset_\path, t'<t\}\big\}
\end{align}
For a given $t\in\hdset^{\up}_\path\cup \V_\path$, let
$X^*$ be the non-increasing rearrangement of the vector $(\mu_{\path(t) i}^2)_{i=1}^n$.
Then we take $\badptclass_\path(t)\in\R^n$ as the value of the standard classifier $\classif(X^*)$
(defined at the beginning of this section).  
Thus, the mapping $\badptclass_\path$ will instruct us which majorizer we should take for a given heavy vertex.
The fact that $\badptclass_\path$ is not defined on the entire set $\hdset$ does not lead to a loss of information
because of our choice of the definition for $\V_\path$; in fact, given values of $\badptclass_\path$ on $\hdset^{\up}_\path\cup \V_\path$,
it is possible to deduce the values of $\classif(\cdot)$ for all heavy vertices (see Proposition~\ref{prop: properties-V}).

The reason why the mapping $\badptclass_\path$ is defined on $\hdset^{\up}_\path\cup \V_\path$
rather than on the entire set $\hdset_\path$ is that the latter would significantly increase the complexity
of our data space, making it ``too large'' to allow a satisfactory upper estimate for sums of the data weights
(we also refer to discussion in Section~\ref{s: overview}).
Specifically, the set $\hdset_\path^\down$ may have cardinality comparable to $k$ meaning that
there are $|\net_{mgr}|^{ck}$ distinct mappings from admissible realizations of $\hdset_\path$
into the net of majorizers $\net_{mgr}$. The data space defined in this manner would have cardinality
much larger than $4^k d^k$, which is unacceptable.
On the contrary, such a problem does not exist for the set $\hdset_\path^\up$
(which may also have cardinality of order $k$) since the increased complexity of large $\hdset_\path^\up$
is overcompensated by the decreased complexity of the weight function (which has a much smaller range on
points from $\hdset_\path^\up$, compared to ``regular'' vertices).
Finally, it can (and will) be shown
that the cardinality of the set $\V_\path$ can be bounded in terms of $|\cset_\path|$, $|\hdset_\path^\up|$, $|\cycleedges_\path|$,
allowing an efficient control of the size of the data space.

\medskip

The data structure associated to a path $\path$ will be written as
$\langle v_\path,H_\path, \cset_\path, \hdset_\path, \cycleedges_\path, \weight_\path,\badptclass_\path\rangle$, or, in the ``reduced'' form,
$\langle v_\path,H_\path, \cset_\path, \hdset_\path, \cycleedges_\path, \weight_\path\rangle$ (we remark here that $\badptclass_\path$ will be used
to count weights of paths, and is not employed for the rest of this section). 

The data structures will be used to estimate the number of distinct paths $\path$
on the graph $G$ as well as their contribution to \eqref{eq: main quantity}. To set up the relation between paths and data structures,
we will prove that the mapping $\path\to \langle v_\path,H_\path, \cset_\path, \hdset_\path, \cycleedges_\path, \weight_\path\rangle$ is injective.
For brevity, for the rest of the section we omit the subscript ``$\vphantom{|}_\path$'' for elements of the structure.
We start with an auxiliary lemma.

\begin{lemma}\label{l: cycle index compar}
Let $\path$ and $\path'$ be two paths mapped to the same structure
$\langle v,H, \cset, \hdset, \cycleedges, \weight\rangle$, and let $T\in[2k]$.
Assume additionally that $\path[T-1]=\path'[T-1]$.
Then
$$\nbr_{\cycle}(\path,T)=\nbr_{\cycle}(\path',T),$$
with the sets $\nbr_{\cycle}(\cdot)$ defined by \eqref{eq: aux 072-50962035350}.
\end{lemma}
\begin{proof}
Since $\path[T-1]=\path'[T-1]$,
any cycle edge of $G_\path$ with discovery time at most $T-1$ is also a cycle edge of $G_{\path'}$ with the same discovery time.
Similarly, every special cycle vertex of $G_{\path}$ with discovery time at most $T-1$ must be
a special cycle vertex of $G_{\path'}$ with the same discovery time. The result follows.
\end{proof}

\begin{prop}
The mapping $\path\to\langle v,H,\cset,\hdset,\cycleedges,\weight\rangle$
constructed above is injective.
\end{prop}
\begin{proof}
Denote by ${\bf Data}$ the mapping of paths to data structures.
Let $\langle v,H, \cset, \hdset,\cycleedges,\weight\rangle$ be a data structure in the range of ${\bf Data}$
and suppose there are two paths $\path$ and $\path'$ such that ${\bf Data}(\path)={\bf Data}(\path')
=\langle v,H, \cset, \hdset,\cycleedges, \weight\rangle$.
Our goal is to show that for any $t\in [0,2k]$, we have $\path(t)=\path'(t)$. 
We will prove the assertion by induction.

Clearly $\path(0)=\path'(0)$. Now, let $T\in[2k]$
and suppose that $\path(t)=\path'(t)$ for any $t\in [0,T-1]$. We will verify that $\path(T)=\path'(T)$,
by considering several cases which mirror the definition of the weight function $\weight$.
First of all, we need to make sure that the four cases in the definition of $\weight$ are matched
for both paths. For a path ${\widetilde\path}$, we say that
\begin{itemize}
\item {\it Condition (A) holds} if either `` $T\in \cset\cup\cycleedges$ and
${\widetilde\path}(T-1)\edg{\widetilde\path}(T)$ is a cycle edge of $G_{\widetilde\path}$
with discovery time at most $T-1$ ''
or `` ${\widetilde\path}(T)$ is a special cycle vertex
in $G_{\widetilde\path}$ with discovery time at most $T-1$ '';

\item {\it Condition (B) holds} if (A) does not hold and $T\in \cset$ and ${\widetilde\path}(T-1)\edg{\widetilde\path}(T)$ is either a non-cycle edge
or a cycle edge of $G_{\widetilde\path}$ with discovery time at least $T$;

\item {\it Condition (C) holds} if (A)--(B) do not hold and $T\in \hdset$ and there is an up arrow from $T-1$ to $T$ in $H$;

\item {\it Condition (D) holds} if (A)--(B)--(C) do not hold.

\end{itemize}

We then have that condition (A) (respectively, B, C or D) holds for $\path$
if and only if the same condition holds for $\path'$.
Indeed, by our convention, (A) holds if and only if the value of the weight function at time $T$ is negative,
and (B) holds if and only if the weight function is zero; similarly, the conditions (C) and (D) are ``path--independent''
i.e.\ are determined completely by the data structure.
Having this in mind, we now consider in detail each of the four conditions.

\begin{itemize}
\item[(A)] 
By the definition of $\weight(\cdot)$, in this case
$-\weight(T)= \cind_{\path}\big(\path(T),T\big)= \cind_{\path'}\big(\path'(T),T\big)$.
In view of Lemma~\ref{l: cycle index compar},
we have $\nbr_{\cycle}(\path,T)=\nbr_{\cycle}(\path',T)$,
and so matching cycle local indices imply that the corresponding vertices coincide: $\path(T)=\path'(T)$.

\item[(B)] In this case, $\weight(T)$ is equal to $0$.
Note that the definition of the set $\cset$ implies that there is a down arrow from $T-1$ to $T$ in the diagram $H$,
whence the edge $\path(T-1)\edg\path(T)$ is traveled before the time $T-1$,
and $\uncov_{\path}\big(\path(T-1)\edg \path(T)\big)=
\path(t)\to \path(T-1)$ for some $t\in [0,T-2]$.
Similarly, $\uncov_{\path'}\big(\path'(T-1)\edg \path'(T)\big)=
\path'(t')\to \path'(T-1)$ for some $t'\in [0,T-2]$.
However, there exists at most one vertex $u\in G_{\path[T-1]}$
such that $\uncov_{\path}(u\edg \path(T-1))=u\to \path(T-1)$ and $u\edg \path(T-1)$
is a non-cycle edge in $G_{\path[T-1]}$, and similarly for $G_{\path'[T-1]}$.
Since $G_{\path[T-1]}=G_{\path'[T-1]}$, we get $t=t'$, that is,
$\path(T)=\path'(T)$.

\item[(C)]
Here, $\weight(T)= \dind_{\path(T-1)}\big(\path(T)\big)=\dind_{\path'(T-1)}\big(\path'(T)\big)$. Since  
the heavy vertex local indexation is independent of a path, this immediately implies $\path(T)=\path'(T)$.

\item[(D)] If there is an up arrow from $T-1$ to $T$ in the diagram,
then the path-independent simple local indexation is used, so necessarily $\path(T)=\path'(T)$.
Finally, assume that there is a down arrow from $T-1$ to $T$, and $\weight(T)=1$.
It follows from the definition of $\cset$ that $\path(T-1)=\path'(T-1)$
cannot be a special cycle vertex in $G_\path$ (and in $G_{\path'}$)
with discovery time $T-1$ or less.
Further, it is clear that the edge $\path(T-1)\edg \path(T)$
exists in the graph $G_{\path[T-1]}=G_{\path'[T-1]}$, and similarly for $\path'(T-1)\edg \path'(T)$.
Let us assume that $\path(T)\neq \path'(T)$, and show that this leads to contradiction.
To emphasize that the paths coincide up to time $T-1$, we will use notation
$G_{T-1}$ for $G_{\path[T-1]}=G_{\path'[T-1]}$, and $p(t)$ instead of $\path(t)$ or $\path'(t)$
whenever $t\leq T-1$.
We have three subcases.

\begin{itemize}
\item[(i)] Assume that $p(T-1)\edg \path(T)$
is a cycle edge of $G_\path$ with a discovery time at most $T-1$,
and that, similarly, $p(T-1)\edg \path'(T)$
is a cycle edge of $G_{\path'}$ with a discovery time at most $T-1$.
Note that if $\uncov_{\path}(p(T-1)\edg \path(T))=p(T-1)\to \path(T)$
and $\uncov_{\path'}(p(T-1)\edg \path'(T))=p(T-1)\to \path'(T)$
then $p(T-1)$ is a splitting point of $G_{T-1}$, which is impossible.
Similarly, if $\uncov_{\path}(p(T-1)\edg \path(T))=\path(T)\to p(T-1)$
and $\uncov_{\path'}(p(T-1)\edg \path'(T))=\path'(T)\to p(T-1)$
then $p(T-1)$ is a cycle completion point in $G_{T-1}$, again leading to contradiction.
Thus, we can assume that $\uncov_{\path}(p(T-1)\edg \path(T))=\path(T)\to p(T-1)$
and $\uncov_{\path'}(p(T-1)\edg \path'(T))=p(T-1)\to \path'(T)$.
Let $t$ be the time such that the edge $p(T-1)\edg \path(T)$ is uncovered by $\path$
at $t+1$; and define time $t'$ for the edge $p(T-1)\edg \path'(T)$ by analogy.
We must have $t'>t$ since otherwise $p(T-1)$ would necessarily be a cycle completion point in $G_{T-1}$.
It is also clear that $t'+2\leq T-1$.
Let $\w t\geq t'+2$ be the smallest time such that $p(\w t)=p(T-1)$. 
Now, if we assume that the edge $p(\w t-1)\edg p(T-1)$
is distinct from both $p(T-1)\edg \path(T)$ and $p(T-1)\edg \path'(T)$
then $p(T-1)$ becomes a cycle meeting point in $G_{T-1}$, which is impossible. Hence, $p(\w t-1)\edg p(T-1)$
must be equal to one of the two edges. But then this (cycle) edge would be traveled by corresponding
path at least three times: $[t,t+1]$ (resp., $[t',t'+1]$ for $\path'$), $[\w t-1,\w t]$ and $[T-1,T]$,
contradicting the condition $T\notin\cycleedges$ (because we are not in situation A).
Thus, we showed that the situation when both $p(T-1)\edg \path(T)$ and $p(T-1)\edg \path'(T)$
are cycle edges of $G_\path$ (resp., $G_{\path'}$) with discovery times at most $T-1$,
is impossible.

\item[(ii)] Assume that $p(T-1)\edg \path(T)$
is either a non-cycle edge of $G_\path$ or a cycle edge with a discovery time at least $T$,
and that the same holds for $p(T-1)\edg \path'(T)$.
We cannot have simultaneously $\uncov_{\path}(p(T-1)\edg \path(T))=\path(T)\to p(T-1)$
and $\uncov_{\path'}(p(T-1)\edg \path'(T))=\path'(T)\to p(T-1)$ since this would
imply that at least one of the edges is a cycle edge in $G_{T-1}$.
Hence, we can assume that $\uncov_{\path}(p(T-1)\edg \path(T))=p(T-1)\to \path(T)$.
But then traveling from $p(T-1)$ to $\path(T)$ should correspond to an up-arrow in the diagram,
and we arrive at contradiction.

\item[(iii)] In the last scenario, we assume that $p(T-1)\edg \path(T)$
is either a non-cycle edge of $G_\path$ or a cycle edge with a discovery time at least $T$,
but $p(T-1)\edg \path'(T)$ is a cycle edge of $G_{\path'}$ with discovery time at most $T-1$.
Clearly, we must have $\uncov_{\path}(p(T-1)\edg \path(T))=\path(T)\to p(T-1)$ as otherwise
we would get an up-arrow in the diagram (we also recall that
the edge must be uncovered by the time $T-1$). 
Since $p(T-1)\edg \path(T)$ is not a cycle edge of $G_{T-1}$, then necessarily the uncovering of $p(T-1)\edg \path'(T)$ 
must happen after the uncovering of $p(T-1)\edg \path(T)$. Now since $p(T-1)$ is not a cycle completion point, 
we must have $\uncov_{\path'}(p(T-1)\edg \path'(T))=p(T-1)\to \path'(T)$. Since $T\not\in \cycleedges$
(otherwise, we would have been in situation A), then 
$p(T-1)\edg \path'(T)$ is only traversed twice by the time $T$:
at the moment of its uncovering and on the time interval $[T-1,T]$. Since there is a path connecting $\path'(T)$ to $p(T-1)$ 
after the uncovering of $p(T-1)\edg \path'(T)$ 
then necessarily $p(T-1)$ is a cycle completion point of $G_{T-1}$, 
leading to contradiction. 
\end{itemize}
\end{itemize}
\end{proof}

\section{Properties of the data structure associated with a path}\label{s: prop-data}

As in the second part of the previous section, in this section we omit the subscript ``$\vphantom{|}_\path$'' for elements of the data structure
corresponding to a path $\path$.
In this part of the paper,
we will establish some properties of the data structure defined in the previous section
by identifying constraints on the sets $\cset,\hdset$ and $\mathcal{C}$. 
First, the number of elements of $\cset\setminus \mathcal{C}$ corresponding to up arrows in the diagram cannot be much
smaller than the number of elements in $\cset\setminus \mathcal{C}$ corresponding to down arrows:
\begin{prop}\label{prop: up minus down}
Let $\path$ be a closed path on an $\ell$-tangle free graph $G$ ($\ell\geq 5$),
and let $\cset, \mathcal{C}$ and $H$ be the associated sets and diagram from the data structure.
Then $|\cset^\up\setminus \mathcal{C}|\geq |\cset^\down\setminus \mathcal{C} |-C|E_{\path}|^4/\ell^4$ for some universal constant $C>0$.
\end{prop}
\begin{proof}
Recall that a time $t$ is added to $\cset^\down$ only if $\path(t-1)$
is a special cycle vertex of $G_\path$ with discovery time at most $t-1$. 
Denote by $M'$ the set of all cycle meeting points, splitting points and cycle completion points in $G_\path$.
Then, by combining Lemmas~\ref{lem: nb-cycle meeting points},~\ref{l: number of splitting}
and~\ref{l: n of cycle compl},
we get that the cardinality of $M'$ is bounded above by $C|E_{\path}|^3/\ell^3$ for some universal constant $C>0$.
Let $S$ be the set of all distinct cycle edges $e$ of $G_\path$ such that there is $t\in \cset^\down$
with $e=\path(t-1)\edg\path(t)$. 
By Corollary~\ref{cor: nb-cycles meeting},
there are at most $C'|E_{\path}|/\ell$ cycle edges incident to a given point in $M'$,
whence the total number of distinct edges in $S$ is at most
$\tilde C|E_{\path}|^4/\ell^4$. For any such edge $e$, 
we have 
$$|\{t\in \cset^\down\setminus \mathcal{C}^{\down}:\;e=\path(t-1)\edg\path(t)\}|
=1,
$$
since $t\not\in \mathcal{C}^{\down}$ automatically implies that $e$ is a cycle edge of multiplicity $2$ and only it's second appearance is recorded with a down arrow. 
Therefore, we deduce that 
$$
|\{t\in \cset^\down\setminus \mathcal{C}^{\down}:\; \path(t-1)\edg\path(t) \text{ is a cycle edge}\}|\leq \tilde C|E_{\path}|^4/\ell^4. 
$$

Next, for any vertex $v\in M'$,
let $T_v$ to be the set of all times
$t\in \cset^\down$ such that $v=\path(t-1)$ and $\path(t-1)\edg\path(t)$ is a non-cycle edge of $G_{\path}$.
Note that the set  of non-cycle edges of $G_{\path}$ incident to $v$ must be traveled in alternating directions
(towards $v$, then from $v$, etc.).
Then necessarily
$$|T_v|\leq \big|\big\{t\in \cset^\up:\;v=\path(t)
\mbox{ and $\path(t-1)\edg\path(t)$ is a non-cycle edge of $G_\path$}\big\}\big|+1.$$
Therefore, we get 
$$
|\{t\in \cset^\down\setminus \mathcal{C}^{\down}:\; \path(t-1)\edg\path(t) \text{ is a non-cycle edge}\}|\leq \vert \cset^{\up}\setminus \mathcal{C}\vert+ C|E_{\path}|^3/\ell^3. 
$$
Combining the above estimates, we get
the result.

\end{proof}

In the next proposition, we show that
within time intervals not containing any elements of $\cset\cup\mathcal{C}$
and such that the corresponding part of the diagram attains its minima at the interval endpoints,
these endpoints correspond to the same vertex of the graph. 
\begin{prop}\label{prop: above-B-point}
Let $G$, $\path$ be as before, and let $\hdset$ be the corresponding set from the data structure for $\path$.
Let $t<t'$ and assume that, first, $[t+1,t']\cap (\cset\cup\mathcal{C})=\emptyset$;
second, $H(t)=H(t')$; third, $H(\tau)\geq H(t)$ for all $\tau\in[t,t']$.
Then $\path(t)=\path(t')$. 
In particular,  either both $t,t'$ belong to $\hdset$ or $\{t,t'\}\cap \hdset=\emptyset$.
\end{prop}
\begin{proof}
We will show that $\path(t)=\path(t')$; then the second assertion of the proposition will follow automatically. 
In turn, to verify this property, it is sufficient to show that whenever
$t<t_1\leq t_2<t'$ are such that $\path(t_1)=\path(t_2)$ and $H(t_1-1)\up H(t_1)$
and  $H(t_2)\down H(t_2+1)$ then necessarily $\path(t_1-1)=\path(t_2+1)$. 
Indeed, an inductive argument based on the latter property will lead to the result. 
 
Assume the opposite i.e. $\path(t_1-1)\neq\path(t_2+1)$ and denote  $u=\path(t_1-1)$, $v=\path(t_1)=\path(t_2)$ and $w= \path(t_2+1)$. Since $H(t_2)\down H(t_2+1)$, then the edge $v\edg w$ appears previously in the path. Let $t_0$ be the time the edge $v\edg w$ was uncovered. 
We will consider several cases: 
\begin{itemize}
\item[(a)] Suppose $t_0<t_1-1$. We have either $\path(t_0-1)=v$ or $\path(t_0)=v$.
Since $H(t_1-1)\up H(t_1)$ then the edge $u\edg v$ wasn't uncovered before the time $t_1-1$.
Therefore, $v$ is connected to $u$ through a path in $G_{\path([t_1-1])}$ which doesn't contain the edge $u\edg v$.
Hence, a new cycle is formed at time $t_1$, making $v=\path(t_1)=\path(t_2)$
a cycle completion point. But then, by the definition of $\cset$, the time $t_2+1$ must belong to $\cset$
which contradicts the hypothesis $[t+1,t']\cap \cset=\emptyset$. 
\item[(b)] Suppose $t_1< t_0$ and that $\uncov(v \edg w)= v\to w$. Therefore $\path(t_0-1)= v$ and $\path(t_0)=w$,
and $v\edg w$ is necessarily a cycle edge (of multiplicity $2$) in $G_\path$.
Let $\tau\geq t_0$ be the first time when a cycle containing the edge $v\edg w$ is completed.
Clearly, $\tau< t_2$ (for $\tau=t_2$, $v$ would become a cycle completion point, and we would get $t_2+1\in\cset$).
Denote by $\bf C$ a cycle in $G_\path$ containing the edge $v\edg w$ and completed
at time $\tau$, and let $q_1$, $q_2$ be the vertices on the cycle which are neighbors of $\path(\tau)$.
Note that traveling $\path(\tau)\to q_1$ and $\path(\tau)\to q_2$
within time interval $[\tau,t_2]$ is prohibited because in that case we would produce a point in $\cset$.
On the other hand, denoting by $\tau'$ the first time in $[\tau+1,t_2]$ when $\path(\tau')=v$,
the edge $\path(\tau'-1)\edg v$ must necessarily be a cycle edge which belongs to $\bf C$
as otherwise $v$ would turn into a cycle meeting point with discovery time at most $t_2$, implying $t_2+1\in\cset$.
These observations, combined together, imply that there is time $\tau''\in[\tau+1,t_2-1]$ such that
two conditions hold simultaneously: first, $\path(\tau'')$ is a cycle meeting point
with discovery time at most $\tau''$; second,
the edge $\path(\tau'')\edg\path(\tau''+1)$ belongs to $\bf C$. But then $\tau''+1\in\cset$ --- a contradiction.
\item[(c)] Suppose $t_1< t_0$ and that $\uncov(v \edg w)= w\to v$. Since $\path(t_1)=v$, $\path(t_0-1)=w$, $\path(t_0)=v$ and $H(t_0-1)\up H(t_0)$, then $v$ is connected to $w$ in $G_{\path([t_0-1])}$ through the path without traversing the edge $v\edg w$. Since $\path(t_0-1)\to \path(t_0)$, then necessarily $\path(t_0)$ is a cycle completion point, 
leading again to a contradiction. 
\end{itemize}
\end{proof}

In the second part of the section, we connect properties of the diagram $H$ with
some structural properties of the sets $\cset$, $\hdset$ and $\mathcal{C}$. 

\begin{prop}\label{prop: nb-levels}
Let $G$ and $\path$ be as before; let $H$, $\cset$, $\hdset$ and $\cycleedges$ be the corresponding
elements of the data structure associated with $\path$, and
let 
$$
\hdset'=\{t\in \hdset:\,  H(t)\down H(t+1)\}.
$$
Then
$$
|\hdset'|\leq 
3|\cycleedges^\down|+|\hdset^\up|+1.
$$
\end{prop}
\begin{proof}
We start the proof by considering an arbitrary heavy vertex visited by the path, and
will estimate some associated quantities.
 
Take any $\stnd$--heavy vertex $v$, and define
$B'_v:=\{t\in \hdset:\, \path(t)=v,\; H(t)\down H(t+1)\}$.
Further, let $t_1<\dots<t_{u}$ ($u\geq 1$) be all the times when $v$ is visited by $\path$.
Define a function $f_v$ on the collection $t_j$, $j\leq u$, by setting
$$f_v(t_j):=\big|\big\{e:\,\mbox{$e$ is an edge of $G_{\path[t_j]}$ incident to $v$ and having multiplicity one}\big\}\big|,\;\;j\leq u.$$
Clearly, $f_v(t_1)\leq 1$ and $|f_v(t_{j+1})-f_v(t_j)|\leq 2$ for all $j\leq u-1$.
Finally, define
\begin{align*}
Q_v^\up&:=\big\{t\in[1,2k]:\,\path(t)=v,\;H(t-1)\up H(t)\big\};\\
Q_v^\down&:=\big\{t\in[1,2k]:\,\path(t)=v,\;H(t-1)\down H(t)\big\}.
\end{align*}
Obviously, $Q_v^\up\subset\hdset^\up$.
Further, assume that $1\leq r\leq u-1$ is such that both $t_r\in B'_v$ and $t_{r+1}\in Q_v^\down$.
Consider several cases:
\begin{itemize}
\item The edge $\path(t_r)\edg \path(t_r+1)$ is either a non-cycle edge or a cycle edge with discovery time at least $t_{r+1}+1$.
Then necessarily $\path(t_{r+1}-1)\edg \path(t_{r+1})=\path(t_r)\edg \path(t_r+1)$,
and one of the two times the edge is traveled within the time interval $[t_r,t_{r+1}]$, it must correspond to an up-arrow
in the diagram, leading to contradiction.
\item The edge $\path(t_r)\edg \path(t_r+1)$ is a cycle edge with discovery time at most $t_{r+1}$,
and $\path(t_r)\edg \path(t_r+1)=\path(t_{r+1}-1)\edg \path(t_{r+1})$.
Since $H(t_r)\down H(t_r+1)$, the cycle edge is traveled at least three times by the time $t_{r+1}$,
so that $t_{r+1}\in \cycleedges^\down$.
\item The edge $\path(t_r)\edg \path(t_r+1)$ is a cycle edge with discovery time at most $t_{r+1}$,
and $\path(t_r)\edg \path(t_r+1)\neq\path(t_{r+1}-1)\edg \path(t_{r+1})$.
Then both are cycle edges with discovery times at most $t_{r+1}$,
and, since $H(t_r)\down H(t_r+1)$ and $H(t_{r+1}-1)\down H(t_{r+1})$,
we have either $\{t_r+1,t_{r+1}\}\cap \cycleedges^\down\neq \emptyset$
or $f(t_{r+1})\leq f(t_r)-2$.
\end{itemize}
To summarize, whenever $t_r\in B'_v$ and $t_{r+1}\in Q_v^\down$, we must have either
$\{t_r+1,t_{r+1}\}\cap \cycleedges^\down\neq \emptyset$ or $f(t_{r+1})\leq f(t_r)-2$.
At the same time, it is not difficult to see that whenever $f(t_{z+1})\geq f(t_z)+1$ for some $z\le u-1$,
we must have $t_{z+1}\in \cycleedges^\down$.
Together with the simple properties of the function $f_v$ mentioned above, this yields
$$
\big|\big\{r\leq u-1:\,\mbox{$t_r\in B'_v$ and $t_{r+1}\in Q_v^\down$}\mbox{ and }
\{t_r+1,t_{r+1}\}\cap \cycleedges^\down= \emptyset\big\}\big|
\leq \big|\{t_j\}_{j=1}^u\cap \cycleedges^\down\big|,
$$
implying that
$$
\big|\big\{r\leq u-1:\,\mbox{$t_r\in B'_v$ and $t_{r+1}\in Q_v^\down$}\big\}\big|\leq 2\big|\{t_j\}_{j=1}^u\cap \cycleedges^\down\big|
+\big|\{t_{j}+1\}_{j=1}^u\cap \cycleedges^\down\big|.
$$
The inclusion $Q_v^\up\subset\hdset^\up$ then gives
$$
|B'_v|-1\leq\big|\big\{r\leq u-1:\,t_r\in B'_v\big\}\big|\leq 2\big|\{t_j\}_{j=1}^u\cap \cycleedges^\down\big|
+\big|\{t_{j}+1\}_{j=1}^u\cap \cycleedges^\down\big|+\big|\{t_j\}_{j=2}^u\cap \hdset^\up\big|.
$$
Finally, whenever $t_1>0$, we necessarily have $t_1\in \hdset^\up$, whence
the above relation can be strengthened to
$$
|B_v'|\leq 2\big|\{t_j\}_{j=1}^u\cap \cycleedges^\down\big|
+\big|\{t_{j}+1\}_{j=1}^u\cap \cycleedges^\down\big|+\big|\{t_j\}_{j=1}^u\cap \hdset^\up\big|.
$$
It remains to apply the estimate for all $\stnd$--heavy vertices to get the result.
\end{proof} 

To end this section, we record some properties on the set $\V$ defined in \eqref{eq: definition-V} and used in the definition of $\badptclass$. 

\begin{prop}\label{prop: properties-V}
Let $G$, $\path$ be as before, and let $H, \cset, \hdset,\cycleedges$ be the corresponding elements from the data structure for $\path$. 
Then we have 
\begin{itemize}
\item[(i)] $\vert \V\vert \leq 
5\vert\cset\cup \hdset^{\up}\cup \cycleedges\vert+1$. 
\item[(ii)] For any $t\in \hdset^{\down}\setminus (\cset\cup \cycleedges)$, we have $[t]\cap \V\neq \emptyset$.
\item[(iii)] For any $t\in \hdset^{\down}\setminus (\cset\cup \cycleedges)$, we have $\path(\kappa_t)=\path(t)$ where $\kappa_t=\max\{t':\, t'\in [t]\cap \V\}$. 
\end{itemize}
\end{prop}
\begin{proof}\ \\
\begin{itemize}
\item[(i)] 
Let $\hdset'$ be defined as in Proposition~\ref{prop: nb-levels}, and
let us define an injective map $f: \V\to (\cset\cup \hdset^{\up}\cup \hdset' \cup \cycleedges)$ as follows: 
Take any $t\in \V$, and
let $s_t:=\max\{t'\in \hdset, t'<t\}$ (note that $s_t$ is well defined). If $[s_t,t]\cap (\cset\cup \hdset^{\up}\cup\cycleedges)\neq \emptyset$, then,
in view of the definition of $\V$, $[s_t,t-1]\cap (\cset\cup \hdset^{\up}\cup\cycleedges)\neq \emptyset$,
and we can set $f(t):= \max\{t'\in [s_t,t-1]:\, t'\in \cset\cup \hdset^{\up}\cup \cycleedges\}$; otherwise, we set $f(t):=s_t$. 
Clearly $f$ is injective since $f(t)\in [s_t,t-1]$ for any $t\in \V$, and we are left to check that 
the range of $f$ is a subset of $\cset\cup \hdset^{\up}\cup \hdset'\cup \cycleedges$. 

Let $t\in \V$. If $[s_t,t]\cap (\cset\cup\hdset^{\up}\cup \cycleedges)\neq \emptyset$, then clearly $f(t)\in \cset\cup \hdset^{\up}\cup \cycleedges$. Now, if $[s_t,t]\cap (\cset\cup \hdset^{\up}\cup\cycleedges)= \emptyset$, then 
since $t\in \V$ we  have $H(s_t)\neq H(t)$. Using Proposition~\ref{prop: above-B-point} together with the definition of $s_t$
and the condition $H(t-1)\down H(t)$,
we get that necessarily $H(s_t)>H(t)$ and $H(s_t)\down H(s_t+1)$ meaning that $s_t\in \hdset'$. 
This proves our assertion about $f$.

Finally, since $f$ is well defined and injective, we deduce that $\vert \V\vert \leq \vert (\cset\cup \hdset^{\up}\cup \hdset' \cup\cycleedges)\vert$. It remains to apply Proposition~\ref{prop: nb-levels}.
 
\item[(ii)] Let $t\in \hdset^{\down}\setminus (\cset\cup \cycleedges)$ and suppose that $[t]\cap \V= \emptyset$. In particular, $t\not\in \V$. Then necessarily $H(s_t)=H(t)$ and $[s_t, t]\cap (\cset\cup \hdset^\up\cup \cycleedges)=\emptyset$,
whence $s_t\in \hdset^{\down}\setminus (\cset\cup \cycleedges)$. By our hypothesis,
$s_t\not\in\V$, so we can repeat the argument.
Continuing the same reasoning, we deduce that for any $t'\in \hdset$ with $t'\leq t$, we have 
$$
t'\in \hdset^{\down},\, [t',t]\cap (\cset\cup \cycleedges)=\emptyset,\, \text{ and } H(t)= H(t').
$$
Proposition~\ref{prop: above-B-point} implies that $\path(t')= \path(t)$ for any $t'<t$ with $t'\in \hdset$. This means that the discovery of the vertex $\path(t)$ was recorded with a down arrow, which contradicts our construction of the diagram $H$. 
\item[(iii)] The proof follows by applying the above procedure between $\kappa_t$ and $t$. 
\end{itemize}
\end{proof}

\section{Summing over the weight functions}

Let $M=(\mu_{ij})\in \matrixset$ be as defined after formula \eqref{eq: definition-matrixset}
(with the parameters satisfying \eqref{eq: condition-d-dmax}) and let $G=G_M=([n],E_M)$ be the associated
simple graph on $[n]$. 
Let $k\geq 1$ and let $\path$ be a closed path on $K_{[n]}$ of length $2k$.
We define the weight of the path as 
\begin{equation}\label{eq: def-path-weight}
\pathw_M(\path):=\prod\limits_{t=1}^{2k}\mu_{\path(t-1),\path(t)}.
\end{equation}
Our goal is to estimate from above the quantity
$$\sum\limits_{\path}\pathw_M(\path),$$
where the sum is taken over all closed paths of length $2k$ on $K_{[n]}$.
Recall that a part of our strategy is to replace the above sum with summation of data structure weights
over the data space (see Section~\ref{s: overview}).
At this stage, we are ready to define precisely the weight of a data structure.
It is given as the right hand side of the relation in the lemma below:

\begin{lemma}\label{l: path major}
Let parameters $n,\ell, d,\dmax\in \N$, $h,\maxnorm\in \R_+$ and $\stnd\in \R_+^n$ satisfy \eqref{eq: condition-d-dmax}, and let
$M\in \matrixset(n,\ell,d, \dmax, h,\maxnorm, \stnd)$ (see \eqref{eq: definition-matrixset} for definition). 
Let $\path$ be a closed path on $G_M$ of length $2k$, and let $\langle v,H, \cset, \hdset, \cycleedges, \weight,\badptclass\rangle$
be the corresponding data structure. Then we have
\begin{align*}
\pathw_M(\path)\leq &\prod\limits_{t\in{\mathcal C}}\sqrt{h}\cdot
\prod\limits_{\substack{t\in (\cset\cup\hdset)\setminus {\mathcal C}:\\ t\geq 1,\,H(t-1)\up H(t)}}h\cdot
\prod\limits_{\substack{t\in [2k]\setminus(\cset\cup\hdset\cup\cycleedges):\\ t-1\in\hdset,\,H(t-1)\up H(t)}}\badptclass(f(t-1))_{\weight(t)}
\cdot\prod\limits_{\substack{t\in [2k]\setminus(\cset\cup\hdset\cup\cycleedges):\\ t-1\notin\hdset,\,H(t-1)\up H(t)}}\stnd_{\weight(t)},
\end{align*}
where $f(t-1):=t-1$ if $t-1\in \cset\cup \hdset^{\up}\cup \cycleedges$ and $f(t-1):=\max\{t':\, t'\in [t-1]\cap \V\}$ otherwise. 
\end{lemma}
\begin{proof}
We will prove the bound by considering every edge individually.
First, if $t\in\cycleedges$ (i.e.\ $\path(t-1)\edg\path(t)$ is a cycle edge of multiplicity not equal to two) then we trivially
bound $|\mu_{\path(t-1),\path(t)}|$ by $\sqrt{h}$.

Next, let $e=i\edg j$ be an edge in $G_\path$ which is not in $\cycleedges$,
so that $\path$ travels along $e$ an even number of times, say, $2m$.
Then, by Lemma~\ref{lem: up-down-total}, we have
\begin{equation}\label{eq: aux 2352308765235}
\mu_{e}^{2m}=\prod\limits_{t\in[2k]:\,\path(t-1)\edg\path(t)=e}\mu_{\path(t-1),\path(t)}
=\prod\limits_{t\in[2k]:\,\path(t-1)\edg\path(t)=e,\,H(t-1)\up H(t)}\mu_{\path(t-1),\path(t)}^2.
\end{equation}
When $t\in (\cset\cup\hdset)\setminus \cycleedges$,  we will bound the quantity above by replacing $\mu_{\path(t-1),\path(t)}^2$ with $h$. 

Further, for any $t\notin\cset\cup\hdset\cup\cycleedges$ such that $t-1\in \hdset$ and $H(t-1)\up H(t)$, we have $\path(t-1)=\path(f(t-1))$.
Indeed, if  $t-1\in \cset\cup \hdset^{\up}\cup \cycleedges$, then $f(t-1)=t-1$ and there is nothing to prove.
Otherwise, if $t-1\in \hdset^{\down}\setminus (\cset\cup \cycleedges)$,
then by Proposition~\ref{prop: properties-V} we have $\path(t-1)= \path(f(t-1))$. 
Therefore, in any case, the non-increasing rearrangement of the vector $(\mu_{\path(t-1),i}^2)_{i=1}^n$
is majorized by the vector $\badptclass(f(t-1))$, by the definition of $\badptclass$.  
Hence, taking into account the definition of the weight function $\weight$, we get
$\mu_{\path(t-1),\path(t)}^2\leq \badptclass(f(t-1))_{\weight(t)}$.

Finally, in the remaining case $t\notin\cset\cup\hdset\cup\cycleedges$, $t-1\notin\hdset$, $H(t-1)\up H(t)$,
we use that $\path(t-1)$ is not $\stnd$--heavy and thus $\stnd$ dominates the rearrangement of the vector
$(\mu_{\path(t-1),i}^2)_{i=1}^n$.

In view of formula \eqref{eq: aux 2352308765235}, this implies the statement of the lemma. 
\end{proof}

As the next (crucial) step, we will compute the contribution of paths sharing the same data structure up to the
realization of the weight function (in a sense, we integrate over the weight function).

\begin{prop}\label{prop: path-weight-1}
Let $n,\ell,d, k,\dmax\in \N$, $h,\maxnorm\in \R_+$ and $\stnd\in \R_+^n$ 
satisfy \eqref{eq: condition-d-dmax}, with $k^3/\ell^3\leq d^{\frac89}$,
and let $M\in \matrixset(n,\ell,d, \dmax, h,\maxnorm, \stnd)$ (see \eqref{eq: definition-matrixset} for definition). 
Fix a vertex $v\in [n]$,
a diagram $H$, subsets $\cset,\hdset, \cycleedges$ and a vector sequence $\badptclass$. 
Denote by $\textbf{P}$ the collection of all closed paths $\path$ of length $2k$ on $G_M$ with corresponding data structures of the form
$\langle v,H, \cset, \hdset, \cycleedges, \cdot, \badptclass\rangle$. 
Then we have 
\begin{align*}
\sum\limits_{\path\in \textbf{P}}\pathw_M(\path)
\leq
&h^{|\cycleedges|/2}(2\dmax)^{|\cycleedges^\up|}\,
\big(C_{\text{\tiny\ref{prop: path-weight-1}}}k^3/\ell^3\big)^{|\cycleedges^\down\cup \cset^\down|}
\big(C_{\text{\tiny\ref{prop: path-weight-1}}}d^{\frac89}h\big)^{|(\cset^\up\cup\hdset^\up)\setminus\cycleedges^\up|}\\
&\cdot\maxnorm^{|\{t\in[2k]:\,t\notin\cset\cup\hdset\cup \cycleedges,\, t-1\in\hdset,\,H(t-1)\up H(t)\}|}
\|\stnd\|_1^{|\{t\in[2k]:\,t\notin\cset\cup\hdset\cup\cycleedges,\, t-1\notin\hdset,\,H(t-1)\up H(t)\}|},
\end{align*}
where $C_{\text{\tiny\ref{prop: path-weight-1}}}>0$ is a universal constant.
\end{prop}
\begin{proof}
Denote by $R^1$ the set of all times $t\geq 1$ with $t\in (\cset\cup\hdset)\setminus\cycleedges$ and $H(t-1)\up H(t)$.
Further, let $R^2$ be the collection of all times $t\geq 1$ such that $t\notin\cset\cup \hdset\cup\cycleedges$,
$t-1\in\hdset$, and $H(t-1)\up H(t)$;
and let $R^3$ be the set of all $t\geq 1$ with $t\notin\cset\cup\hdset\cup\cycleedges$, $t-1\notin\hdset$, and $H(t-1)\up H(t)$.
For each $\path\in\textbf{P}$, denote by $\weight_\path$ the weight function for $\path$.
By Lemma~\ref{l: path major}, for any path $\path\in \textbf{P}$ we have
\begin{align*}
\pathw_M(\path)\leq \prod\limits_{t\in{\cycleedges}}\sqrt{h}\cdot\prod\limits_{t\in R^1}h\cdot
\prod\limits_{t\in R^2}\badptclass(f(t-1))_{\weight_\path(t)}
\cdot
\prod\limits_{t\in R^3}\stnd_{\weight_\path(t)}.
\end{align*}
For each $\path\in \textbf{P}$, let $S_\path$ be the sequence $(\weight_{\path}(t))_{t\in R^2\cup R^3}$.
The mapping $\path\to S_{\path}$ on $\textbf{P}$ is not injective in general. Let $\alpha$ be the maximal cardinality of the preimage of
a sequence under this mapping, i.e.\ let $\alpha:=\max\limits_{S}|\{\path\in \textbf{P}:\;S_\path=S\}|$.
Then it is not difficult to verify that 
\begin{align*}
\sum\limits_{\path\in \textbf{P}}\pathw_M(\path)\leq \alpha\prod\limits_{t\in{\cycleedges}}\sqrt{h}\cdot\prod\limits_{t\in R^1}h\cdot 
\sum\limits_{S_\path}
\prod\limits_{t\in R^2}\badptclass(f(t-1))_{\weight(t)}
\cdot
\prod\limits_{t\in R^3}\stnd_{\weight(t)},
\end{align*}
where the summation is taken over all {\it admissible} sequences $S_\path:R^2\cup R^3\to\N_0$
(i.e.\ subsequences of weight functions of some paths in $\textbf{P}$).
Taking the summation inside the product, we get
\begin{align*}
\sum\limits_{\path\in \textbf{P}}\pathw_M(\path)&\leq \alpha\,h^{|\cycleedges|/2+|R^1|} \cdot
\prod\limits_{t\in R^2}\Big(\sum\limits_{i=1}^n\badptclass(f(t-1))_i\Big)
\cdot
\prod\limits_{t\in R^3}\Big(\sum\limits_{i=1}^{n}\stnd_{i}\Big)\\
&\leq \alpha\,h^{|\cycleedges|/2+|R^1|}
\maxnorm^{|R^2|}\|\stnd\|_1^{|R^3|},
\end{align*}
where the last relation follows from the definition of $\badptclass$ and 
$\matrixset(n,\ell,d, \dmax, h,\maxnorm, \stnd)$.

Thus, in order to prove the result it remains to estimate $\alpha$.
Fix any sequence $S$ indexed over $R^2\cup R^3$.
Since the mapping of paths to data structures is injective,
it is enough to obtain an upper bound on the cardinality of the set $\{(W_\path(t))_{t\in[0,2k]\setminus (R^2\cup R^3)}:\;\path\in\textbf{P}\}$.
Since $G_M$ is $\ell$-tangle free, then using the definition of a weight function, we get that
for any $t\in\cycleedges^\down$, there are at most $C_{\text{\tiny\ref{l: neighborhoods}}}k^3/\ell^3$
admissible realizations of $\weight_\path(t)$,
so the total number of admissible realizations of $(\weight_\path(t))_{t\in\cycleedges^\down}$, $\path\in\textbf{P}$,
is at most $(C_{\text{\tiny\ref{l: neighborhoods}}}k^3/\ell^3)^{|\cycleedges^\down|}$
(we recall that, for $t\in\cycleedges^\down$, the weight function at $t$ is equal to the negative of the local cycle index
of the vertex $\path(t)$, which can be bounded using Lemma~\ref{l: neighborhoods}).
Further, for every
$t\in\cycleedges^\up$ there can be at most $2\dmax$ realizations of $\weight_\path(t)$,
since each row in $M$ has at most $\dmax$ non-zero entries, and here we count for the possibility
of $\weight_\path(t)$ taking a negative value if the local cycle indexation is used.
Further, since for any $v\in [n]$, the number of its $\stnd$--heavy neighbors is at most $d^{\frac89}$, then there are at most $(d^{\frac89}
+C_{\text{\tiny\ref{l: neighborhoods}}}k^3/\ell^3)^{|R^1|}$ possible realizations of the sequence
$(\weight_\path(t))_{t\in R^1}$, $\path\in \textbf{P}$, where we count for the possibility of the weight function
taking negative values if the local cycle indexation is used, and apply Lemma~\ref{l: neighborhoods}. Next,
since for any $t\in\cset^\down\setminus \cycleedges^{\down}$
we either use the local cycle indexation or assign weight $0$ (see the definition of the weight function),
then, in view of Lemma~\ref{l: neighborhoods},
there are at most $1+C_{\text{\tiny\ref{l: neighborhoods}}}k^3/\ell^3$ realizations of
$\weight_\path(t)$.  Therefore, the total number of admissible
realizations of $(\weight_\path(t))_{t\in\cset^\down\setminus \cycleedges^\down}$
is at most $(1+C_{\text{\tiny\ref{l: neighborhoods}}}k^3/\ell^3)^{|\cset^\down\setminus \cycleedges^\down|}$.
Finally, note that for all $t\in[2k]\setminus (\cset^\down\cup \cycleedges^\down)$ with $H(t-1)\down H(t)$
we have $\weight_\path(t)=1$. Therefore,
$$\alpha\leq (2\dmax)^{|\cycleedges^\up|}\, (1+C_{\text{\tiny\ref{l: neighborhoods}}}k^3/\ell^3)^{|\cycleedges^\down\cup \cset^\down|}
(d^{\frac89}+C_{\text{\tiny\ref{l: neighborhoods}}}k^3/\ell^3)^{|R^1|},$$
and
$$\sum\limits_{\path\in \textbf{P}}\pathw_M(\path)\leq
h^{|\cycleedges|/2+|R^1|}(2\dmax)^{|\cycleedges^\up|}\,
\big(1+C_{\text{\tiny\ref{l: neighborhoods}}}k^3/\ell^3\big)^{|\cycleedges^\down\cup \cset^\down|}
\big(d^{\frac89}+C_{\text{\tiny\ref{l: neighborhoods}}}k^3/\ell^3\big)^{|R^1|}
\maxnorm^{|R^2|}\|\stnd\|_1^{|R^3|}.$$
Using the definition of $R^1,R^2,R^3$, we get the result.
\end{proof}

\medskip

Crucially, the above proposition allows to estimate the sum of paths' weights over all admissible
paths via the sum over {\it admissible} data structures;
specifically, we can write
\begin{align*}
\sum\limits_{\path}\pathw_M(\path)
\leq
&\sum\limits_{v,H,\cset,\hdset,\cycleedges,\badptclass}h^{|\cycleedges|/2}(2\dmax)^{|\cycleedges^\up|}\,
\big(C_{\text{\tiny\ref{prop: path-weight-1}}}k^3/\ell^3\big)^{|\cycleedges^\down\cup \cset^\down|}
\big(C_{\text{\tiny\ref{prop: path-weight-1}}}d^{\frac89}h\big)^{|(\cset^\up\cup\hdset^\up)\setminus\cycleedges^\up|}\\
&\cdot\maxnorm^{|\{t\in[2k]:\,t\notin\cset\cup\hdset\cup \cycleedges,\, t-1\in\hdset,\,H(t-1)\up H(t)\}|}
\|\stnd\|_1^{|\{t\in[2k]:\,t\notin\cset\cup\hdset\cup\cycleedges,\, t-1\notin\hdset,\,H(t-1)\up H(t)\}|},
\end{align*}
where the structures $\langle v,H,\cset,\hdset,\cycleedges,\badptclass\rangle$
must satisfy the relations established in Section~\ref{s: prop-data}.
Roughly, this is the strategy we are going to take, so that the estimation
could be reduced to combinatorial computations of which a major part is done in the next section and Section~\ref{s: few mult1}.
However, the summation over paths having many edges of multiplicity one requires another approach
and will be carried out later in Section~\ref{s: many mult1}.

\section{Summing over the admissible diagrams}

In Proposition~\ref{prop: path-weight-1}, we estimated the sum of path weights over paths sharing the same (part of)
data structure of the form $(v, \height, \cset, \hdset, \cycleedges,\badptclass)$.
The goal of this section is to further sum the path weights obtained in Proposition~\ref{prop: path-weight-1} over all possible realizations
of diagrams of a special type. 
To this aim, we will split a diagram into a sequence of sub-diagrams defined
in accordance with the structure of the sets $\cset$, $\hdset$ and $\cycleedges$,
and count the number of possible sub-diagrams separately. 
The domains of the sub-diagrams will be integer intervals,
with the constraint of not containing any of the points $\{t-1,t\}$ for $t\in\cset\cup \cycleedges$, nor any points $t\in \hdset$ for which
there is an up arrow from $t-1$ to $t$ in the sub-diagram.
One can keep in mind that these sub-diagrams comprise traveling on tree subgraphs and
on cycle edges of multiplicity $2$ since we omit all the times in $\cset\cup\cycleedges$ corresponding to travelling on other cycle edges.
The splitting will be made precise later in this section. 
As a first step, we restrict our attention to counting the number of possible choices for each sub-diagram viewed
as a new diagram with additional properties. Let us start by defining the class of diagrams.

Let $m\in \N$, $0\leq u\leq r\leq m$. 
We define $\Hset(m,r,u)$ as the collection of all couples $(\height,\w\hdset)$
where the function $\height:[0,m]\to \Z$ and the subset $\w\hdset\subseteq [0,m]$ 
satisfy the following properties:
\begin{enumerate}
\item\label{p1} $\height$ is a diagram starting at $(0,0)$  i.e. 
$\height(0)=0$ and $\height(t)-\height(t-1)\in \{\pm 1\}$  for any $t\in [m]$. 
\item\label{p6} We have $\vert \w\hdset\vert =r$ and $|\{t\in [m]:\; t-1\in \w\hdset \text{ and } \height(t-1)\up\height(t) \}|=u$.
\item\label{p2} For any $0\leq t'<t\leq m$ such that $\height(t)=\height(t')$ and $\height(t)\leq \height(T)$ for all $T\in [t',t]$, we have
 either $\{t',t\}\subseteq  \w\hdset$ or $\{t',t\}\cap \w\hdset=\emptyset$.
\item\label{p3} For any $t\in \w\hdset\setminus\{0\}$, there is a down arrow from $t-1$ to $t$ in $\height$.
\end{enumerate}

The set $\w\hdset$ will be taken to be a subset of the set $\hdset$ of the data structure.
Let us note that the above properties ``match'' properties of the data structure established previously.  
The fourth property above only asserts that there are no up arrows in the diagram leading to points in $\w\hdset$
while the second property implies that the number of up arrows from points in $\w\hdset$ is $u$.
The third property above will be automatically implied by Proposition~\ref{prop: above-B-point}. 
Let us record the following simple consequences of the above definition. 

\begin{claim}\label{claim0}
If the set $\Hset(m,r,u)$ is non-empty then necessarily $2u\leq m+1$.
\end{claim}
\begin{proof}
The claim immediately follows from the properties \ref{p6} and \ref{p3} above.
\end{proof}

\begin{claim}\label{claim1}
Let $(\height,\w\hdset)\in \Hset(m,r,u)$. Then, whenever $t<t'$
belong to $\w\hdset$ and $\height(t)=\height(t')$, we have $\height(\tau)\geq \height(t)$ for any $\tau\in [t,t']$.
\end{claim}
\begin{proof}
Assume the opposite, and let $\tau\in [t,t'-1]$ be the largest number such that $\height(\tau)<\height(t')$.
Then $\height(\tau+1)=\height(t')$ and $\height(x)\geq \height(t')$ for any $x\in [\tau+1,t']$. 
Whence, by property~\eqref{p2}, we have $\tau+1\in \w\hdset$ which contradicts property~\eqref{p3}.
\end{proof}

\begin{claim}\label{claim2}
Let $(\height,\w\hdset)\in \Hset(m,r,u)$
and let $t_1<t_2<\ldots<t_{r}$
be the elements of $\w\hdset$ arranged in increasing order.
Then $\height(t_1)\geq \height(t_2)\geq \ldots \geq \height(t_{r})$, and,
moreover, for any $i\in [r-1]$ and $x\in [t_i,t_{i+1}]$ we have
$\height(x)\geq \height(t_{i+1})$.
\end{claim}
\begin{proof}
It is sufficient to prove the second assertion only.
Assume that there exists $i<r$ such that 
$\height(t)<\height(t_{i+1})$ for some $t\in[t_i,t_{i+1}]$.
Let $\tau\in [t+1, t_{i+1}]$ be the largest number such that
$\height(\tau-1)<\height(t_{i+1})$. Then $\height(\tau)=\height(t_{i+1})$ and
$\height(x)\geq \height(t_{i+1})$ for any $x\in[\tau, t_{i+1}]$, whence by property~\eqref{p2}, we have $\tau\in \w\hdset$.
But $\height(\tau-1)<\height(\tau)$, contradicting property~\eqref{p3}.
\end{proof}

\begin{claim}\label{claim3}
Let  $Q:=(\height,\w\hdset)\in \Hset(m,r,u)$ and 
 $\hdset_Q':=\{t\in \w\hdset:\, H(t)\up H(t+1)\}$.
If $t<t'$ are two points in $\w\hdset$ with $t\in \hdset_Q'$ and 
$[t+1,t'-1]\cap\w\hdset=\emptyset$, then
necessarily $\height(t)=\height(t')$
and $\height(\tau)>\height(t)$ for all $t<\tau<t'$.
\end{claim}
\begin{proof}
Assume the opposite.
Then either there exists a point $x\in[t+1,t'-1]$ with $\height(x)=\height(t)$
and $\height(y)\geq \height(x)$ for all $y\in[t,x]$, or there exists a point $x'\in [t+1,t'-1]$ with $\height(x')=\height(t')$
and $\height(y')\geq \height(x')$ for all $y'\in[x',t']$. By property~\eqref{p2}, this would imply the existence of a point in $[t+1,t'-1]\cap \w\hdset$ contradicting the hypothesis. 
\end{proof}

Before proceeding with bounding the cardinality of $\Hset(m,r,u)$, we need the following auxiliary lemma. 
\begin{lemma}\label{lem: dyck}
Let $s,p\in\N$ and let $S$ be the set consisting of all sequences of non-empty Dyck paths $(H_i)_{i\leq s}$,
where each $H_i$ is a Dyck path with no returns to zero, except for the right end point,
and the total length of domains of the Dyck paths $(H_i)_{i\leq s}$ is equal to $2p$.
Then 
$$
|S|\leq\frac{s}{2p-s} {2p-s\choose p}.
$$
\end{lemma}
\begin{proof}
We will define a mapping 
$f$ from $S$ to the set of diagrams on $[0,2p]$,
where $D=f(\{H_i\}_{i\leq s})$ is obtained simply by concatenating $(H_i)_{i\leq s}$. 
Observe that $f$ is injective as the Dyck paths $H_i$ are non-zero in the interior of their respective domains;
hence, every sequence $(H_j)_{j\leq s}$ can be reconstructed from the diagram $D$
by splitting it into separate Dyck paths at the points where the diagram $D$ takes value zero.  

Further, the total number of diagrams in the range of $f$
is the number of Dyck paths on $[0,2p]$ having $s$ returns to zero, which is given by 
$$
\frac{s}{2p-s} {2p-s\choose p}.
$$
\end{proof}

The following short calculation will be used later in the section. 
\begin{lemma}\label{lem: calculation-binomial}
Let $L>1$. Then for any $p\in \N$, we have 
$$
\sum_{s=1}^p \frac{s}{2p-s} {2p-s\choose p} L^s \leq \frac{\max(L,2)^{2p-1}}{(\max(L,2)-1)^{p-1}}. 
$$
\end{lemma}
\begin{proof}

We start by setting $\w L:= \max(L,2)$ and  defining 
$$
\alpha(p):=\frac{1}{\w L^p} \sum_{s=1}^p \frac{s}{2p-s} {2p-s\choose p} \w L^s=\sum_{s=0}^{p-1} \frac{p-s}{p+s}{p+s \choose p}\w L^{-s},
$$
where the second equality follows by a change of variables. 
To estimate $\alpha(p)$, we will establish a recursive formula. To this aim, we calculate 
\begin{align*}
\alpha(p+1)-\alpha(p)&=\frac{1}{2p+1}{2p+1\choose p+1} \w L^{-p}+ \sum_{s=1}^{p-1} \Big[ \frac{p+1-s}{p+1+s}{p+1+s \choose p+1}-\frac{p-s}{p+s}{p+s \choose p}\Big] \w L^{-s}\\
&= \frac{1}{p+1}{2p\choose p} \w L^{-p}+ \sum_{s=1}^{p-1} \frac{s(p+2-s)}{(p+1)(p+s)} {p+s\choose p} \w L^{-s}\\
&= \sum_{s=1}^{p} \frac{s(p+2-s)}{(p+1)(p+s)} {p+s\choose p} \w L^{-s}\\
&= \sum_{s=1}^{p} \frac{(p+2-s)}{(p+s)} {p+s\choose p+1} \w L^{-s}=\w L^{-1}\sum_{s=0}^{p-1} \frac{(p+1-s)}{(p+1+s)} {p+1+s\choose p+1} \w L^{-s}\\
&=\w L^{-1}\, \alpha(p+1) - \frac{1}{2p+1}{2p+1\choose p+1} \w L^{-(p+1)}.
\end{align*}
Therefore, we deduce that $$\alpha(p+1)\leq \frac{\w L}{\w L-1} \alpha(p).$$
Since $\alpha(1)=1$, we deduce that 
$$\alpha(p)\leq \frac{\w L^{p-1}}{(\w L-1)^{p-1}} .$$ 
Replacing $\alpha(p)$ by its definition, we finish the proof. 
\end{proof}

We are now ready to estimate the cardinality of $\Hset(m,r,u)$.

\begin{prop}\label{prop: window}
Let $m\in \N$, $0\leq  u\leq r\leq m$. Then we have 
$$\vert \Hset(m,r,u)\vert\leq \Big(\frac{C_{\text{\tiny\ref{prop: window}}}\, m}{r-u+1}\Big)^{2(r-u+2)}\, \beta(u),$$
where $C_{\text{\tiny\ref{prop: window}}}$ is a universal constant and $\beta(u)$ is given by 
\[\beta(u)= \begin{cases} 
      2^m & \text{ if } u\leq 1,\\
      & \\
    \displaystyle\sum_{p=u-1}^{\lfloor m/2\rfloor}\frac{u-1}{2p-u+1} {2p-u+1\choose p} 2^{m-2p}& \text{ otherwise.} 
   \end{cases}
\]
\end{prop}
\begin{proof}
The first assertion follows since there are $2^m$ ways to form a diagram on $[m]$ and ${m\choose r}$ choices for the set $\w\hdset$. Therefore, we always have 
$$
\vert \Hset(m,r,u)\vert\leq 2^m {m\choose r}\leq 2^m \Big(\frac{em}{r}\Big)^r.
$$
When $u\leq 1$, we can replace $r$ by $r-u+1$ at the expense of changing the constant to deduce the first estimate of the lemma.  
In the remainder, we suppose $u>1$.

Fix for a moment an element $Q:=(\height,\w\hdset)\in \Hset(m,r,u)$,
and let $t_1<t_2<\ldots<t_{r}$
be the ordered elements of $\w\hdset$. Additionally, we set $t_0:=0$ and $t_{r+1}:=m$.
We define a collection $(H_i)_{i=0}^{r}$ of integer-valued functions as follows. 
Take $i\in [0,r]$ and let $H_i$ be a mapping from
$[0,t_{i+1}-t_i]$ to $\Z$ with $H_i(t):=\height(t+t_i)-\height(t_i)$, $t\in[0,t_{i+1}-t_i]$ i.e. we take the part of the diagram $\height$ lying between $t_i$ and $t_{i+1}$ and shift it so that it starts at $(0,0)$. 
Note that with this construction, the total length of the $H_i$'s, $i\in [0, r]$,  is equal to $m$. 
In view of Claim~\ref{claim3}, 
whenever $t_i\in \hdset_Q'$, $i=1,\ldots, r-1$, the function $H_i$ is non-negative and equal zero only at the  endpoints of its domain, thus $H_i$ is a Dyck path of positive length with no returns to zero except for the right-end point. 

Now, we define a mapping $f$ from $\Hset(m,r,u)$ to the sequences of the form $(H_i)_{i\in [0,r]}$ where for each $Q=(\height,\w\hdset)\in \Hset(m,r,u)$, the sequence $f(Q)=(H_i)_{i\in [0,r]}$ is constructed as above. To prove that $f$ is injective, observe that $\height$ can be uniquely reconstructed from the sequence $(H_i)_{i\in [0,r]}$ by simply concatenating the diagrams, and
that $\w\hdset$ can be uniquely reconstructed by taking the end points of $H_j$, $j\in  [1,r-1]$.

The injectivity of $f$ implies that $\vert \Hset(m,r,u)\vert$ is equal to the total number of sequences $(H_i)_{i\in [0,r]}$ in the range of $f$. 
It follows from a previous observation that in this sequence, there are at least $u-1$ non-empty Dyck paths with no returns to zero except for the right-end points. 

Fix a subset $I\subseteq [0,r]$ and suppose that $H_i$ is a non-empty Dyck path with no returns to zero (except for the right-end point) for every $i\in I$, and the total length of the $H_i$'s, $i\in I$, is equal to $2p$ with $\vert I\vert\leq p\leq \lfloor m/2\rfloor$. Then, by Lemma~\ref{lem: dyck}, the  number of admissible $(H_i)_{i\in I}$ of total length $2p$ is at most 
\begin{equation}\label{eq: dyck-path}
\frac{\vert I\vert}{2p-\vert I\vert} {2p-\vert I\vert\choose p}.
\end{equation}
Further, the number of admissible $(H_i)_{i\in I^c}$ of total length $m-2p$ is at most 
$$
2^{m-2p} {m-2p +\vert I^c\vert -1 \choose \vert I^c\vert-1}.
$$ 
Indeed, the above binomial coefficient corresponds to splitting an integer interval of length $m-2p$ into $\vert I^c\vert$
sub-intervals and $2^{m-2p}$ bounds the total number of choices of
$(H_i)_{i\in I^c}$ given their domains of total length $m-2p$.

Combining these estimates and using that there are at least $u-1$ non-empty Dyck paths with no returns to zero, we get 
$$
\vert \Hset(m,r,u)\vert\leq \sum_{\underset{ \vert I\vert\geq u-1}{I\subseteq [0,r]}} \sum_{p=\vert I\vert}^{\lfloor m/2\rfloor} 
\frac{\vert I\vert}{2p-\vert I\vert} {2p-\vert I\vert\choose p} 2^{m-2p}{m-2p +\vert I^c\vert-1 \choose \vert I^c\vert-1}.
$$
Now using that $\vert I\vert \geq u-1\geq 1$, that \eqref{eq: dyck-path} is decreasing in $\vert I\vert$, and interchanging the sums, we get  
\begin{align*}
\vert \Hset(m,r,u)\vert&\leq \sum_{p=u-1}^{\lfloor m/2\rfloor}\frac{u-1}{2p-u+1} {2p-u+1\choose p} 2^{m-2p} \sum_{\underset{ \vert I\vert\geq u-1}{I\subseteq [0,r]}} 
{m-2p +\vert I^c\vert \choose \vert I^c\vert-1}\\ &\\
&\leq\sum_{p=u-1}^{\lfloor m/2\rfloor}\frac{u-1}{2p-u+1} {2p-u+1\choose p} 2^{m-2p} \Big(\frac{Cm}{r-u+1}\Big)^{r-u+2} \sum_{s\leq r-u+2}
{m-2p +s \choose s-1},
\end{align*}
where we trivially bounded the number of choices of the set $I$ (or $I^c$) of a given admissible cardinality
by $\Big(\frac{Cm}{r-u+1}\Big)^{r-u+2}$, for an appropriate constant $C$. 
It remains to bound similarly the last binomial coefficient to finish the proof. 
\end{proof}

The last proposition will serve as
a tool for counting number of diagrams within a given family of data structures,
by splitting them into sub-diagrams belonging to some $\Hset(m, r, u)$. 
To this aim, let us first introduce the class of diagrams which will be of interest to us.
Let $k\in \N\cup\{0\}$, $v\in [n]$, and assume we
are given sets $\cset, \cycleedges\subset[2k]$, $\hdset^\up\subset [0,2k]$,
and each of the first two sets is partitioned into two subsets: $\cset=\cset^\up\cup\cset^\down$,
$\cycleedges=\cycleedges^\up\cup\cycleedges^\down$.
So far, we {\it do not assign} to the sets the same meaning as when we defined a data structure corresponding
to a path; here, these are just some abstract sets.
In what follows, these sets (and the numbers $k,v$) are assumed to be fixed.
The sets $\cset, \hdset^{\up}, \mathcal{C}$ define a {\it minimal} (i.e.\ of smallest possible cardinality) partition $(K_j)_{j\leq \eta}$
of $[0,2k]$ into integer subintervals where each subinterval $K_j$ satisfies 
$$
K_j\cap (\cset\cup \hdset^{\up}\cup \mathcal{C})\subseteq \{\min K_j\}.
$$
In what follows, we will view the number $\eta$ and the collection $(K_j)_{j\leq \eta}$
as functions of $k,\cset, \hdset^{\up}, \mathcal{C}$.

For parameters $U\leq R\leq 2k$, let $\w\Hset(2k, R, U)$ be  the collection of all couples  $(\height,\w\hdset)$
where the function $\height:[0,2k]\to \Z$ and the subset $\w\hdset\subseteq [0,2k]$ 
satisfy the following properties:

\begin{enumerate}
\item\label{p1'} $\height$ is a diagram starting at $(0,0)$  i.e. 
$\height(0)=0$ and $\height(t)-\height(t-1)\in \{\pm 1\}$  for any $t\in [2k]$. 
\item\label{p2'} For any $t\in  \cset^{\up}\cup \hdset^{\up}\cup \mathcal{C}^\up$ (resp. $t\in \cset^{\down}\cup \mathcal{C}^\down$), we have $\height(t-1)\up \height(t)$ (resp. $\height(t-1)\down\height(t)$). 
\item\label{p6'} We have $\vert \w\hdset\vert =R$ and $|\{t\in [2k]:\; t-1\in \w\hdset \text{ and } \height(t-1)\up\height(t) \}|=U$.
\item\label{p2''}
We require that for each $j\leq \eta$,
$(\height_{\mid K_j}, \w\hdset\cap K_j)$ belong (up to an
appropriate shifting of the coordinate system) to the set $\Hset(m_j, r_j, u_j)$
defined at the beginning of the section with $m_j:=\vert K_j\vert -1$ and 
$$
\vert \w\hdset\cap K_j\vert=r_j\, \text{ and } \,  |\{t\leq \max K_j:\; t-1\in \w\hdset\cap K_j \text{ and } \height(t-1)\up\height(t) \}|=u_j,
$$
so that $U-\eta\leq \sum_{j\leq \eta}u_j\leq U$.
\end{enumerate}

Let us first note that from the definition of the splitting above, we have $\eta\leq \vert \cset\cup \hdset^{\up}\cup \mathcal{C}\vert +1$. 
Formally, the collection shifted sub-diagrams $(\height_{\mid K_j}, \w\hdset\cap K_j)$, $j\leq \eta$,
does not contain full information about the entire diagram $H$ because of the ``gaps'' between $\max K_j$ and $\min K_{j+1}
=\max K_j+1$,
$j<\eta$. However, since $\min K_{j+1}\in \cset\cup \hdset^{\up}\cup \mathcal{C}$
(in view of minimality of the partition), the differences $H(\min K_{j+1})-H(\max K_j)$
can be reconstructed by checking if $\min K_{j+1}\in\cset^\up\cup \hdset^{\up}\cup \mathcal{C}^\up$
or $\min K_{j+1}\in\cset^\down\cup  \mathcal{C}^\down$.
This simple observation will be very important since it implies that
the cardinality of the set $\w\Hset(2k, R, U)$ can be bounded by estimating the number of admissible sequences of the sub-diagrams.
In turn, 
we will make use of Proposition~\ref{prop: window} where the number of possible sub-diagrams on any given interval $K_j$ is bounded.

\begin{prop}\label{prop: number-diag}
Let $k\geq 2$, $\Gamma=\{(R,U):\, 0\leq U\leq R\leq 2k \text{ and } R-U\leq 4(\vert \cset\cup \hdset^{\up}\cup \mathcal{C}\vert +1)\}$. Then for any real number $L\in (1,\infty)$, we have 
$$
\sum_{(R,U)\in \Gamma}L^U\, \vert \w\Hset(2k, R, U)\vert  \leq  k\Big[\frac{C_{\text{\tiny\ref{prop: number-diag}}}\,L\,  k\log k}{\vert \cset\cup \hdset^{\up}\cup \mathcal{C}\vert+1}\Big]^{ C_{\text{\tiny\ref{prop: number-diag}}}(\vert \cset\cup \hdset^{\up}\cup \mathcal{C}\vert+1)} \Big(\frac{\max(L,2)^2}{\max(L,2)-1}\Big)^{k-\frac{\vert \cset\cup \hdset^{\up}\cup \mathcal{C}\vert}{2}}, 
$$
where $C_{\text{\tiny\ref{prop: number-diag}}}>0$ is a universal constant. 
\end{prop}
\begin{proof}
Let $(R,U)\in \Gamma$. 
Clearly, we have 
$$
\vert \w\Hset(2k, R, U)\vert \leq \sum_{\underset{u_1,\ldots,u_\eta}{r_1,\ldots,r_\eta}}\, \prod_{j=1}^{\eta} \vert \Hset(m_j, r_j, u_j)\vert, 
$$
where $\sum_{j=1}^{\eta} m_j\leq 2k- \vert \cset\cup \hdset^{\up}\cup \mathcal{C}\vert+1$,
and the integer sequences $(u_j)$ and $(r_j)$ satisfy $\sum_{j=1}^\eta r_j=R$, $U-\eta\leq\sum_{j=1}^\eta u_j\leq U$,
$0\leq u_j\leq r_j\leq m_j$ and $u_j\leq (m_j+1)/2$ (see Claim~\ref{claim0}) for any $j\leq \eta$.

Fix any two admissible sequences $(u_j)$ and $(r_j)$.
Using Proposition~\ref{prop: window}, we get 
\begin{equation}\label{eq: prod-diagram}
 \prod_{j=1}^{\eta} \vert \Hset(m_j, r_j, u_j)\vert\leq \prod_{j=1}^\eta \beta_{m_j}(u_j)
\Big(\frac{C_{\ref{prop: window}}\, m_j}{r_j-u_j+1}\Big)^{2(r_j-u_j+2)}.
\end{equation}
To bound the product of the second terms, let $\Phi$ be the set of all indices $j\leq \eta$ such that 
$$
\frac{m_j}{r_j-u_j+2}\leq \frac{ k\log k}{R-U+\eta+1}.
$$
Then 
$$
\prod_{j\in \Phi}\Big(\frac{C_{\ref{prop: window}}\,m_j}{r_j-u_j+1}\Big)^{2(r_j-u_j+2)}\leq \Big(\frac{2C_{\ref{prop: window}}\,  k \log k}{R-U+\eta+1}\Big)^{2(R-U+3\eta)}
$$
where we have used that $\vert \Phi\vert \leq \eta$ and $\sum_{j\in \Phi} (r_j-u_j)\leq R-U+\eta$. On the other hand, when $j\not\in \Phi$, we have 
$$
r_j-u_j+2\leq  \frac{R-U+\eta+1}{ k\log k}\, m_j,
$$ which together with the bound $\frac{C_{\ref{prop: window}}\, m_j}{r_j-u_j+1}\leq 2C_{\ref{prop: window}}\, k$ imply 
$$
\prod_{j\not\in \Phi}\Big(\frac{C_{\ref{prop: window}}\,m_j}{r_j-u_j+1}\Big)^{2(r_j-u_j+2)}\leq (2C_{\ref{prop: window}}\, k)^{\frac{4(R-U+\eta+1)}{ \log k}},
$$
where we have used that $\sum_{j\not\in \Phi} m_j\leq 2k$. 
Putting together the above estimates and \eqref{eq: prod-diagram},
and using that $\max(R-U,\eta)\leq 4(\vert \cset\cup \hdset^{\up}\cup \mathcal{C}\vert +1)$,
we get 
$$
\prod_{j=1}^{\eta} \vert \Hset(m_j, r_j, u_j)\vert  \leq \Big(\frac{C\, k\log k}
{\vert \cset\cup \hdset^{\up}\cup \mathcal{C}\vert +1}\Big)^{C(\vert \cset\cup \hdset^{\up}\cup \mathcal{C}\vert +1)}\,  \prod_{j=1}^\eta \beta_{m_j}(u_j).
$$
for some universal constant $C>0$. 
Next, note that the number of admissible choices $\alpha$ of the sequences $(r_j)_{j\leq \eta}$ can be estimated by 
$$
\alpha\leq {\eta+R-1 \choose \eta-1}. 
$$
Bounding the above binomial coefficient, we get
$$
\alpha\leq \bigg(C'\Big(1+\frac{R}{\eta}\Big)\bigg)^{\eta},
$$
for some universal constant $C'>0$. Further, a short calculation shows that 
$$
\bigg(C'\Big(1+\frac{R}{\eta}\Big)\bigg)^{\eta}\leq \Big(\frac{C''\, k\log k}{R-U+\eta+1}\Big)^{R-U+\eta},
$$
for some appropriate constant $C''>0$.
Putting together the above estimates, and using that $\max(R-U,\eta)\leq 4(\vert \cset\cup \hdset^{\up}\cup \mathcal{C}\vert +1)$, we can write 
$$
\sum_{(R,U)\in \Gamma}L^U\, \vert \w\Hset(2k, R, U)\vert  \leq \Big(\frac{\w C\, k\log k}{\vert \cset\cup \hdset^{\up}\cup \mathcal{C}\vert +1}\Big)^{\w C(\vert \cset\cup \hdset^{\up}\cup \mathcal{C}\vert +1)}\, \sum_{U=0}^{2k}\,  \sum_{u_1,\ldots,u_\eta}  \, \prod_{j=1}^\eta \beta_{m_j}(u_j)\, L^{u_j}
$$
for some appropriate constant $\w C>0$, where the sequences $(u_j)_{j=1}^\eta$
must satisfy $2u_j\leq m_j+1$, $j\leq \eta$.
Interchanging the sum and the product, we get
$$
\sum_{(R,U)\in \Gamma}L^U \vert \w\Hset(2k, R, U)\vert  \leq \Big(\frac{\w C\, k\log k}{\vert \cset\cup \hdset^{\up}\cup \mathcal{C}\vert +1}\Big)^{\w C(\vert \cset\cup \hdset^{\up}\cup \mathcal{C}\vert +1)}
(2k+1)\prod\limits_{j=1}^\eta \bigg(\sum\limits_{u=0}^{\lfloor(m_j+1)/2\rfloor}\beta_{m_j}(u)L^u\bigg).
$$ 
Hence,
\begin{align*}
\sum_{(R,U)\in \Gamma}L^U\, \vert \w\Hset(2k, R, U)\vert & \leq \Big(\frac{\w C\, k\log k}{\vert \cset\cup \hdset^{\up}\cup \mathcal{C}\vert +1}\Big)^{\w C(\vert \cset\cup \hdset^{\up}\cup \mathcal{C}\vert +1)}\, (2k+1)\, \\ 
&\, \cdot\prod_{j=1}^\eta \Big( 2^{m_j}+L2^{m_j} +\sum_{p=1}^{\lfloor m_j/2\rfloor} \sum_{u=2}^{p+1} \frac{u-1}{2p-u+1} {2p-u+1\choose p} 2^{m_j-2p} \, L^{u}\Big)\\
&\leq \Big(\frac{\w C\, k\log k}{\vert \cset\cup \hdset^{\up}\cup \mathcal{C}\vert +1}\Big)^{\w C(\vert \cset\cup \hdset^{\up}\cup \mathcal{C}\vert +1)}\, (2k+1)\, \\ 
&\, \cdot\prod_{j=1}^\eta \Big( 2^{m_j}+L2^{m_j} +L\sum_{p=1}^{\lfloor m_j/2\rfloor} \sum_{u=1}^{p} \frac{u}{2p-u} {2p-u\choose p} 2^{m_j-2p} \, L^{u}\Big)\\
\end{align*}
where we have interchanged the sums over $u$ and $p$ from the definition of $\beta$ in the first inequality, and made a change of variables in the second.  
Using Lemma~\ref{lem: calculation-binomial} and the estimate $\sum_{j=1}^{\eta} m_j\leq
2k- \vert \cset\cup \hdset^{\up}\cup \mathcal{C}\vert+1$, we get
\begin{align*}
\sum_{(R,U)\in \Gamma}L^U\, \vert \w\Hset(2k, R, U)\vert & \leq k\, 2^{2k- \vert \cset\cup \hdset^{\up}\cup \mathcal{C}\vert}\,  \Big(\frac{\bar C\, k\log k}{\vert \cset\cup \hdset^{\up}\cup \mathcal{C}\vert +1}\Big)^{\bar C(\vert \cset\cup \hdset^{\up}\cup \mathcal{C}\vert +1)}\, \\ 
&\, \cdot\prod_{j=1}^\eta \Big[ 1+L+L\sum_{p=1}^{\lfloor m_j/2\rfloor} \frac{\max(L,2)^{2p-1}}{4^p(\max(L,2)-1)^{p-1}}\Big],\\
\end{align*}
for some appropriate constant $\bar C>0$. Finally, since $\max(L,2)^2\geq 4(\max(L,2)-1)$, a short calculation implies that 
\begin{align*}
\sum_{(R,U)\in \Gamma}L^U\, \vert \w\Hset(2k, R, U)\vert & \leq k\, 2^{2k- \vert \cset\cup \hdset^{\up}\cup \mathcal{C}\vert}\,  \Big(\frac{\hat C\, L\, k\log k}{\vert \cset\cup \hdset^{\up}\cup \mathcal{C}\vert +1}\Big)^{\hat C(\vert \cset\cup \hdset^{\up}\cup \mathcal{C}\vert +1)}\, \\ 
&\, \cdot\prod_{j=1}^\eta \Big[ (1+m_j)\Big(\frac{\max(L,2)^2}{4(\max(L,2)-1)}\Big)^{\lfloor m_j/2\rfloor}\Big].\\
\end{align*}
Using that $\sum_{j=1}^{\eta} m_j\leq 2k- \vert \cset\cup \hdset^{\up}\cup \mathcal{C}\vert+1$, and the Arithmetic-Geometric mean inequality with $\sum_{j=1}^{\eta} (1+m_j)\leq 4k$ and 
$\eta\leq \vert \cset\cup \hdset^{\up}\cup \mathcal{C}\vert+1$, we finish the proof. 
\end{proof}

We now have everything in place in order to ``integrate'' the estimate of Proposition~\ref{prop: path-weight-1}
over admissible choices of diagrams and of the mappings $\badptclass$. Below is the 
main result of this section. 

\begin{prop}\label{prop: contribution-diagram}
Let parameters $n,\ell,d, \dmax,k\in \N$, $h,\maxnorm\in \R_+$ and $\stnd\in \R_+^n$
satisfy relations \eqref{eq: condition-d-dmax}, with $k^3/\ell^3\leq d^{\frac89}$,
and let $M\in \matrixset(n,\ell, d,\dmax, h,\maxnorm, \stnd)$.  
Fix a vertex $v\in [n]$, subsets $\cset,\hdset^{\up},\cycleedges$, partitions $(\cset^{\up}, \cset^{\down})$ and  $(\cycleedges^{\up}, \cycleedges^{\down})$ of $\cset$ and $\mathcal{C}$ respectively. 
Denote by $\textbf{P}'$ the collection of all paths $\path$ of length $2k$ on $G_M$
whose corresponding data structures satisfy $v_\path=v$,
$\cset^\up_\path=\cset^\up$, $\cset^\down_\path=\cset^\down$, $\hdset^\up_\path=\hdset^\up$,
$\cycleedges_\path^\up=\cycleedges^\up$, $\cycleedges_\path^\down=\cycleedges^\down$, 
and the other elements of the data structure --- $H_\path$, $\hdset^\down_\path$, $\weight_\path,\badptclass_\path$
--- take any admissible values (i.e.\ compatible with the fixed part of the data structure). 
Then
\begin{align*}
\sum\limits_{\path\in \textbf{P}'}\pathw_M(\path)\leq
&\, 
\|\stnd\|_1^k\Big(\frac{\max(L,2)^2}{\max(L,2)-1}\Big)^{k}
 \Big(\frac{h}{\|\stnd\|_1}\Big)^{|\mathcal C|/2}\, (2\dmax)^{|\mathcal C^\up|}\,
\big(C_{\text{\tiny\ref{prop: contribution-diagram}}}k^3/\ell^3\big)^{|\mathcal C^\down\cup \cset^\down|}\\
&\cdot 
\Big(\frac{C_{\text{\tiny\ref{prop: contribution-diagram}}}d^{\frac89}h}{\|\stnd\|_1}\Big)^{|(\cset^\up\cup\hdset^\up)\setminus\mathcal C^\up|}
\Big(\frac{C_{\text{\tiny\ref{prop: contribution-diagram}}}\, L\, k\, \vert \net_{mjr}\vert\, \log k}
{\vert \cset\cup \hdset^{\up}\cup \mathcal{C}\vert +1} \Big)^{C_{\text{\tiny\ref{prop: contribution-diagram}}}
(\vert \cset\cup \hdset^{\up}\cup \mathcal{C}\vert+1)},
\end{align*}
where $C_{\text{\tiny\ref{prop: contribution-diagram}}}>$ is a universal constant and $L:=  \maxnorm/\|\stnd\|_1$. 
\end{prop}
\begin{proof}
Let us first fix a diagram $H$ and a set $\hdset$ compatible with the known part of the data structure,
and consider the collection $\w P(H,\hdset)$ of all closed paths $\path$
of length $2k$ on $G_M$ with corresponding data structures of the form 
$\langle v,H, \cset, \hdset,\cycleedges, \weight_\path,\badptclass_\path\rangle$,
with arbitrary weight functions $\weight_\path$ and mappings $\badptclass_\path$.
Note that by Proposition~\ref{prop: properties-V} and the definition of the mapping
$\badptclass$, the number of possible choices for $\badptclass_\path$ is bounded above by 
$$
\vert \net_{mjr}\vert^{6(\vert\cset\cup \hdset^{\up}\cup \cycleedges\vert+1)}.
$$
Therefore, applying Proposition~\ref{prop: path-weight-1}, we get 
\begin{align*}
\sum\limits_{\path\in \w P(H,\hdset)}\pathw_M(\path)&\leq
h^{|\cycleedges|/2}(2\dmax)^{|\cycleedges^\up|}\,
\big(C_{\text{\tiny\ref{prop: path-weight-1}}}k^3/\ell^3\big)^{|\cycleedges^\down\cup \cset^\down|}
\big(C_{\text{\tiny\ref{prop: path-weight-1}}}d^{\frac89}h\big)^{|(\cset^\up\cup\hdset^\up)\setminus\cycleedges^\up|}
\vert \net_{mjr}\vert^{6(\vert\cset\cup \hdset^{\up}\cup \cycleedges\vert+1)}\\
&\cdot\maxnorm^{|\{t\in[2k]:\,t\notin\cset\cup\hdset\cup \cycleedges,\, t-1\in\hdset,\,H(t-1)\up H(t)\}|}
\|\stnd\|_1^{|\{t\in[2k]:\,t\notin\cset\cup\hdset\cup\cycleedges,\, t-1\notin\hdset,\,H(t-1)\up H(t)\}|}\\
&\leq h^{|\cycleedges|/2}(2\dmax)^{|\cycleedges^\up|}\,
\big(C_{\text{\tiny\ref{prop: path-weight-1}}}k^3/\ell^3\big)^{|\cycleedges^\down\cup \cset^\down|}
\big(C_{\text{\tiny\ref{prop: path-weight-1}}}d^{\frac89}h\big)^{|(\cset^\up\cup\hdset^\up)\setminus\cycleedges^\up|}
\vert \net_{mjr}\vert^{6(\vert\cset\cup \hdset^{\up}\cup \cycleedges\vert+1)}\\
&\cdot\big(\maxnorm/\|\stnd\|_1\big)^{|\{t\in[2k]:\,\, t-1\in\hdset,\,H(t-1)\up H(t)\}|}
\|\stnd\|_1^{|\{t\in[2k]:\,t\notin\cset\cup\hdset\cup\cycleedges,\,H(t-1)\up H(t)\}|}.
\end{align*}
The sets $\cset, \hdset^{\up}, \mathcal{C}$ define a {\it minimal} partition
$(K_j)_{j\leq \eta}$ of $[0,2k]$ into integer subintervals where each $K_j$ satisfies 
$$
K_j\cap (\cset\cup \hdset^{\up}\cup \mathcal{C})\subset\{\min K_j\}, \;\;j\leq \eta.
$$
Note that for any $j\leq \eta$, the part of the diagram restricted to $K_j$ satisfies
\begin{itemize}
\item Whenever $t\leq t'$ belong to $K_j$ and $H(T)\geq H(t)=H(t')$ for all $T\in[t,t']$,
we have either $\{t,t'\}\subset \hdset$ or $\{t,t'\}\cap \hdset=\emptyset$.
This property follows from the condition $[t+1,t']\cap (\cset\cup\cycleedges)=\emptyset$ and Proposition~\ref{prop: above-B-point}.
\item For any $t\in K_j\cap \hdset\setminus\{\min K_j\}$, we have $H(t-1)\down H(t)$.
\end{itemize}
Denote $R:=\vert \hdset\vert$ and
$U:=\vert\{t\in [2k]:\, t-1\in \hdset \text{ and } H(t-1)\up H(t)\}\vert$.
Therefore, applying the definition of $\w\Hset(2k, R, U)$ and the above observations,
we get that the pair $(H,\hdset)$ belongs to the collection $\w\Hset(2k, R, U)$.
Moreover, it follows from Proposition~\ref{prop: nb-levels} that 
$R-U\leq 4\vert \hdset^{\up}\cup \mathcal{C}^\down\vert +1\leq 4(\vert \cset\cup \hdset^{\up}\cup \mathcal{C}\vert +1)$.

Thus, any pair $(H,\hdset)$ compatible with the fixed part of the data structure, belongs to
$\bigcup_{(R,U)\in\Gamma}\w\Hset(2k, R, U)$, where $\Gamma$ is defined in Proposition~\ref{prop: number-diag}.
Further, 
$$
\vert\{ t\in [2k]\setminus\cset\cup\hdset\cup\mathcal{C}:\, H(t-1)\up H(t)\}\vert
\leq  k-\frac{1}{2}\vert \mathcal{C}\vert- \vert (\cset^{\up}\cup \hdset^{\up})\setminus \mathcal{C}^{\up}\vert,
$$
where we have used that the total number of up-arrows in the diagram is
$k+\frac{\vert\mathcal{C}^{\up}\vert-\vert\mathcal{C}^{\down}\vert}{2}$ (see Lemma~\ref{lem: up-down-total}). Therefore, using the above, we can write 
\begin{align*}
\sum\limits_{\path\in \textbf{P}'}\pathw_M(\path)\leq
&h^{|\cycleedges|/2}(2\dmax)^{|\cycleedges^\up|}\,
\big(C_{\text{\tiny\ref{prop: path-weight-1}}}k^3/\ell^3\big)^{|\cycleedges^\down\cup \cset^\down|}
\big(C_{\text{\tiny\ref{prop: path-weight-1}}}d^{\frac89}h\big)^{|(\cset^\up\cup\hdset^\up)\setminus\cycleedges^\up|}
\vert \net_{mjr}\vert^{6(\vert\cset\cup \hdset^{\up}\cup \cycleedges\vert+1)}\\
&\cdot\sum_{U} \sum_{0\leq R-U\leq 4\vert \hdset^{\up}\cup \mathcal{C}^\down\vert +1} \big(\maxnorm/\|\stnd\|_1\big)^{U}\cdot 
\|\stnd\|_1^{ k-\frac{1}{2}\vert \mathcal{C}\vert- \vert (\cset^{\up}\cup \hdset^{\up})\setminus \mathcal{C}^{\up}\vert
} \cdot \vert \w\Hset(2k, R, U)\vert.
\end{align*}
Applying Proposition~\ref{prop: number-diag} with $L:= \maxnorm/\|\stnd\|_1$, we get 
\begin{align*}
\sum\limits_{\path\in \textbf{P}'}\pathw_M(\path)\leq
&\,  \Big(\frac{h}{\|\stnd\|_1}\Big)^{|\mathcal C|/2}\, (2\dmax)^{|\mathcal C^\up|}\,
\big(C_{\text{\tiny\ref{prop: path-weight-1}}}k^3/\ell^3\big)^{|\mathcal C^\down\cup \cset^\down|}
\Big(\frac{C_{\text{\tiny\ref{prop: path-weight-1}}}d^{\frac89}h}{\|\stnd\|_1}\Big)^{|(\cset^\up\cup\hdset^\up)\setminus\mathcal C^\up|}\\
&\cdot \|\stnd\|_1^k\Big(\frac{C\, L\, k\, \vert \net_{mjr}\vert\, \log k}
{\vert \cset\cup \hdset^{\up}\cup \mathcal{C}\vert +1} \Big)^{C(\vert \cset\cup \hdset^{\up}\cup \mathcal{C}\vert+1)} 
\Big(\frac{\max(L,2)^2}{\max(L,2)-1}\Big)^{k-\frac{\vert \cset\cup \hdset^{\up}\cup \mathcal{C}\vert}{2}},
\end{align*}
for some appropriate constant $C>0$. The proof is finished after rearranging the terms in the above expression. 
\end{proof}

\section{Contribution of paths with few multiplicity one edges}\label{s: few mult1}

In the previous two sections, we bounded
the sum of path weights over all paths sharing the same starting vertex and sets
$\cset^\up$, $\cset^\down$, $\hdset^{\up}$, $\cycleedges^\up$, $\cycleedges^\down$.
In other words, we have integrated the expression in \eqref{eq: def-path-weight} over all admissible choices of the weight function $\weight$,
the diagram, the set $\hdset^\down$ and the mapping $\badptclass$.
We now count the contribution of the remaining quantities forming the data structure.
We will be able to do so in terms of the number of multiplicity one edges.

Given a complete graph $K_{[n]}$ on $n$ vertices and a closed path $\path$ on $K_{[n]}$ of length $2k$,
let $\multone(\path)$ be the cardinality of a largest subset $S$ of edges of $K_{[n]}$ 
such that each edge from $S$ is travelled by $\path$ exactly once, and no two edges from $S$ are incident. 
We have the following observation. 

\begin{claim}
Let $\path$ be a path of length $2k$ on an $\ell$-tangle graph $G$ (with $5\leq\ell\leq k$) and denote
$$\cycleone:=\big\{t\in [2k]:\, \path(t-1)\edg \path(t) \text{ is a cycle edge of multiplicity one in } G_{\path}\big\}.$$
Then 
$$
\vert \cycleone\vert \leq 64\multone(\path)\, \Big(1+\frac{k}{\ell}\Big).
$$
\end{claim}
\begin{proof}
Let $S$ be a subset of edges of $G$ of cardinality $\multone(\path)$ such that each edge from $S$ is travelled by $\path$ exactly once, and no two edges from $S$ are incident. 
Note that any cycle edge in $G_\path$ travelled once by $\path$ is either in $S$ or incident to $S$, otherwise it would contradict the definition of $\multone(\path)$.  Since $G_{\path}$ is $\ell$-tangle free, then by Corollary~\ref{cor: nb-cycles meeting} any edge in $S$ has at most $64\big(1+\frac{k}{\ell}\big)$ cycle edges incident to it.  The claim follows. 
\end{proof}

The above claim provides information
on the cardinality of $\cycleedges^\up_\path$ in terms of the parameter $\multone(\path)$. We record this in the next statement. 

\begin{claim}\label{claim: cycle-up}
Let $\path$ be a path of length $2k$ on an $\ell$-tangle graph $G$ ($5\leq \ell\leq k$)
and let $H_\path$ and $\cycleedges_\path$ be the corresponding elements from the data structure for $\path$. Then 
$$
\vert \cycleedges_\path^{\up}\vert \leq \frac{\vert \cycleedges_\path\vert}{3} +64 \multone(\path)\, \Big(1+\frac{k}{\ell}\Big)
$$
\end{claim}
\begin{proof}
Let $\cycleone$ be defined as in the above claim. Clearly, we have 
$$
\vert \cycleedges^{\up}_\path\vert =\vert \cycleone\vert +\vert \cycleedges^{\up}_\path\setminus \cycleone\vert. 
$$
Now note that for any $t\in \cycleedges_\path\setminus \cycleone$, $\path(t-1)\edg\path(t)$ is a cycle edge of multiplicity at least $3$. Since for such an edge, we associated one up arrow and at least $2$ down arrows, then we deduce that 
$$
\vert \cycleedges^{\up}_\path\setminus \cycleone\vert\leq \frac{\vert \cycleedges_\path\setminus \cycleone\vert}{3}.
$$
Putting together the above estimates and using the previous claim, we finish the proof. 
\end{proof}

\begin{prop}
Let $n,d$ be large natural numbers, and let parameters $n, d,\dmax\in \N$, $h,\maxnorm\in \R_+$ and $\stnd\in \R_+^n$
satisfy \eqref{eq: condition-d-dmax}. Suppose further that  $\multon, k\in\N$ and
\begin{equation}\label{eq: conditions}
h\, d^{0.9} \leq \Vert \stnd\Vert_1<\maxnorm\leq 2h\dmax\quad \text{ and } \quad d\leq k^2,
\end{equation}
and denote $L= \frac{\maxnorm}{\| \stnd\|_1}$. 
Then for any $M\in \matrixset(n,k/\log^2 d, d, \dmax, h,\maxnorm, \stnd)$, we have 
$$
\sum\limits_{\path:\, \multone(\path)=\multon }\pathw_M(\path)\leq n\,
\Big(\frac{\|\stnd\|_1\max(L,2)^2}{\max(L,2)-1}\Big)^{k}\, e^{C  \max(\multon, 1)\log^{10} k}   \cdot 
\exp\Big( \frac{C k L \vert \net_{mjr}\vert \log k}{\sqrt[C]{d}}
\Big),
$$
where the sum is over closed paths of length $2k$ on $K_{[n]}$, and $C>1$ is a universal constant. 
\end{prop}

\begin{proof}
Let $\Gamma$ be the collection of $5$--tuples
of the form
$(\cset^{\up}, \cset^{\down}, \hdset^{\up}, \mathcal{C}^{\up}, \mathcal{C}^{\down})$,
where $\cset^{\up}, \cset^{\down}, \mathcal{C}^{\up}, \mathcal{C}^{\down}$ are subsets of $[2k]$,
$\hdset^{\up}$ is a subset of $[0,2k]$; the sets $\cset^{\up}, \cset^{\down}$ are disjoint
(and similarly for $\mathcal{C}^{\up}, \mathcal{C}^{\down}$) and, additionally,
$$
\vert \cycleedges^{\up}\setminus \cycleone\vert\leq \frac{\vert \cycleedges\vert}{3} +128\multon\log^2d
\text{ \ and \ } \vert \cset^{\up}\setminus \cycleedges\vert \geq \vert \cset^{\down}\setminus \cycleedges\vert - C\log^8d,
$$
for a large universal constant $C>0$. Note that by Claim~\ref{claim: cycle-up} and Proposition~\ref{prop: up minus down} applied with $\ell:= k/\log^2d$, $\Gamma$ contains all admissible realizations
of $\cset^{\up}_\path, \cset^{\down}_\path, \hdset^{\up}_\path, \mathcal{C}^{\up}_\path, \mathcal{C}^{\down}_\path$
for paths $\path$ with $\multone(\path)=\multon$.
Therefore, it follows from Proposition~\ref{prop: contribution-diagram} applied with $\ell=k/\log^2d$ that 
\begin{align*}
 \sum_{\path:\,  \multone(\path)=\multon }&\pathw_M(\path)\leq
 n\, \Big(\frac{\|\stnd\|_1\max(L,2)^2}{\max(L,2)-1}\Big)^{k}\cdot  \sum_{(\cset^{\up}, \cset^{\down}, \hdset^{\up}, \mathcal{C}^{\up}, \mathcal{C}^{\down})\in
\Gamma}\Bigg[ \Big(\frac{h}{\|\stnd\|_1}\Big)^{|\mathcal C|/2} \, (2\dmax)^{|\mathcal C^\up|}\\
&\cdot \big(C'\log^6d\big)^{|\mathcal C^\down\cup \cset^\down|}\Big(\frac{C'd^{\frac89} h}{\|\stnd\|_1}\Big)^{|(\cset^\up\cup\hdset^\up)\setminus\mathcal C^\up|}\cdot  \Big(\frac{C'\, L\,  k\, \vert \net_{mjr}\vert\log k}{\vert \cset\cup \hdset^{\up}\cup \mathcal{C}\vert+1}\Big)^{C'(\vert \cset\cup \hdset^{\up}\cup \mathcal{C}\vert+1)}\Bigg],
\end{align*}
where $C'>0$ a universal constant. Since $\vert \cycleedges^{\up}\vert \leq  \frac{\vert \cycleedges\vert}{3} +128\multon\log^2d$, we get
\begin{align*}
 \sum_{\path:\,  \multone(\path)= \multon }\pathw_M(\path)\leq
& n\, (2\dmax)^{128  \multon \log^2d} \Big(\frac{\|\stnd\|_1\, \max(L,2)^2}{\max(L,2)-1}\Big)^{k}\\
&\cdot\sum_{(\cset^{\up}, \cset^{\down}, \hdset^{\up}, \mathcal{C}^{\up}, \mathcal{C}^{\down})\in\Gamma}
\Bigg[ \Big(\frac{C'\sqrt{h}\sqrt[3]{ 2\dmax}\, \log^6d}{\sqrt{\|\stnd\|_1}}\Big)^{|\mathcal C|}
\cdot \Big(C'\log^6 d\Big)^{\vert \cset^{\down}\setminus \mathcal{C}\vert}
\\
&\cdot \Big(\frac{C'd^{\frac89} h}{\|\stnd\|_1}\Big)^{|(\cset^\up\cup\hdset^\up)\setminus\mathcal C^\up|}\cdot  \Big(\frac{C'\, L\,  k\, \vert \net_{mjr}\vert\log k}{\vert \cset\cup \hdset^{\up}\cup \mathcal{C}\vert+1}\Big)^{C'(\vert \cset\cup \hdset^{\up}\cup \mathcal{C}\vert+1)}\Bigg].
\end{align*}
 Now using that $\vert \cset^{\up}\setminus \cycleedges\vert \geq \frac12 \vert \cset\setminus\cycleedges\vert -\frac{C}{2}\log^8 d$, $hd^{\frac89}\leq \|\stnd\|_1$ and $\dmax\leq d^{\frac43}$, we get from the above that 
\begin{align*}
 \sum_{\path:\,  \multone(\path)= \multon }\pathw_M(\path)\leq
& n\, e^{\w C  \multon \log^3d+ \w C \log^9 d}
\Big(\frac{\|\stnd\|_1\max(L,2)^2}{\max(L,2)-1}\Big)^{k}\cdot \sum_{\Gamma}  \Bigg[
\Big(\frac{\sqrt{h}\sqrt[3]{\dmax}\, \log^6d}{\sqrt{\|\stnd\|_1}}\Big)^{|\mathcal C|} 
\\
&\cdot  \Big(\frac{\sqrt{d^{\frac89} h}\, \log^6d}{\sqrt{\|\stnd\|_1}}\Big)^{|\cset\cup \hdset^{\up}\setminus\mathcal C|}\cdot
\Big(\frac{\w C\, L\,  k\, \vert \net_{mjr}\vert\log k}{\vert \cset\cup \hdset^{\up}\cup \mathcal{C}\vert+1}\Big)^{\w C(\vert \cset\cup \hdset^{\up}\cup \mathcal{C}\vert+1)}\Bigg].
\end{align*}
Regrouping the terms and using that $\dmax\leq d^{\frac43}$, we deduce 
\begin{align*}
 \sum_{\path:\,  \multone(\path)= \multon }\pathw_M(\path)&\leq
n\, e^{\w C \multon \log^3 d+ \w C\log^9d}\Big(\frac{\|\stnd\|_1\max(L,2)^2}{\max(L,2)-1}\Big)^{k} \\
&\quad \cdot \sum_{\Gamma}   \Big(\frac{\sqrt{d^{\frac89}h}\, \log^6 d}{\sqrt{\|\stnd\|_1}}\Big)^{|\cset\cup \hdset^{\up}\cup\mathcal C|}\cdot  \Big(\frac{\w C\, L\,  k\, \vert \net_{mjr}\vert\log k}{\vert \cset\cup \hdset^{\up}\cup \mathcal{C}\vert+1}\Big)^{\w C(\vert \cset\cup \hdset^{\up}\cup \mathcal{C}\vert+1)}.
\end{align*}
Now note that given any non-negative number $s\leq 2k$,
there are less than $\big(\frac{2ek}{s}\big)^{5s}$ possible realizations
of the $5$--tuple $(\cset^{\up}, \cset^{\down}, \hdset^{\up}, \mathcal{C}^{\up}, \mathcal{C}^{\down})$
so that $|\cset\cup \hdset^{\up}\cup \mathcal{C}|=s$.
Therefore, we obtain
\begin{align*}
 \sum_{\path:\,  \multone(\path)= \multon }\pathw_M(\path)&\leq
n\, e^{\w C \multon \log^3 d+ \w C\log^9d} \Big(\frac{\|\stnd\|_1\max(L,2)^2}{\max(L,2)-1}\Big)^{k}\\ 
&\cdot   \sum_{s\leq 2k}\Big(\frac{\sqrt{d^{\frac89} h}\, \log^2d}{\sqrt{\|\stnd\|_1}}\Big)^{s} \Big(\frac{\hat C \, L\, k\, \vert \net_{mjr}\vert\log k }{s+1} \Big)^{\hat C(s+1)},
\end{align*}
for some universal constant $\hat C>0$. 
Using \eqref{eq: conditions}, a short calculation finishes the proof. 
\end{proof}

The above proposition provides a satisfactory bound
for paths $\path$ with sufficiently small $\multone(\path)$.
We record the following corollary which follows by an easy calculation using the estimate \eqref{eq: aux 5239565039853-50} on $\vert \net_{mjr}\vert$. 

\begin{cor}\label{cor: contribution-few multone}
There exist two universal constants $c_{\text{\tiny\ref{cor: contribution-few multone}}}<1$ and $C_{\text{\tiny\ref{cor: contribution-few multone}}}>1$ such that the following holds. 
Let $n$ be a large natural number, and let
$n,d, \dmax , k\in \N$, $h,\maxnorm\in \R_+$ and $\stnd\in \R_+^n$ satisfy \eqref{eq: condition-d-dmax}. Suppose that 
$$
h\, d^{0.9}\leq \Vert \stnd\Vert_1<\maxnorm\leq 2h\dmax,\quad \sqrt{\log n} \leq \sqrt{k}\leq d\leq k^2,\quad  \text{ and } \quad  L\leq d^{c_{\text{\tiny\ref{cor: contribution-few multone}}}},
$$
where $L= \frac{\maxnorm}{\Vert \stnd\Vert_1}$. 
Then for any $M\in \matrixset(n,k/\log^2 d,\dmax, \dmax, h,\maxnorm, \stnd)$, we have 
$$
\sum\limits_{\path:\, \multone(\path)<  \frac{k}{\log^{11} k} }\pathw_M(\path)\leq n\, \Big(\frac{\|\stnd\|_1\max(L,2)^2}{\max(L,2)-1}\Big)^{k}\, e^{C_{\text{\tiny\ref{cor: contribution-few multone}}}  k/\log k},
$$
where the sum is over closed paths of length $2k$ on $K_{[n]}$. 
\end{cor}


\section{Auxiliary probabilistic constructions}

Unlike the first part of the paper where everything was carried in a deterministic setting, the remainder of the paper will heavily rely on probabilistic facts. We gather in this section some of the tools we will be using. 


The following is the well known Bernstein inequality (see \cite{BLM}). 

\begin{lemma}[Bernstein's inequality]\label{l:bernstein}
Let $X_1,X_2,\dots,X_m$ be i.i.d.\ mean zero random variables, and assume that $|X_i|\leq K$ a.e
for some $K>0$. Then for any $t>0$ we have
$$
\Prob\Big\{\sum_{i=1}^m X_i>t\Big\}\leq \exp\bigg(-\frac{ct^2}{m\E X_1^2+Kt}\bigg),
$$
where $c>0$ is a universal constant.
\end{lemma}

The next theorem follows from Talagrand's concentration inequality for product measures (see \cite{Talagrand}).
\begin{theorem}[{for example, \cite[Corollary~4.10]{Ledoux}}]\label{th: tal}
Let $X_1,\dots,X_m$ be i.i.d.\ random variables with $|X_i|\leq 1$ a.e.
Then for any convex $1$--Lipschitz function $f:\R^m\to\R$ we have
$$\Prob\big\{\big|f(X_1,\dots,X_m)-\E f(X_1,\dots,X_m)\big|\geq t\big\}\leq 2e^{-ct^2},\quad t\geq 0,$$
where $c>0$ is a universal constant.
\end{theorem}

We will need the following estimate on the probability that the Erd\H os--Renyi random graph is tangle free. 

\begin{lemma}[{for example, \cite[Lemma~6.2]{MNS}}]\label{lem: tangle-free-prob}
Let $G$ be the Erd\H os--Renyi random graph on $n$ vertices, with no loops, and with parameter $p=d/(n-1)$. Then 
for any $\ell\geq 1$, we have 
$$
\Prob\{ G \text{ is $\ell$-tangle free}\}\geq 1-\frac{C\ell^3(2d)^{4\ell}}{n},
$$
where $C>0$ is a universal constant. 
\end{lemma}

\subsection{Construction of majorizers}\label{subs: major}

In the previous sections, we used vector majorizers in a rather abstract form, without discussing
whether the majorizers can be efficiently constructed in our random setting.
This subsection gives a probabilistic viewpoint to the notion.

We start with a general probabilistic construction. 
Let $\psi$ be a non-negative variable with unit expectation and an absolutely continuous distribution, uniformly
bounded above by a number $h\geq 1$, and let
$\kappa\geq 2$, $\tau\in(0,1]$.
For each $0< a< 1$, let $q_a$ be the quantile of $\psi$ of order $a$, that is the unique number satisfying
$\Prob\{\psi\leq q_a\}=a$.
Now, define an $n$-dimensional vector $\stnd=\stnd(\psi,h,\kappa,\tau)$ by setting
$$
\stnd_i: = \left\{
    \begin{array}{lll}
        h & \mbox{if }\quad  i\leq \tau\kappa, \\
        q_{(\kappa-i)/\kappa+\tau} & \mbox{if }\quad  \tau\kappa< i\leq (1+\tau)\kappa,\\
        0& \mbox{if }\quad  i>(1+\tau)\kappa.
    \end{array}
\right.
$$

\begin{lemma}\label{lem: bound-Y}
With the above definition, we have
$$
\kappa-h\leq \|\stnd\|_1\leq \kappa+(1+\tau\kappa) \, h.
$$
\end{lemma}
\begin{proof}
First note that 
$$
1=\E\psi=\int\limits_{0}^1 q_a\,da \geq \frac{1}{\kappa}
\sum\limits_{i\in [\tau\kappa+1,(1+\tau)\kappa]}  q_{(\kappa-i)/\kappa+\tau}.
$$
Hence, we get
$$
\|\stnd\|_1= \sum\limits_{i\in (\tau \kappa,(1+\tau)\kappa]}
q_{(\kappa-i)/\kappa+\tau} + \lfloor \tau\kappa\rfloor \, h
\leq \kappa+(1+\lfloor\tau\kappa\rfloor) \, h,
$$ 
and finish the proof of the upper bound. 

Similarly, using that $\psi$ is bounded by $h$, we can write 
$$
1=\E\psi=\int\limits_{0}^1 q_a\,da\leq \frac{1}{\kappa}\sum_{\ell=1}^{\kappa-1} q_{\ell/\kappa} + \frac{h}{\kappa}\leq \frac{1}{\kappa} \Vert \stnd\Vert_1+ \frac{h}{\kappa},
$$
and get the lower bound. 

\end{proof}

\begin{lemma}\label{lem: prob-heavy}
Let $n>\kappa\geq 2$ two integers, and $\tau\in(0,1]$.
Further, let
$X=(b_{i}\psi_i)_{i=1}^{n-\lfloor \tau\kappa/2\rfloor}$ be a random vector in $\R^{n-\lfloor \tau\kappa/2\rfloor}$,
such that $b_i,\psi_i$, $i=1,\dots,n-\lfloor \tau\kappa/2\rfloor$, are jointly independent; $\psi_i$ are equidistributed with $\psi$,
and $b_i$ are $0/1$ random variables with probability of success $\kappa/(n-1)$.
Define $\widetilde X:=X\oplus (h\,{\bf 1}_{\lfloor \tau\kappa/2\rfloor})\in\R^n$,
and let $\widetilde X^*$ be the non-increasing rearrangement of $\widetilde X$.
Then with probability at least $1-\kappa\exp(-c_{\text{\tiny\ref{lem: prob-heavy}}}\tau^2\kappa)$ we have
$\stnd\geq \widetilde X^*$
coordinate-wise.
Here, $c_{\text{\tiny\ref{lem: prob-heavy}}}>0$ is a universal constant.
\end{lemma}
\begin{proof}
Without loss of generality, suppose that $\tau^2\kappa$ is bounded from below by a large universal constant.
First, we recall the following consequence of Chernoff's inequality.
Let $m\in\N$ and let $W=(W_i)_{i=1}^m$ be a random vector with i.i.d coordinates equidistributed with $\psi$.
Then for any $a\in(0,1)$ with $i\geq (1-a)m$ we have
$$\Prob\big\{W^*_i> q_a\big\}\leq \exp\bigg(-\frac{(i-(1-a)m)^2}{2i}\bigg).$$
We will apply this relation to the vector $\widetilde X$ conditioned on an upper bound for the sum $\sum\limits_i b_i$.
Namely, define the event
$$\Event:=\bigg\{\sum\limits_{i=1}^{n-\lfloor \tau\kappa/2\rfloor} b_i\leq \kappa+\lfloor \tau\kappa/4\rfloor\bigg\}.$$
Then from the above we get that, conditioned on $\Event$, we have for $m:=\kappa+\lfloor \tau\kappa/4\rfloor$
and any $\tau\kappa< i\leq (1+\tau)\kappa$ and $i\geq \lfloor \tau\kappa/2\rfloor+(1-a)m$,
$$
\Prob\big\{\widetilde X^*_i> q_a\mid\Event\big\}
\leq \Prob\big\{X^*_{i-\lfloor \tau\kappa/2\rfloor}> q_a\mid\Event\big\}
\leq \exp\bigg(-\frac{(i-\lfloor \tau\kappa/2\rfloor-(1-a)m)^2}{2(i-\lfloor \tau\kappa/2\rfloor)}\bigg).
$$
Taking $a:=(\kappa-i)/\kappa+\tau$ and using the definition of $\stnd$, we obtain
$$
\Prob\big\{\widetilde X^*_i> \stnd_i\mid\Event\big\}\leq 2\exp(-c\tau^2\kappa),\quad \tau\kappa< i\leq (1+\tau)\kappa,
$$
for some universal constant $c>0$.
It remains to note that, by the Bernstein inequality, the probability of $\Event$ can be bounded from below by $1-2\exp(-c'\tau^2 \kappa)$.
The result follows.
\end{proof}

\bigskip

We are now ready to define majorizers for our random variables of interest. 
Fix parameters $n$ and $d$ and $h\geq 2$. Let $b$ be a Bernoulli random variable with probability of success $d/(n-1)$ and let 
$\xi$ be a random variable independent of $b$,
with an absolutely continuous distribution, of zero mean, unit variance and with $\Prob\{\xi^2\leq h\}=1$. 
Denote $\varepsilon_0:= (h\log\log\log n)^{-1}-d^{-1}$. 
In what follows, we call 
\begin{equation}\label{eq: definition-Y}
\stnd:= \stnd\big(\xi^2,h, d, \varepsilon_0\big)
\end{equation}
{\it the standard majorizer} (with respect to $\xi$) and we will suppress parameters whenever they are clear from the context. 
Note that by Lemma~\ref{lem: bound-Y}, we have 
\begin{equation}\label{eq: bound-Y}
d-h\leq \Vert \stnd\Vert_1\leq \big(1+(\log\log\log n)^{-1}\big) d. 
\end{equation}

Let $M=(\mu_{ij})$ be an $n\times n$ symmetric random matrix with zero diagonal and i.i.d.\ entries
(up to the symmetry constraint).
Assume that each off-diagonal entry of $M$ has the form $\mu_{ij}=b_{ij}\xi_{ij}$,
where $b_{ij}$ is $0/1$ random variable
with probability of success $d/(n-1)$, and $\xi_{ij}$ is equidistributed with $\xi$.

Fix for a moment any realization of $M$. We denote by $G_M=([n],E_M)$ the graph
with the edge set $E_M:=\{i\edg j:\;b_{ij}=1\}$.
With some abuse of terminology, we will say that
a vertex $v$ of $G_M$ is {\it majorized} by a vector $y\in\R^n_+$ if
the non-increasing rearrangement of
the sequence $(\mu_{vi}^2)_{i=1}^n$ is majorized (coordinate-wise) by $y$.
We say that a vertex $v$ of the graph $G_M$ is {\it heavy}
if it is not majorized by $\stnd$.
Note that, in particular, every non-heavy vertex has at most $\w d$ neighbors in $G_M$, where $\w d= \big( 1+ \varepsilon_0\big)d$. 

In the next lemma, we show that with a large probability any given vertex has a relatively
small number of heavy neighbors in $G_M$. 

\begin{lemma}\label{lem: nb-heavy-neighbors}
Assume that $2d^2\leq \exp(c_{\text{\tiny\ref{lem: prob-heavy}}}\varepsilon_0^2d/2)$.
Then for any integer $q\leq \varepsilon_0 d/2-1$ we have
$$
\big\vert\big\{ j\edg v:\,  j \text{ is {\it heavy}}\big\}\big\vert \leq q\quad \mbox{for all }v\in[n]
$$
with probability at least $1-n\,\exp(-c_{\text{\tiny\ref{lem: prob-heavy}}}\varepsilon_0^2d q/2)$. 
\end{lemma}
\begin{proof}
Fix for a moment any distinct indices $i,j_1,j_2,\dots,j_q\in[n]$.
We will estimate probability of the event
$$
\Event:=\big\{\mbox{$j_\ell$ is adjacent to $i$ for all $\ell\leq q$, and all vertices $j_1,\dots,j_q$ are heavy}\big\}.
$$
For every $\ell\leq q$, let $X_\ell^*\in\R^{n-q-1}$ be the non-increasing rearrangement of the sequence
$(\mu_{uj_\ell}^2)_{u\in[n]\setminus\{i,j_1,\dots,j_q\}}$.
It is then not difficult to see that a necessary condition for $j_\ell$ to be heavy is
$$\stnd\;\;\mbox{ is {\it not} a majorizer for }\;\;(h\,{\bf 1}_{q+1})\oplus X_\ell^*.$$
This latter condition allows to estimate the probability of $\Event$ via a decoupling:
we get
\begin{align*}
\Prob(\Event)&\leq 
\Prob\big\{\mbox{$j_\ell$ is adjacent to $i$ for all $\ell\leq q$}\big\}
\prod\limits_{\ell=1}^q
\Prob\big\{\stnd\;\;\mbox{ is {\it not} a majorizer for }\;\;(h\,{\bf 1}_{q+1})\oplus X_\ell^*\big\}\\
&\leq \big(d/(n-1)\big)^q\big(d\exp(-c_{\text{\tiny\ref{lem: prob-heavy}}}\varepsilon_0^2d)\big)^{q},
\end{align*}
where at the second step we applied Lemma~\ref{lem: prob-heavy}.

It remains to take the union bound over all possible choices of indices $i,j_1,j_2,\dots,j_q\in[n]$:
we have from the above
\begin{align*}
\Prob\big\{
\big\vert\big\{ j\edg v:\,  j \text{ is {\it heavy}}\big\}\big\vert \leq q\mbox{ for all }v\in[n]
\big\}
&\geq 1-n^{q+1}\big(d/(n-1)\big)^q\big(d\exp(-c_{\text{\tiny\ref{lem: prob-heavy}}}\varepsilon_0^2d)\big)^{q}\\
&\geq 1-n\,\exp(-c_{\text{\tiny\ref{lem: prob-heavy}}}\varepsilon_0^2d q/2).
\end{align*}
The result follows.
\end{proof}

\medskip
We summarize the results above in the next proposiiton. 

\begin{prop}[Majorizers]\label{prop: heavyn}
Let $\delta\in (0,1/3)$. Assume that $n\geq C$ and $d\geq C\log^{\frac{1}{1+\delta}} n$,
and that the random matrix $M=(b_{ij}\xi_{ij})$ and the graph $G_M$ are as above. 
Define the event
\begin{align*}
\Event_{mjr}(\delta):=\Big\{
&\forall\,i\leq n,\;\mbox{the vector $(b_{ij})_{j=1}^n$ has at most $d^{1+\delta}$ non-zero components
AND}\\
&\mbox{for any vertex $v\in[n]$ the number of its heavy neighbors is at most $d^{\frac89}$}\Big\}.
\end{align*}
 Then we have 
$\Prob(\Event_{mjr})\geq 1-\exp(-cd^{1+\delta})$.
Here, $C,c>0$ are universal constants. 
\end{prop}
\begin{proof}
Applying the Bernstein inequality (Lemma~\ref{l:bernstein}), we get for any $i\leq n$:
$$
\Prob\big\{\mbox{the vector $(b_{ij})_{j=1}^n$ has at least $d^{1+\delta}$ non-zero components}\big\}
\leq 2\exp(-cd^{1+\delta} )
$$
for a universal constant $c>0$.
Together with Lemma~\ref{lem: nb-heavy-neighbors},
this gives the result.

\end{proof}

\section{Contribution of paths with a large number of multiplicity one edges}\label{s: many mult1}

The goal of this section is to provide a bound on the sum of paths weights over all paths having a large number of multiplicity one edges. More precisely, let $M=(\mu_{ij})$ be an $n\times n$ random symmetric matrix
with zero diagonal and with i.i.d entries above the diagonal of the form $b_{ij}\xi_{ij}$, where $b_{ij}$ are i.i.d
Bernoulli random variables with probability of success $d/(n-1)$, and
$\xi_{ij}$ are i.i.d random variables (independent from $b_{ij}$)
of mean zero, variance one, and uniformly bounded above by $\sqrt{h}$.  
Let
\begin{equation}\label{eq: maxnorm act def}
\maxnorm:=\big(1+(\log\log\log n)^{-1}\big)^3\, \E\max\limits_{i\leq n}\sum\limits_{j=1}^n \mu_{ij}^2,
\end{equation}
and 
\begin{equation}\label{eq: dmax act def}
\dmax:=\big(1+(\log\log\log n)^{-1}\big)\, \E\max\limits_{i\leq n}\sum\limits_{j=1}^n b_{ij}.
\end{equation}
Note that
\begin{equation}\label{eq: aux 043209275-395}
\begin{split}
\maxnorm&=\big(1+(\log\log\log n)^{-1}\big)^3\, \E\bigg(\E
\Big[\max\limits_{i\leq n}\sum\limits_{j=1}^n b_{ij}\xi_{ij}^2\;\big|\;(b_{ij})\Big]\bigg)\\
&\geq \big(1+(\log\log\log n)^{-1}\big)^3\, \E\max\limits_{i\leq n}\sum\limits_{j=1}^n b_{ij}
=(1+(\log\log\log n)^{-1})^2\dmax.
\end{split}
\end{equation}
Define two events
\begin{equation}\label{eq: event maxnorm def}
\Event_{\maxnorm}:=\bigg\{\sum\limits_{j=1}^n \mu_{ij}^2\leq \frac{\maxnorm}{1+(\log\log\log n)^{-1}}\mbox{ for all $i\in[n]$}\bigg\},
\end{equation}
and
\begin{equation}\label{eq: event g def}
\Event_g:=\big\{\mbox{$G$ is $(k/\log^2 d)$--tangle free and }\deg_i(G)\leq d_{\max}\mbox{ for all $i\in[n]$}\big\},
\end{equation}
where $G$ is the random graph on $[n]$ with the adjacency matrix $(b_{ij})$. Our goal in this section is to bound the quantity 
\begin{equation}\label{eq: path-weight-large-m1}
\E\Big(\sum\limits_{\multone(\path)\geq \frac{k}{\log^{11}k}}\, 
\prod\limits_{\ell=1}^{2k}\mu_{\path(\ell-1),\path(\ell)}\mathbf{1}_{\Event_{\maxnorm}\cap\Event_g}\Big),
\end{equation}
where the summation is taken over all closed paths on $K_{[n]}$ of length $2k$
with $\multone(\path)\geq \frac{k}{\log^{11}k}$. Recall that $\multone(\path)$ is the cardinality of a largest subset $S$ of edges of $K_{[n]}$
such that each edge from $S$ is travelled by $\path$ exactly once, and no two edges from $S$ are incident.

In the classical applications of the moment method for random matrices,
one takes advantage of the fact that the entries are independent and
centered in order to eliminate all paths having edges of multiplicity one.
In our case, multiplying by the indicator $\mathbf{1}_{\Event_{\maxnorm}\cap\Event_g}$
produces complex dependencies; and the quantity $\prod\limits_{\ell=1}^{2k}\mu_{\path(\ell-1),\path(\ell)}\mathbf{1}_{\Event_{\maxnorm}\cap\Event_g}$
may be not centered. 
Nevertheless, in some sense the event $\Event_{\maxnorm}$ only affects the distribution of entries $\mu_{i,j}$ when 
one of the rows $i$ or $j$ is heavy i.e.\ its squared Euclidean norm is close to $\maxnorm$.
Informally, when the Euclidean norm of rows $i$ and $j$ is significantly below $\sqrt{\maxnorm}$,
any choice of value for the $(i,j)$-th entry will keep the matrix realization in the event $\Event_{\maxnorm}$.  This idea is made precise in the next lemma. 

\begin{lemma}\label{l: pre-cancel}
Let $\Xi=(\xi_{ij})$ be an $n\times n$ random symmetric matrix with zero diagonal
whose non-diagonal entries are independent (up to the symmetry constraint) centered random variables
uniformly bounded by $\sqrt{h}$ (the entries do not need to be identically distributed). Let $r>0$ and let $\Event$ be the event 
$$
\Event:=\big\{ \Vert \Row_i\Vert_2^2\leq r \text{ for all } i\in [n]\big\},
$$
where $\Row_i$ stands for the $i$-th row of $\Xi$. Assume that $\Prob(\Event)>0$. 
Let $S\subset {[n]\choose 2}$ be a set of unordered couples such that for any
$\{i,j\}, \{i',j'\}\in S$ we have $\{i,j\}\cap \{i',j'\}=\emptyset$. Denote by $\Event_S$ the event 
$$
\Event_S:=\Big\{\text{For every $\{i,j\}\in S$},\,  \max\big(\Vert \Row_i\Vert_2^2, \Vert \Row_j\Vert_2^2\big)\geq r+\xi_{ij}^2-h\Big\}.
$$
Then for any multiset $E\subset {[n]\choose 2}$ of unordered couples 
such that $S\subset E$ and all elements of $S$ in $E$ have multiplicity one,
$$
\E\, \Big[ \prod_{\{i,j\}\in E} \xi_{ij}\mathbf{1}_{\Event}\Big] = \E\, \Big[\prod_{\{i,j\}\in E} \xi_{ij} \mathbf{1}_{\Event_{S}\cap\Event}\Big]
$$
\end{lemma}
\begin{proof}
Let us enumerate the elements of $S$ as $\{i_1,j_1\},\ldots, \{i_s,j_s\}$. 
For any $u\leq s$, denote 
$$
\Event_u:=\Big\{ \max\big(\Vert \Row_{i_u}\Vert_2^2, \Vert \Row_{j_u}\Vert_2^2\big)\geq r+\xi_{i_uj_u}^2-h\Big\}
=\bigg\{ \max\Big(\sum\limits_{j\neq j_u}\xi_{i_u j}^2, \sum\limits_{i\neq i_u}\xi_{i j_u}^2\Big)\geq r-h\bigg\}.
$$
Obviously,
$$
\E\, \Big[ \prod_{\{i,j\}\in E} \xi_{ij}\mathbf{1}_{\Event}\Big]= \E\, \Big[ \prod_{\{i,j\}\in E} \xi_{ij}\mathbf{1}_{\Event_s}\mathbf{1}_{\Event}\Big] + \E\, \Big[ \prod_{\{i,j\}\in E} \xi_{ij}\mathbf{1}_{\Event_s^c}\mathbf{1}_{\Event}\Big].
$$
Since $\xi_{i_sj_s}^2\leq h$, everywhere on $\Event_s^c$ we have
$\max\big(\Vert \Row_{i_s}\Vert_2^2, \Vert \Row_{j_s}\Vert_2^2\big)<r$.
Thus, if we denote by $\Event^s$ the event
$$
\Event^s:=\big\{ \Vert \Row_i\Vert_2^2\leq r \text{ for all } i\in [n]\setminus\{i_s,j_s\}\big\}, 
$$
then $\mathbf{1}_{\Event}\mathbf{1}_{\Event_s^c}=\mathbf{1}_{\Event_s^c} \mathbf{1}_{\Event^s}$
everywhere on the probability space.
The crucial observation is that $\mathbf{1}_{\Event_s^c} \mathbf{1}_{\Event^s}$ depends only
on the variables $\{\xi_{ij}\}_{\{i,j\}\neq \{i_s,j_s\}}$, implying that
the product $\prod_{\{i,j\}\in E\setminus \{i_s,j_s\}} \mathbf{1}_{\Event_s^c} \mathbf{1}_{\Event^s}$
is independent from $\xi_{i_sj_s}$.
Thus, we can write
\begin{align*}
\E\, \Big[ \prod_{\{i,j\}\in E} \xi_{ij}\mathbf{1}_{\Event}\mathbf{1}_{\Event_s^c}\Big]
&=\E\, \Big[ \prod_{\{i,j\}\in E} \xi_{ij}\mathbf{1}_{\Event_s^c} \mathbf{1}_{\Event^s}\Big]\\
&=
\E\,\Big[\prod_{\{i,j\}\in E\setminus \{i_s,j_s\}} \xi_{ij} \mathbf{1}_{\Event}\mathbf{1}_{\Event_s^c}\Big]\cdot
\E\xi_{i_sj_s} =0.
\end{align*}
We deduce that 
$$
\E\, \Big[ \prod_{\{i,j\}\in E} \xi_{ij}\mathbf{1}_{\Event}\Big]= \E\, \Big[ \prod_{\{i,j\}\in E} \xi_{ij}\mathbf{1}_{\Event_s}\mathbf{1}_{\Event}\Big] .
$$
Next, we write 
$$
\E\, \Big[ \prod_{\{i,j\}\in E} \xi_{ij}\mathbf{1}_{\Event_s}\mathbf{1}_{\Event}\Big]= \E\, \Big[ \prod_{\{i,j\}\in E} \xi_{ij}\mathbf{1}_{\Event_s}\mathbf{1}_{\Event_{s-1}}\mathbf{1}_{\Event}\Big]+ \E\, \Big[ \prod_{\{i,j\}\in E} \xi_{ij}\mathbf{1}_{\Event_s}\mathbf{1}_{\Event_{s-1}^c}\mathbf{1}_{\Event}\Big].
$$
Similarly to the above,
$\mathbf{1}_{\Event_{s-1}^c}\mathbf{1}_{\Event}=\mathbf{1}_{\Event_{s-1}^c} \mathbf{1}_{\Event^{s-1}}$
everywhere on the probability space, where $
\Event^{s-1}:=\big\{ \Vert \Row_i\Vert_2^2\leq r \text{ for all } i\in [n]\setminus\{i_{s-1},j_{s-1}\}\big\}
$. Since $\mathbf{1}_{\Event_s}\mathbf{1}_{\Event_{s-1}^c} \mathbf{1}_{\Event^{s-1}}$
is independent from $\xi_{i_{s-1}j_{s-1}}$, 
by repeating the above argument we get
$$
\E\, \Big[ \prod_{\{i,j\}\in E} \xi_{ij}\mathbf{1}_{\Event}\Big]= \E\, \Big[ \prod_{\{i,j\}\in E} \xi_{ij}\mathbf{1}_{\Event_s}\mathbf{1}_{\Event_{s-1}}\mathbf{1}_{\Event}\Big].
$$
It remains to re-run this procedure and note that $\prod_{\ell=1}^s \mathbf{1}_{\Event_\ell}= \mathbf{1}_{\Event_S}$ to finish the proof. 
\end{proof}

The next lemma will allow us to estimate the probability of the event $\Event_S$ appearing in Lemma~\ref{l: pre-cancel}. 

\begin{lemma}\label{l: concentration-mult1}
Let $\Xi$ be an $n\times n$ random symmetric matrix with zero diagonal whose entries above the diagonal are
independent centered random variables uniformly bounded (in absolute value) by $\sqrt{h}$. 
Let $\widetilde r, \delta >0$ be such that for any $i\in [n]$ 
$$
\Prob\Big\{ \Vert \Row_i\Vert_2^2\geq \widetilde r\Big\}\leq \delta,
$$
where $\Row_i$ stands for the $i$-th row of $\Xi$. Let $S\subset {[n] \choose 2}$
be a set of unordered couples such that for any $\{i,j\}, \{i',j'\}\in S$ we have $\{i,j\}\cap \{i',j'\}=\emptyset$. Then, we have 
$$
\Prob\Big\{ \text{For any }\{i,j\}\in S,\, \max\big(\Vert \Row_i\Vert_2^2, \Vert \Row_j\Vert_2^2\big)\geq \widetilde r + 2\vert S\vert h\Big\}\leq (2\delta)^{\vert S\vert}. 
$$
\end{lemma}
\begin{proof}
Denote by $\widetilde \Event$ the event 
$$
\widetilde \Event:=\Big\{ \text{For any }\{i,j\}\in S,\, \max\big(\Vert \Row_i\Vert_2^2, \Vert \Row_j\Vert_2^2\big)\geq \widetilde r + 2\vert S\vert h\Big\}.
$$
Let us enumerate elements of $S$ as $\{i_1,j_1\},\ldots, \{i_s,j_s\}$. 
For any $i\in [n]$ and any $J\subset [s]$,
we denote by $\Row_i\setminus J$ the vector obtained from $\Row_i$ by removing the coordinates indexed by $i_\ell, j_\ell$ with $\ell \in J$. 
With these notations, and using that the entries are uniformly bounded by $\sqrt{h}$, we can write 
$$
\Prob\big(\widetilde \Event\big) \leq \Prob\Big(\bigcap_{\ell=1}^s \big\{\max\big(\Vert \Row_{i_\ell}\setminus \{\ell+1,\ldots,s\}\Vert_2^2, \Vert \Row_{j_\ell}\setminus \{\ell+1,\ldots,s\}\Vert_2^2\big) \geq \widetilde r+ 2\ell h\big\}\Big).
$$
It follows from the independence of the entries of $\Xi$ that the events indexed by $\ell$ in the above intersection are independent. Therefore, we deduce 
\begin{align*}
\Prob\big(\widetilde \Event\big)& \leq \prod_{\ell=1}^s \Prob\Big\{ \max\big(\Vert \Row_{i_\ell}\setminus \{\ell+1,\ldots,s\}\Vert_2^2, \Vert \Row_{j_\ell}\setminus \{\ell+1,\ldots,s\}\Vert_2^2\big) \geq \widetilde r+ 2\ell h\Big\}\\
&\leq \prod_{\ell=1}^s \Prob\Big\{ \max\big(\Vert \Row_{i_\ell}\Vert_2^2, \Vert \Row_{j_\ell}\Vert_2^2\big) \geq \widetilde r\Big\}\leq (2\delta)^s,
\end{align*}
where in the last step we used the union bound together with the hypothesis of the lemma.
\end{proof}

We are now ready to provide an upper bound for the quantity in \eqref{eq: path-weight-large-m1}. Recall that we restricted our attention to paths with
$\multone(\path)\geq k/\log^{11}k$ as the complementary regime was treated previously in Section~\ref{s: few mult1}
for a class of deterministic matrices.
The main statement of this section is the following. 
\begin{prop}\label{prop: mult-one}
Let $h\geq 2$, let $n\geq n_0(h)$ be a large integer, 
let $M=(\mu_{ij})$ be an $n\times n$ random symmetric matrix
with zero diagonal and with i.i.d entries above the diagonal, with the $(i,j)$--th
entry of the form $b_{ij}\xi_{ij}$, where $b_{ij}$ are i.i.d
Bernoulli random variables with probability of success $d/(n-1)$, and
$\xi_{ij}$ are i.i.d random variables (independent from $b_{ij}$)
of mean zero, variance at most one, and uniformly bounded above by $\sqrt{h}$.
Let $\maxnorm$ and $\dmax$ be defined by \eqref{eq: maxnorm act def}
and \eqref{eq: dmax act def}, respectively,
and events $\Event_{\maxnorm}$ and $\Event_g$
--- by \eqref{eq: event maxnorm def} and \eqref{eq: event g def}, respectively,
with $G$ being the random graph on $[n]$ with the adjacency matrix $(b_{ij})$. 
Then for any positive integer $k$ satisfying 
$$
\sqrt{\log n}\leq \sqrt{k}\leq \frac{\maxnorm}{32 \big(h \log\log\log n)^2 \sqrt{\log (h\maxnorm^2)}},
$$
we have 
$$\E\Big(\sum\limits_{\multone(\path)\geq \frac{k}{\log^{11}k}}\, 
\prod\limits_{\ell=1}^{2k}\mu_{\path(\ell-1),\path(\ell)}\mathbf{1}_{\Event_{\maxnorm}\cap\Event_g}\Big)
\leq n
,$$
where the summation is taken over all closed paths on $K_{[n]}$ of length $2k$
with $\multone(\path)\geq \frac{k}{\log^{11}k}$. 
\end{prop}
We note that ``$n$'' on the right hand side of the estimate can be replaced, without affecting the rest of the argument,
with anything ``small enough'', of order $(2+o(1))^{2k}d^k$.
\begin{proof}
We start by writing 
$$\E\Big(\sum\limits_{\multone(\path)\geq \frac{k}{\log^{11}k}}\, 
\prod\limits_{\ell=1}^{2k}\mu_{\path(\ell-1),\path(\ell)}\mathbf{1}_{\Event_{\maxnorm}\cap\Event_g}\Big)
= \E\Big(\sum\limits_{\multone(\path)\geq \frac{k}{\log^{11}k}}\, 
\E\Big[\prod\limits_{\ell=1}^{2k}\mu_{\path(\ell-1),\path(\ell)}{\bf 1}_{\Event_{\maxnorm}\cap\Event_g}
\;\Big\vert\;G\Big]\Big),
$$
and aim first at bounding $
\E\Big[\prod\limits_{\ell=1}^{2k}\mu_{\path(\ell-1),\path(\ell)}{\bf 1}_{\Event_{\maxnorm}\cap\Event_g}
\;\Big\vert\;G\Big]$
for a fixed path $\path$ on $K_n$ with $\multone(\path)\geq k/\log^{11}k$. 
Fix a realization $G_0$ of $G$ so that $\Event_g$ holds, and let $\Xi$ be the $n\times n$
random matrix obtained as the Hadamard product of the adjacency matrix of $G_0$ and the matrix $(\xi_{ij})_{1\leq i,j\leq n}$. 
Note that if the path $\path$ is not contained in $G_0$ then the expectation is zero.
Below, we assume that all edges traversed by $\path$ are contained in $G_0$.
Denoting by $\Row_i(\Xi)$ the $i$-th row of $\Xi$ and setting
$$
\widetilde\Event:=\Big\{\Vert \Row_i(\Xi)\Vert_2^2\leq \frac{\maxnorm}{1+(\log\log\log n)^{-1}} \text{ for all } i\in [n]\Big\},
$$
we can write 
$$
\E\Big[\prod\limits_{\ell=1}^{2k}\mu_{\path(\ell-1),\path(\ell)}{\bf 1}_{\Event_{\maxnorm}\cap\Event_g}\;\Big\vert\;G=G_0\Big]
= \E\, \Big[ \prod_{\{i,j\}\in E} \xi_{ij}{\bf 1}_{\widetilde \Event}\Big], 
$$
where $E$ denotes  the multiset of edges traversed by $\path$.
Applying Lemma~\ref{l: pre-cancel} and using that the $\xi_{ij}$'s are bounded by $\sqrt{h}$, we get 
$$
\E\Big[\prod\limits_{\ell=1}^{2k}\mu_{\path(\ell-1),\path(\ell)}{\bf 1}_{\Event_{\maxnorm}\cap\Event_g}\;\Big\vert\;G=G_0\Big]
\leq h^{k}\,  \Prob\big(\Event_S\cap \widetilde \Event\big),
$$
where $S=S(\path,G_0)$
is some fixed subset of edges of $G_0$ of size
$$s:=\big\lceil\min \big(k/\log^{11}k, \maxnorm/(32h \log\log\log n)\big)\big\rceil$$
such that each edge from $S$ is travelled by $\path$ exactly once, no two edges from $S$ are incident, and 
$$
\Event_S:=\Big\{\text{For every $\{i,j\}\in S$},\,  \max\big(\Vert \Row_i(\Xi)\Vert_2^2, \Vert \Row_j(\Xi)\Vert_2^2\big)\geq \maxnorm+\xi_{ij}^2-h\Big\}.
$$
Note that such $S$ exists since $\multone(\path)\geq k/\log^{11}k$.
The choice of $s$ ensures that 
$h+ 2sh \leq (8\log\log\log n)^{-1}\maxnorm$, 
so
we can write  for any $i\leq n$ 
\begin{align*}
\Prob&\Big\{ \Vert \Row_i(\Xi)\Vert_2^2\geq \frac{\maxnorm}{1+(\log\log\log n)^{-1}} -h-2\vert S\vert h\Big\}\\
&
\leq \Prob\Big\{ \Vert \Row_i(\Xi)\Vert_2^2-\E \Vert \Row_i(\Xi)\Vert_2^2\geq (8\log\log\log n)^{-1}\maxnorm \Big\}\\
&\leq \exp\Big(- \frac{ \maxnorm^2}{32\deg_{G_0}(i) \big(h\log\log\log n)^2}\Big),
\end{align*}
where we used that $\E \Vert \Row_i(\Xi)\Vert_2^2\leq \dmax\leq\big(1+(\log\log\log n)^{-1}\big)^{-2} \maxnorm$
(see \eqref{eq: aux 043209275-395})
for the first inequality, and Hoeffding's inequality in the last step.
Applying Lemma~\ref{l: concentration-mult1} and using the previous relation,
we get
$$
\Prob\big(\Event_S\big) \leq 2^{s}\exp\Big(- \frac{s \, \maxnorm^2}{32\dmax \big(h\log\log\log n)^2}\Big),
$$
whence
$$
\E\Big[\prod\limits_{\ell=1}^{2k}\mu_{\path(\ell-1),\path(\ell)}{\bf 1}_{\Event_{\maxnorm}\cap\Event_g}\;\Big\vert\;G=G_0\Big]
\leq h^{k}2^{s}\exp\Big(- \frac{s \, \maxnorm^2}{32\dmax \big(h\log\log\log n)^2}\Big).
$$
Since there are at most $n\dmax^{2k}$ distinct closed paths of length $2k$ on $G_0$, we get 
\begin{align*}
\sum\limits_{\multone(\path)\geq \frac{k}{\log^{11} k}}\, 
&\E\Big[\prod\limits_{\ell=1}^{2k}\mu_{\path(\ell-1),\path(\ell)}{\bf 1}_{\Event_{\maxnorm}\cap\Event_g}\;\Big\vert\;G=G_0\Big]\\
&\leq n(h\dmax^{2})^k 2^{s}\exp\Big(- \frac{s \, \maxnorm^2}{32\dmax \big(h\log\log\log n)^2}\Big).
\end{align*}
Since $\dmax\leq \maxnorm$, we deduce that 
$$
\E\Big(\sum\limits_{\multone(\path)\geq \frac{k}{\log^{11}k}}\, 
\prod\limits_{\ell=1}^{2k}\mu_{\path(\ell-1),\path(\ell)}\mathbf{1}_{\Event_{\maxnorm}\cap\Event_g}\Big)
\leq n(h\maxnorm^{2})^k 2^{s}\exp\Big(- \frac{s \, \maxnorm}{32 \big(h\log\log\log n)^2}\Big).
$$
It remains to use the condition on $k$ to finish the proof. 
\end{proof}

\section{Upper bound for the operator norm}

We now have everything in place in order to complete the  proof of the upper bound in Theorem~A. 
We will start by showing that ``small perturbations'' of the quantity $\rho_n$:
replacing the maximum with its expectation or adding/removing the matrix diagonal ---
have no effect on the final (asymptotic) result. 
We make this precise in the next simple lemma. 

\begin{lemma}\label{lem: upper-zero-diag}
Let $h\geq 1$.
For each $n$, let $W_n$ be $n\times n$ symmetric random matrix with i.i.d.\ entries
above and on the main diagonal,
with each entry equidistributed with the product $b_n\xi_n$, where
$\xi_n$ is a real random variable with $\E\xi_n^2=1$ and $\xi_n^2\leq h$ a.e., and
$b_n$ is $0/1$ (Bernoulli) random variable independent of $\xi_n$,
with probability of success equal to $p_n$. 
For each $n$, set $M_n= W_n-{\rm Diag}(W_n)$ and $d=d(n)=(n-1)p_n$. 
Assume further that $n p_n\to\infty$ with $n$ and  denote
\begin{equation*}
\begin{split}
\rho_n&:=\theta_n+\frac{n p_n}{\theta_n},\quad \theta_n:=\sqrt{\max\big(\max\limits_{i\leq n}\|\Row_i(W_n)\|_2^2-np_n,n p_n\big)},\\
\rho_n'&:=\theta_n'+\frac{n p_n}{\theta_n'},\quad \theta_n':=\sqrt{\max\big(\E\max\limits_{i\leq n}\|\Row_i(W_n)\|_2^2-np_n,n p_n\big)},
\end{split}
\end{equation*}
and 
\begin{equation*}
\begin{split}
\w\rho_n&:=\w\theta_n+\frac{d}{\w\theta_n},\quad \w\theta_n:=\sqrt{\max\big(\max\limits_{i\leq n}\|\Row_i(M_n)\|_2^2-d,d\big)},\\
\w\rho_n'&:=\w\theta_n'+\frac{d}{\w\theta_n'},\quad \w\theta_n':=\sqrt{\max\big(\E\max\limits_{i\leq n}\|\Row_i(M_n)\|_2^2-d,d\big)}.
\end{split}
\end{equation*} 
Then the sequences $\big(\frac{\Vert W_n\Vert}{\Vert M_n\Vert}\big)_{n\geq 1}$,
$\big(\frac{\rho_n}{{\vphantom{\w A}}\rho_n'}\big)_{n\geq 1}$,
$\big(\frac{\w\rho_n}{{\vphantom{\w A}}\w\rho_n'}\big)_{n\geq 1}$,
$\big(\frac{\rho_n}{{\vphantom{\w A}}\w\rho_n}\big)_{n\geq 1}$ converge to one in probability. 
\end{lemma}
\begin{proof}
We start by noticing that 
$$
\Vert M_n\Vert-\sqrt{h}\leq \Vert W_n\Vert \leq \Vert M_n\Vert +\sqrt{h},
$$
where we have used that the absolute values of the
entries are uniformly bounded by $\sqrt{h}$. Therefore, for any fixed $\varepsilon\in (0,1)$, we have 
$$
\Prob\Big\{\Big\vert \frac{\Vert W_n\Vert}{\Vert M_n\Vert}-1\Big\vert \geq \varepsilon\Big\} 
\leq \Prob\Big\{ \Vert M_n\Vert \leq \frac{\sqrt{h}}{\varepsilon}\Big\}\leq \Prob\Big\{ \Vert \Row_1(M_n)\Vert_2 \leq \frac{\sqrt{h}}{\varepsilon}\Big\}.
$$
Using that $\E\Vert \Row_1(M_n)\Vert_2^2= d\underset{n\to \infty}{\longrightarrow} \infty$ together with Bernstein's inequality (Lemma~\ref{l:bernstein}), we deduce the first assertion of the lemma. 

Further, since
$$
\max\limits_{i\leq n}\|\Row_i(M_n)\|_2^2\leq \max\limits_{i\leq n}\|\Row_i(W_n)\|_2^2\leq \max\limits_{i\leq n}\|\Row_i(M_n)\|_2^2+h,
$$
since the function $(x,y)\to \sqrt{\max(x-y,y)}+\frac{y}{\sqrt{\max(x-y,y)}}=\frac{\max(x,2y)}{\sqrt{\max(x-y,y)}}$
is coordinate-wise increasing, and $\sqrt{\max(x'-y,y)}\geq\sqrt{\max(x-y,y)}-\frac{x-x'}{2\sqrt{y}}$ for all $0<x'\leq x$,
we have 
$$
\w\rho_n\leq \rho_n\leq (1+C/n)\w\rho_n+ \frac{h}{2\sqrt{np_n}};
\quad \w\rho_n'\leq \rho_n'\leq (1+C/n)\w\rho_n'+ \frac{h}{2\sqrt{np_n}}.
$$
Since $\lim_n np_n=\infty$, we deduce that $\big(\frac{\w\rho_n}{\rho_n}\big)_{n\geq 1}$
and $\big(\frac{\w\rho_n'}{\rho_n'}\big)_{n\geq 1}$
converge to one in probability. 
Therefore, if we show that $\big(\frac{\w\rho_n}{\w\rho_n'}\big)_{n\geq 1}$ converges to one in probability,
the statement will be proved. 

If we define $f_n(x):= \frac{\max(x^2,2d)}{\sqrt{ \max(x^2-d,d)}}$,
then, $f_n$ is increasing on $\R_+$ and $\vert f_n(x)-f_n(x')\vert \leq 2\vert x-x'\vert$
for any $x,x'\in \R_+$.
Using this and noticing that $\w\rho_n= f_n\big(\max\limits_{i\leq n}\|\Row_i(M_n)\|_2\big)$
and $\w\rho_n'= f_n\Big(\sqrt{\E\max\limits_{i\leq n}\|\Row_i(M_n)\|_2^2}\Big)$, we can write for any $\varepsilon\in (0,1)$ 
$$
\Prob\Big\{\Big\vert\frac{\w\rho_n'}{\w\rho_n}-1\Big\vert \geq \varepsilon\Big\} 
\leq \Prob\Big\{ \Big\vert \max\limits_{i\leq n}\|\Row_i(M_n)\|_2-\sqrt{\E \max\limits_{i\leq n}\|\Row_i(M_n)\|_2^2}
\Big\vert \geq \varepsilon \w \rho_n/2\Big\}.
$$
Now using that 
$$
\Big\vert \E \max\limits_{i\leq n}\|\Row_i(M_n)\|_2-\sqrt{\E \max\limits_{i\leq n}\|\Row_i(M_n)\|_2^2}\Big\vert \leq \sqrt{{\rm Var}(\max\limits_{i\leq n}\|\Row_i(M_n)\|_2)},
$$
we deduce that 
\begin{align*}
\Prob&\Big\{\Big\vert\frac{\w\rho_n'}{\w\rho_n}-1\Big\vert \geq \varepsilon\Big\}\\
&\leq \Prob\Big\{ \Big\vert \max\limits_{i\leq n}\|\Row_i(M_n)\|_2-\E \max\limits_{i\leq n}\|\Row_i(M_n)\|_2\Big\vert \geq \varepsilon \w\rho_n/2
-\sqrt{{\rm Var}(\max\limits_{i\leq n}\|\Row_i(M_n)\|_2)}\Big\}.
\end{align*}
Finally, note that $\max\limits_{i\leq n}\|\Row_i(M_n)\|_2$ is a convex $1$-Lipschitz function of i.i.d random variables uniformly bounded by $\sqrt{h}$. Therefore, by Talagrand's inequality (see Theorem~\ref{th: tal}), we have ${\rm Var}(\max\limits_{i\leq n}\|\Row_i(M_n)\|_2)\leq C\, h$ for some appropriate constant $C$. Using that $\w\rho_n\geq 2\sqrt{d}$ and that
$ d\underset{n\to \infty}{\longrightarrow} \infty$, we get that for $n$ large enough 
$$
\Prob\Big\{\Big\vert\frac{\rho_n'}{\w\rho_n}-1\Big\vert \geq \varepsilon\Big\} 
\leq \Prob\Big\{ \Big\vert \max\limits_{i\leq n}\|\Row_i(M_n)\|_2-\E \max\limits_{i\leq n}\|\Row_i(M_n)\|_2\Big\vert \geq
\varepsilon \sqrt{d}/4\Big\}.
$$
It remains to apply Talagrand's inequality again to finish the proof. 
\end{proof}

While in the regime
$\frac{np_n}{\log n}\to \infty$ we have $\frac{\rho_n}{2\sqrt{np_n}}\to 1$
thanks to concentration, the situation is completely different in the sparse regime 
$\lim \frac{np_n}{\log n}\to 0$. The lack of strong concentration implies existence 
of a row with squared Euclidean norm significantly above the average $np_n$.
This, in turn, implies that $\rho_n$ is significanly larger than $2\sqrt{np_n}$. We will formally verify this fact in the next lemma. 

\begin{lemma}\label{l: sublog convergence}
Let $\xi$ be a real uniformly bounded random variable with $\E\xi^2=1$, and for each $n$, let
$M_n$ be an $n\times n$ symmetric random matrix with zero diagonal and i.i.d.\ entries above the diagonal
equidistributed with the product $b_n\xi$, where $b_n$ is a Bernoulli random variable independent from $\xi$
and with probability of success $p_n$. Assume further that $\lim\limits_{n\to \infty} n p_n=\infty$
and that $\lim\limits_{n\to \infty} \frac{n p_n}{\log n}=0$.
Then 
$$\frac{\max\limits_{i\leq n}\|\Row_i(M_n)\|_2}{\sqrt{n p_n}}\underset{n\to\infty}{\overset{\Prob}{\longrightarrow}} \infty.$$
\end{lemma}
\begin{proof} 
We will assume that $\xi^2\leq h$ everywhere on the probability space.
Let $\alpha \geq 3$. Our goal is to show that 
$$
\lim\limits_{n\to \infty}\Prob\Big\{\max\limits_{i\leq n}\|\Row_i(M_n)\|_2\geq \alpha\sqrt{n p_n}\Big\} =1.
$$
We write $M_n$ as the Hadamard product of two independent symmetric random matrices
$B=(b_{ij})_{1\leq i,j\leq n}$ and $\Xi= (\xi_{ij})_{1\leq i,j\leq n}$ where $b_{ij}$ are equidistributed 
with $b_n$ and $\xi_{ij}$ equidistributed with $\xi$. Denote
\begin{align*}
\Event_1&:=\Big\{ \max\limits_{i\leq n}\|\Row_i(M_n)\|_2\geq \E\max\limits_{i\leq n}\|\Row_i(M_n)\|_2- \sqrt{np_n}\Big\},\\
\Event_2&:=\Big\{ \E\max\limits_{i\leq n}\sqrt{\deg(i)} \geq \max\limits_{i\leq n}\sqrt{\deg(i)}- \sqrt{np_n}\Big\},\\
\Event_3&:=\Big\{ \max\limits_{i\leq n}\deg(i) \geq 4\alpha^2 np_n\Big\},
\end{align*}
where $\deg(i)$ refers to the degree of vertex $i$ in the random graph with the adjacency matrix $B$. 
Since $\lim\limits_{n\to \infty} n p_n=\infty$, then it follows from Theorem~\ref{th: tal} that 
$
\lim\limits_{n\to \infty}\Prob\big(\Event_1\cap \Event_2\big)=1.
$
 Moreover, Theorem~\ref{th: tal} also implies that ${\rm Var}\big(\max\limits_{i\leq n}\|\Row_i(M_n)\|_2\big)\leq C(h)$, where $C(h)$ is a constant depending only on $h$. 
Using that $\xi$ has unit second moment together with Jensen's inequality, we can write  
 $$
\E\max\limits_{i\leq n}\sqrt{\deg(i)}\leq \sqrt{\E\max\limits_{i\leq n}\deg(i)}
\leq \sqrt{\E\max\limits_{i\leq n}\|\Row_i(M_n)\|_2^2}\leq  \E\max\limits_{i\leq n}\|\Row_i(M_n)\|_2 + \sqrt{C(h)}.
 $$ 
If $n$ is large enough, we have $ \sqrt{C(h)}\leq \sqrt{np_n}$, and the above implies that 
$$
\Event_1\cap\Event_2\cap\Event_3\subseteq \{\max\limits_{i\leq n}\|\Row_i(M_n)\|_2\geq \alpha\sqrt{n p_n}\Big\}.
$$
Since $\lim\limits_{n\to \infty}\Prob\big(\Event_1\cap \Event_2\big)=1$,
our remaining task is to show that $\lim\limits_{n\to \infty}\Prob(\Event_3)=1$. 
Estimates on the maximum degree of an Erd\H{o}s--Renyi graph are available in the literature. We will use the following estimate \cite[Theorem~3.1]{book-bollobas} asserting that if $\kappa:=\kappa(n)$ is an integer satisfying 
$$
\lim_{n\to \infty} n{n-1 \choose \kappa} p_n^{\kappa}(1-p_n)^{n-1-\kappa} =\infty,
$$
then we have
$$
\lim_{n\to \infty} \Prob\Big\{\max\limits_{i\leq n}\deg(i) \geq \kappa\Big\}=1
$$
It remains to check that we could apply the above fact with $\kappa= \lceil 4\alpha^2 np_n\rceil$. Using that $\lim\limits_{n\to \infty} \frac{n p_n}{\log n}=0$, an easy calculation finishes the proof. 
\end{proof}

The following lemma shows that the quantity $\max\limits_{1\leq i\leq n}\|\Row_i(W)\|_2$ 
is stable if we eliminate a fraction of the rows. This will help us affirm that for a fixed $k$, the $k$ largest rows are of the same order. 

\begin{lemma}\label{l: maxnorm-restrict}
Let $\xi$ be a uniformly bounded random variable with $\E\xi^2=1$. 
For each $n$, let $W_n=(w_{ij})$ be an $n\times n$ random symmetric matrix with independent
(up to the symmetry constraint) entries equidistributed with $b_n\xi$, where $b_n$ is Bernoulli ($0/1$)
random variable independent from $\xi$, with probability of success $p_n$ and assume that $np_n\to \infty$. Then for 
any $r\in [1/n,1]$ and any $\varepsilon\in (0,1)$, we have 
$$
\Prob\Big\{ \max\limits_{1\leq i\leq n}\|\Row_i(W_n)\|_2 \geq (1+\varepsilon) \max\limits_{1\leq i\leq \lfloor rn\rfloor}\|\Row_i(W_n)\|_2\Big\} 
\leq 5r^{-1} e^{-c\varepsilon^2 np_n},
$$ 
where $c>0$ may only depend on the distribution of $\xi$. 
In particular, for any fixed integer $k$, we have
$$\frac{\Vert \Row_k(W_n)\Vert_2^*}{\max\limits_{i\leq n}\|\Row_i(W_n)\|_2}\underset{n\to\infty}{\overset{\Prob}{\longrightarrow}} 1,$$
where we denoted by $\Vert \Row_k(W_n)\Vert_2^*$ the $k$-th largest element in the sequence $(\Vert \Row_i(W_n)\Vert_2)_{i\leq n}$. 
\end{lemma}
\begin{proof}
We will assume that $\xi^2\leq h$ for some number $h$, and that $n p_n$ is bounded from below by a large constant.
Let $r\in[1/n,1]$ and $\varepsilon\in (0,1/3]$. 
Since $\max\limits_{1\leq i\leq \lfloor r  n\rfloor}\|\Row_i(W_n)\|_2$
is a $1$--Lipschitz function of the entries of $W_n$, Theorem~\ref{th: tal} implies that 
$$
{\rm Var}\big(\max\limits_{1\leq i\leq \lfloor r  n\rfloor}\|\Row_i(W_n)\|_2\big)\leq C,
$$
where $C:=C(h)$ is a constant depending only on $h$. Therefore, since $\E\, \max\limits_{1\leq i\leq \lfloor r  n\rfloor}\|\Row_i(W_n)\|_2^2\geq np_n$ is large enough, we deduce that 
$\E\, \max\limits_{1\leq i\leq \lfloor r  n\rfloor}\|\Row_i(W_n)\|_2\geq \sqrt{np_n/2}$. This, together with Theorem~\ref{th: tal}, implies
$$
\Prob\Big\{ \max\limits_{1\leq i\leq \lfloor r  n\rfloor}\|\Row_i(W_n)\|_2
\leq \Big(1-\frac{\varepsilon}{2}\Big) \E\, \max\limits_{1\leq i\leq \lfloor r  n\rfloor}\|\Row_i(W_n)\|_2\Big\} \leq 
e^{-c\varepsilon^2 np_n},
$$
for some appropriate constant $c>0$. Therefore, we can write 
\begin{align*}
\Prob\Big\{ \max\limits_{1\leq i\leq  n}&\|\Row_i(W_n)\|_2\geq (1+\varepsilon) \max\limits_{1\leq i\leq \lfloor r  n\rfloor}\|\Row_i(W_n)\|_2\Big\} \\
&\leq \Prob\Big\{ \max\limits_{1\leq i\leq   n}\|\Row_i(W_n)\|_2\geq \Big(1+\frac{\varepsilon}{3}\Big) \E\, \max\limits_{1\leq i\leq \lfloor r  n\rfloor}\|\Row_i(W_n)\|_2\Big\}
+e^{-c\varepsilon^2 np_n}.
\end{align*}
Now, define random variables 
$$
\beta_u= \max\limits_{u< i\leq u+\lfloor r n\rfloor}\|\Row_i(W_n)\|_2,
\quad u=0,1,\dots,n-\lfloor r  n\rfloor,
$$
and note that these variables have the same distribution. 
Clearly, we can choose $\big\lceil n/\lfloor r n\rfloor\big\rceil$ indices $u$ such that
$\max\limits_{1\leq i\leq n}\|\Row_i(W_n)\|_2$ equals the maximum of $\beta_u$ over those indices.
Therefore, applying the union bound, we get 
\begin{align*}
\Prob\Big\{\max\limits_{1\leq i\leq n}\|\Row_i(W_n)\|_2\geq \Big(1+\frac{\varepsilon}{3}\Big)\E\, \max\limits_{1\leq i\leq \lfloor r  n\rfloor}\|\Row_i(W_n)\|_2\Big\}
&\leq \frac{4}{r}\,  \Prob\Big\{ \beta_1\geq \Big(1+\frac{\varepsilon}{3}\Big)\E\, \beta_1\Big\}.
\end{align*}
Using once again Theorem~\ref{th: tal} and putting together the previous estimates with the fact that $n$ is large enough, we finish the proof of the first part. 

To prove the last claim, consider the sets $I_s:=\{(s-1) \lfloor n/k\rfloor+1,\ldots, s \lfloor n/k\rfloor\}$ for any $s\leq k-1$ and $I_k:= \{(k-1) \lfloor n/k\rfloor+1,\ldots,n\}$. Note that, by the above, we have 
$$\lim\limits_{n\to\infty}\Prob\Big\{\min_{s\leq k}\max\limits_{i\in I_s}\|\Row_i(W_n)\|_2\geq (1-\varepsilon)\max\limits_{i\leq n}\|\Row_i(W_n)\|_2\Big\}=1\quad \mbox{for all }\varepsilon>0.$$  
Therefore, we deduce that 
$$\lim\limits_{n\to\infty}\Prob\Big\{\Vert \Row_k(W_n)\Vert_2^*\geq (1-\varepsilon)\max\limits_{i\leq n}\|\Row_i(W_n)\|_2\Big\}=1\quad \mbox{for all }\varepsilon>0.$$  
It remains to note that we always have $\Vert \Row_k(W_n)\Vert_2^*\leq \max\limits_{i\leq n}\|\Row_i(W_n)\|_2$ to finish the proof. 
\end{proof}

In the previous sections, we restricted our attention to the regime when $np_n\geq \log^{c} n$ for some $c\in(0,1)$. The reason is that when the matrix is very sparse ($\lim_n \frac{np_n}{\log n}=0$), one can use available results \cite{BGBK sparse} where an upper bound is provided for the norm of the centered adjacency matrix of an Erd\H{o}s--Renyi graph. 
Since our model is slightly different, we indicate the necessary changes to similarly obtain an adequate bound. We summarize this in the next proposition.

\begin{prop}\label{prop: bound-very sparse}
Let $\xi$ be a real centered uniformly bounded random variable of unit variance. 
For each $n$, let $W_n$ be an $n\times n$ symmetric random matrix with i.i.d.\ entries (up to the symmetry constraint),
with each entry equidistributed with the product $b_n\xi$, where $b_n$ is $0/1$ (Bernoulli) random variable independent of $\xi$,
with probability of success equal to $p_n$. Assume further that $n p_n\to\infty$ and $\frac{np_n}{\log n}\to 0$ with $n$. 
Then for any fixed integer $k\geq 1$, we have
$$\frac{\vert \lambda_{\vert k\vert}(W_n)\vert}{\max\limits_{i\leq n}\|\Row_i(W_n)\|_2}\underset{n\to\infty}{\overset{\Prob}{\longrightarrow}} 1,$$
where $\lambda_{\vert k\vert}(W_n)$ denotes the $k$-th largest (in absolute value) eigenvalue of $W_n$. 
\end{prop}
\begin{proof}
We will assume that $\xi^2\leq h$ for some $h\geq 1$.
By Lemmas~\ref{lem: upper-zero-diag} and \ref{l: sublog convergence}, we may (and will) suppose that $W_n$ has zero diagonal. 
We write $W_n:= \Xi_n\bullet A_n$ as the Hadamard product of two independent symmetric random matrices (with zero diagonals) $A_n=(a_{ij})_{1\leq i,j\leq n}$ and $\Xi_n= (\xi_{ij})_{1\leq i,j\leq n}$ where $a_{ij}$ are equidistributed with $b_n$ and $\xi_{ij}$ equidistributed with $\xi$ (for $i\neq n$). Let us denote $\w W_n:= \Xi_n \bullet (A_n-\E\, A_n)$. 

Fix an integer $k$. By Weyl's perturbation inequality, we have 
$$
\Big\vert\,  \vert \lambda_{\vert k\vert}(W_n)\vert- \vert \lambda_{\vert k\vert}(\w W_n)\vert\, \Big\vert\leq  p_n\Vert \Xi_n\Vert,
$$
where we have used that $\E\, A_n$ is the $n\times n$ zero diagonal matrix having all its non zero entries equal to $p_n$. 
It is known that $\lim_{n\to\infty}\Prob\big\{ \Vert \Xi_n\Vert \leq C \sqrt{n}\big\} =1$, where $C:=C(h)$ is a constant depending only on $h$ (see for example \cite[Corollary 4.4.8]{vershynin-book}). Using this together with Lemma~\ref{l: sublog convergence}, we deduce that in order to obtain the statement of the proposition, it is sufficient to prove that 
$$\frac{\vert \lambda_{\vert k\vert}(\w W_n)\vert}{\max\limits_{i\leq n}\|\Row_i(W_n)\|_2}\underset{n\to\infty}{\overset{\Prob}{\longrightarrow}} 1.$$

We follow the strategy of the proof of \cite[Theorem~1.2]{BGBK sparse} (see page 15 there). 
Let $0<\delta<(2h)^{-1}$ be fixed, and set $t=\delta L_1$ with 
$$
L_1:=\frac{\log n}{\log \big((\log n)/(np_n)\big)}.
$$
We denote by $G$ the graph with adjacency matrix $A_n$ and note that the largest degree in $G$ is $(1+o(1))L_1$ with
probability going to one with $n$ (see \cite[Corollary~1.13]{BGBK sparse}). Define $G_{\star}$ as the subgraph of
$G$ with the vertex set $[n]$
obtained by keeping edges $i\edg j$ whenever $i\in \mathcal{V}_{\geq t}$ and $j\not\in \mathcal{V}_{\geq t} \cup \mathcal{N}_G\big(\mathcal{V}_{\geq t}\setminus \{i\}\big)$, where 
$$
\mathcal{V}_{\geq t}:=\{v\in [n]:\, \deg_G(v)\geq t\}
$$
and
$$
\mathcal{N}_G\big(\mathcal{V}_{\geq t}\setminus \{i\}\big):=\{v\in [n]:\, \exists u\in \mathcal{V}_{\geq t}\setminus \{i\} \text{ such that } v\edg u\}.
$$
Combining Lemma~\ref{l: maxnorm-restrict}, Talagrand's inequality and that $\E \max\limits_{i\leq n}\|\Row_i(W_n)\|_2^2 \geq \E \max\limits_{i\leq n} \degree_G(i)$, and using that $L_1/(2h)\geq \delta L_1$, we deduce that for any fixed $k$,
the random index $i_k$ of the row with the $k$-th largest Euclidean norm in $W_n$
belongs to $\mathcal{V}_{\geq t}$ with probability going to one with $n$. 
It follows from \cite[Lemma~2.5]{BGBK sparse} that with probability going to $1$ with $n$, we have for any $i\in \mathcal{V}_{\geq t}$ 
$$
\Big\vert\mathcal{N}_G(\{i\})\cap \big(\mathcal{V}_{\geq t} \cup \mathcal{N}_G\big(\mathcal{V}_{\geq t}\setminus \{i\}\big)\big)\Big\vert \leq \frac{c}{\delta},
$$
for some constant $c$. Therefore, with probability going to $1$ with $n$, we have that for any $i\in \mathcal{V}_{\geq t}$ 
$$
\Big\vert \Vert \Row_i(W_n)\Vert_2- \Vert \Row_i(\Xi_n\bullet A_{\star})\Vert_2\Big\vert \leq \frac{c\, \sqrt{h}}{\delta},
$$
where we denoted $A_{\star}$ the adjacency matrix of $G_{\star}$. In view of the above, with probability going to one with $n$, we have 
\begin{equation}\label{eq: block-diag}
\Big\vert \Vert \Row_k(W_n)\Vert_2^{*}- \Vert \Row_k(\Xi_n\bullet A_{\star})\Vert_2^*\Big\vert \leq \frac{c\, \sqrt{h}}{\delta},
\end{equation}
where we denoted by $\Vert \Row_k(W_n)\Vert_2^*$ (resp. $\Vert \Row_k(\Xi_n\bullet A_{\star})\Vert_2^*$) the $k$-th largest element in the sequence $(\Vert \Row_i(W_n)\Vert_2)_{i\leq n}$ (resp. $(\Vert \Row_i(\Xi_n\bullet A_{\star})\Vert_2)_{i\leq n}$). 
 Note that by construction, $\Xi_n \bullet A_{\star}$ is formed (up to a permutation)
by disjoint block diagonal matrices where each block diagonal matrix has only its first row and column non-zero.
Each such submatrix has two opposite non-zero eigenvalues, whose absolute value is equal to the Euclidean norm of its non-zero row. 
Using this, \eqref{eq: block-diag} and 
Lemma~\ref{l: maxnorm-restrict}, we deduce that 
$$\frac{\vert \lambda_{\vert k\vert}(\Xi_n \bullet A_{\star})\vert}{\max\limits_{i\leq n}\|\Row_i(W_n)\|_2}\underset{n\to\infty}{\overset{\Prob}{\longrightarrow}} 1.$$

In view of the Weyl perturbation inequality, the remaining task is to show that for any $\varepsilon >0$ 
\begin{equation}\label{eq: bound-very sparse}
\lim\limits_{n\to\infty}\Prob\big\{\Vert \Xi_n\bullet (A'_n-\E\, A_n) \Vert \leq \varepsilon \max\limits_{i\leq n}\|\Row_i(W_n)\|_2\big\}=1,
\end{equation}
where we denoted $A'_n=A_n-A_{\star}$. To this aim, as in \cite{BGBK sparse}, we will make use of the results in \cite{LLV}. First, we note that by \cite[Proposition~1.11]{BGBK sparse}, for all large $n$
the cardinality of $\mathcal{V}_{\geq t}$ is at most $10/p_n$
with probability going to one.
Moreover, with probability going to one with $n$, all vertices in the graph generated by $A'_n$ have degrees bounded by $t$. If 
our matrix $\Xi_n$ was a matrix of all ones then, applying directly \cite[Theorem~2.1]{LLV} like it is stated in \cite{LLV}, we would get 
$$
\lim_{n\to\infty}\Prob\Big\{\Vert A'_n-\E\, A_n\Vert \leq C' \big(\sqrt{np_n}+\sqrt{\delta L_1}\big)\Big\}=1.
$$
We need a weighted version of this result in our setting to obtain a similar bound for $\Xi_n\bullet (A'_n-\E\, A_n)$.
The proof of \cite[Theorem~2.1]{LLV} relies on a special decomposition of the Erd\H{o}s--Renyi graphs (see \cite[Theorem~2.6]{LLV}) which can be modified 
as to serve our needs. Namely, the first conclusion of \cite[Theorem~2.6]{LLV} stating that the adjacency matrix of the Erd\H{o}s--Renyi graph concentrates well, can be replaced by its weighted version i.e. $\Vert \Xi\bullet (A-\E A)_{\mathcal{N}}\Vert$ following their notations. Indeed, to view this, it can be checked that the same decomposition procedure can be carried over the graph and one only needs to update \cite[Lemma~3.3]{LLV} as to allow a weighted version of it. This, in turn, can be easily checked 
by carrying almost the same proof, and modifying the variables $X_i$ appearing in formula (3.3) in \cite{LLV} by introducing $\xi_{ij}$ in the corresponding sum. The rest of that proof follows the same lines by using Bernstein's inequality which produces the same bounds up to a constant depending only on $h$. 

In view of this, one gets 
$$
\lim_{n\to\infty}\Prob\Big\{\Vert \Xi_n\bullet (A'_n-\E\, A_n)\Vert \leq C'' \big(\sqrt{np_n}+\sqrt{\delta L_1}\big)\Big\}=1,
$$
where $C'':=C''(h)$ is a constant depending only on $h$. Now using \cite[Corollary~1.13]{BGBK sparse}, we can replace $L_1$ by $\max_i \deg_G(i)$ in the above expression to obtain 
$$
\lim_{n\to\infty}\Prob\Big\{\Vert \Xi_n\bullet (A'_n-\E\, A_n)\Vert \leq \w C \big(\sqrt{np_n}+\sqrt{\delta \max_i \deg_G(i)}\big)\Big\}=1,
$$
for an appropriate constant $\w C$ depending only on $h$. Using Talagrand's inequality (Theorem~\ref{th: tal}) similarly to what is done with the events $\Event_1$ and $\Event_2$ in the previous lemma, we deduce that 
$$
\lim_{n\to\infty}\Prob\Big\{\Vert \Xi_n\bullet (A'_n-\E\, A_n)\Vert \leq \bar{C} \big(\sqrt{np_n}+\sqrt{\delta}\max\limits_{i\leq n}\|\Row_i( W_n)\|_2\big)\Big\}=1,
$$
for an appropriate constant $\w C$ depending only on $h$. It remains to use Lemma~\ref{l: sublog convergence} and choose $\delta$ appropriately to deduce \eqref{eq: bound-very sparse} and finish the proof. 
\end{proof}

We are now ready to state and prove an upper bound on the operator norm. 
The main statement of this section is the following. 

\begin{theorem}\label{th: centered upper}
Let $\xi$ be a real centered uniformly bounded random variable of unit variance.
For each $n$, let $W_n$ be $n\times n$ symmetric random matrix with i.i.d.\ entries above and on the main diagonal,
with each entry equidistributed with the product $b_n\xi$, where $b_n$ is $0/1$ (Bernoulli) random variable independent of $\xi$,
with probability of success equal to $p_n$. Assume further that $n p_n\to\infty$ with $n$ and  denote 
$$\rho_n:=\theta_n+\frac{n p_n}{\theta_n},\quad \theta_n:=\sqrt{\max\big(\max\limits_{i\leq n}\|\Row_i(W_n)\|_2^2-np_n,n p_n\big)},$$
Then for any $\varepsilon>0$ we have
$$\lim\limits_{n\to\infty}\Prob\big\{\Vert W_n\Vert \leq (1+\varepsilon)\rho_n\big\}=1.$$
\end{theorem}
\begin{proof}
Let us first note that when $\lim\limits_n \frac{np_n}{\log n}=\infty$, then standard concentration estimates show that
$\frac{\rho_n}{2\sqrt{np_n}}$ converges to one in probability as $n$ goes to infinity. 
On the other hand, known results (see \cite{Khorunzhy,BGBK,LVY}) imply that in this regime
$$
\lim\limits_{n\to\infty}\Prob\big\{\Vert W_n\Vert \leq (1+\varepsilon)2\sqrt{n p_n}\big\}=1\quad \mbox{ for any $\varepsilon>0$.} 
$$

Further, note that deterministically $\rho_n\geq \max\limits_{i\leq n}\|\Row_i(W_n)\|_2$,
so in view of Proposition~\ref{prop: bound-very sparse}, we get 
$$
\lim\limits_{n\to\infty}\Prob\big\{\Vert W_n\Vert \leq (1+\varepsilon)\rho_n\big\}=1\quad\mbox{ for any $\varepsilon>0$},
$$
whenever $\lim\limits_{n}\frac{np_n}{\log n}=0$.

In view of the above remarks, we can (and will) assume that $\log^{\frac{1}{1+\delta/2}} n\leq np_n\leq \log^2 n$,
where $\delta:=\min(c_{\text{\tiny\ref{cor: contribution-few multone}}}/2,1/100)$
(note that we could use a much stronger assumption $c\log n\leq n p_n\leq C\log n$ for all $n$,
but prefer to work under weaker conditions, which show that our argument developed in the previous sections,
covers a wider range of parameters). Assume that $\xi^2\leq h$ for some $h\geq 2$ everywhere on the probability space.

An approximation argument shows that for every $n$ there is a random variable $\xi_n$ with an absolutely continuous distribution,
of zero mean and unit variance and bounded by the absolute value by $\sqrt{h}$ with the following property:
first, denoting by $\w W_n$ the random $n\times n$ symmetric matrix with i.i.d.\ random variables above and on the main diagonal
equidistributed with $\xi_n b_n$, where $b_n$ is independent from $\xi_n$, we have that $\frac{\|\w W_n\|}{\|W_n\|}$
converges to one in probability; second, the ratio 
$$\frac{\max(\max_{i}\|\Row_i(\w W_n)\|_2^2,2n p_n)}
{\max(\max_{i}\|\Row_i(\w W_n)\|_2^2- np_n,np_n)^{1/2}}\;
\Big\slash\;
\frac{\max(\max_{i}\|\Row_i(W_n)\|_2^2,2n p_n)}
{\max(\max_{i}\|\Row_i(W_n)\|_2^2- np_n,np_n)^{1/2}}$$
converges to one in probability.
Thus, in our proof we can ``replace'' the matrices $W_n$ with $\w W_n$
and quantities $\rho_n$ -- with $\frac{\max(\max_{i}\|\Row_i(\w W_n)\|_2^2,2n p_n)}
{\max(\max_{i}\|\Row_i(\w W_n)\|_2^2- np_n,np_n)^{1/2}}$.

Let $\varepsilon \in (0,1)$, $k=k(n):=\lceil \log n\,\log\log n\rceil$ and let $n_0:= n_0(\varepsilon)$ be large enough. 
Assume that $n\geq n_0$.
Let us denote $M_n:=\widetilde W_n-{\rm Diag}(\widetilde W_n)$, $d=d(n):=(n-1)p_n$ and 
$$\w\rho_n':=\w\theta_n'+\frac{d}{\w\theta_n'},\quad \w\theta_n'
:=\sqrt{\max\big(\E\max\limits_{i\leq n}\|\Row_i(M_n)\|_2^2-d,d\big)}.$$
In view of Lemma~\ref{lem: upper-zero-diag}, it is sufficient to prove that 
$$\lim\limits_{n\to\infty}\Prob\big\{\Vert M_n\Vert \leq (1+\varepsilon)\w\rho_n'\big\}=1.$$
We will define parameters $\maxnorm$ and $\dmax$ the same way as in \eqref{eq: maxnorm act def}
and \eqref{eq: dmax act def}, respectively:
$$
\maxnorm:=\big(1+(\log\log\log n)^{-1}\big)^3\, \E\max\limits_{i\leq n}\sum\limits_{j=1}^n \mu_{ij}^2,
$$
where $\mu_{ij}$ are entries of $M_n$,
and 
$$
\dmax:=\big(1+(\log\log\log n)^{-1}\big)\, \E\max\limits_{i\leq n}\sum\limits_{j=1}^n b_{ij}. 
$$ 
Let us note that $d\leq \dmax\leq \maxnorm\leq 2h\dmax$. 
Moreover let $\stnd$ be defined as in \eqref{eq: definition-Y} and note that by \eqref{eq: bound-Y} we have 
$$
d-h\leq \Vert \stnd\Vert_1\leq \maxnorm,
$$
and define events 
$$\Event_{\maxnorm}:=\Big\{\sum\limits_{j=1}^n \mu_{ij}^2\leq \frac{\maxnorm}{1+(\log\log\log n)^{-1}}\mbox{ for all $i\in[n]$}\Big\},$$
and
$$\Event_g:=\big\{G \text{ is $(k/\log^2d)$-tangle free and }\deg_i(G)\leq d_{\max}\mbox{ for all $i\in[n]$}\big\},$$
where $G$ is the random graph on $[n]$ with the adjacency matrix $(b_{ij})$. Finally let 
\begin{align*}
\Event_{mjr}:=\Big\{
&\forall\,i\leq n,\;\mbox{the vector $(b_{ij})_{j=1}^n$ has at most $d^{1+\delta}$ non-zero components
AND}\\
&\mbox{for any vertex $v\in[n]$ the number of its heavy neighbors is at most $d^{\frac89}$}\Big\}.
\end{align*}
We start by writing 
$$
\Prob\big\{\Vert M_n\Vert \geq (1+\varepsilon)\w\rho_n'\big\}\leq 
\Prob\big\{\Vert M_n\Vert \mathbf{1}_{\Event_{\maxnorm}\cap \Event_g\cap \Event_{mjr}}
\geq (1+\varepsilon)\w\rho_n'\big\}+ \Prob\big\{(\Event_{\maxnorm}\cap \Event_g\cap \Event_{mjr})^c\big\}.
$$
We will use Markov's inequality in order to estimate the first term above. We have for any $k\geq 1$:
$$
\Prob\big\{\Vert M_n\Vert \mathbf{1}_{\Event_{\maxnorm}\cap \Event_g\cap \Event_{mjr}}\geq (1+\varepsilon)\w\rho_n'\big\}\leq 
\frac{\E [\Vert M\Vert^{2k} \mathbf{1}_{\Event_{\maxnorm}\cap \Event_g\cap \Event_{mjr}}]}{(1+\varepsilon)^{2k}{(\w\rho_n')}^{2k}},
$$
where $\mu_{ij}$ are the entries of $M_n$.
Now using that $\Vert M\Vert^{2k}\leq {\rm Tr} (M^{2k})$ and expressing the trace of the $2k$-th power in terms of the entries if the matrix, we deduce that 
$$
\Prob\big\{\Vert M_n\Vert \mathbf{1}_{\Event_{\maxnorm}\cap \Event_g\cap \Event_{mjr}}\geq (1+\varepsilon)\w\rho_n'\big\}\leq 
\frac{\E\Big(\sum\limits_{\path}\, 
\prod\limits_{\ell=1}^{2k}\mu_{\path(\ell-1),\path(\ell)}\mathbf{1}_{\Event_{\maxnorm}\cap\Event_g\cap \Event_{mjr}}\Big)}{(1+\varepsilon)^{2k}(\w\rho_n')^{2k}},
$$
where the summation is taken over all closed paths on $K_{[n]}$ of length $2k$.
Further, we write
\begin{align*}
\E\Big(&\sum\limits_{\path}\, 
\prod\limits_{\ell=1}^{2k}\mu_{\path(\ell-1),\path(\ell)}\mathbf{1}_{\Event_{\maxnorm}\cap\Event_g\cap \Event_{mjr}}\Big)\\
&\leq \E\Big(\sum\limits_{\multone(\path)< k/\log^{11}k}\, 
\prod\limits_{\ell=1}^{2k}\mu_{\path(\ell-1),\path(\ell)}\mathbf{1}_{\Event_{\maxnorm}\cap\Event_g\cap \Event_{mjr}}\Big)\\
&+\E\Big(\sum\limits_{\multone(\path)\geq k/\log^{11}k}\, 
\prod\limits_{\ell=1}^{2k}\mu_{\path(\ell-1),\path(\ell)}\mathbf{1}_{\Event_{\maxnorm}\cap\Event_g}\Big)
+\E\Big(\sum\limits_{\multone(\path)\geq k/\log^{11}k}\, 
\prod\limits_{\ell=1}^{2k}|\mu_{\path(\ell-1),\path(\ell)}|\,\mathbf{1}_{\Event_{mjr}^c}\Big).
\end{align*}

Since $\maxnorm\geq d\gg \log^{2/3} n$, then choosing $n_0$ large enough, we get, in view of Proposition~\ref{prop: mult-one},
\begin{align*}
\E\Big(\sum\limits_{\multone(\path)\geq k/\log^{11}k}\, 
\prod\limits_{\ell=1}^{2k}\mu_{\path(\ell-1),\path(\ell)}\mathbf{1}_{\Event_{\maxnorm}\cap\Event_g}\Big)
\leq 
n
\end{align*}
for any $n\geq n_0$.
On the other hand, since $d\geq \log^{\frac{1}{1+\delta/2}}n$, we can apply Proposition~\ref{prop: heavyn} to get
$
\Prob(\Event_{mjr}^c)\leq \exp(-c d^{1+\delta}),
$
whence
$$
\E\Big(\sum\limits_{\multone(\path)\geq k/\log^{11}k}\, 
\prod\limits_{\ell=1}^{2k}|\mu_{\path(\ell-1),\path(\ell)}|\,\mathbf{1}_{\Event_{mjr}^c}\Big)
\leq n(\dmax)^{2k}h^k\exp(-c d^{1+\delta}),
$$
where we used a trivial bound on the path weights. 
In view of our assumption on $k$,
$$
\lim_{n\to\infty} \Prob\big\{\Vert M_n\Vert \mathbf{1}_{\Event_{\maxnorm}\cap \Event_g\cap \Event_{mjr}}\geq (1+\varepsilon)\w\rho_n'\big\}
\leq 
\lim_{n\to \infty} \frac{\E\Big(\sum\limits_{\multone(\path)< \frac{k}{\log^{11}k}}\, 
\prod\limits_{\ell=1}^{2k}\mu_{\path(\ell-1),\path(\ell)}
\mathbf{1}_{\Event_{\maxnorm}\cap\Event_g\cap \Event_{mjr}}\Big)}{(1+\varepsilon)^{2k}(\w\rho_n')^{2k}}.
$$
Note that almost everywhere on the event $\Event_{\maxnorm}\cap\Event_g\cap \Event_{mjr}$
the matrix $M_n$ belongs to the set $\matrixset(n,k/\log^2d,d,\dmax, \maxnorm, \stnd)$
defined in \eqref{eq: definition-matrixset} (in particular, we use here that the distribution of $\xi_n$ is absolutely continuous);
moreover, the parameters $d,\dmax,\maxnorm,h$ and $\|\stnd\|_1$ satisfy \eqref{eq: condition-d-dmax}.
Since on this event we also have
$$
\frac{\maxnorm}{\Vert \stnd\Vert_1}\leq \frac{h\dmax}{d-h}\leq \frac{hd^{1+\delta}}{d-h}\leq d^{c_{\text{\tiny\ref{cor: contribution-few multone}}}},
$$
then we can apply Corollary~\ref{cor: contribution-few multone} to get 
$$
\lim_{n\to\infty} \Prob\big\{\Vert M_n\Vert \mathbf{1}_{\Event_{\maxnorm}\cap \Event_g\cap \Event_{mjr}}\geq (1+\varepsilon)\w\rho_n'\big\}\leq 
\lim_{n\to \infty} \frac{n\, f(\maxnorm, \Vert \stnd\Vert_1)^{2k}\, e^{C\frac{k}{\log k}}}{(1+\varepsilon)^{2k}(\w\rho_n')^{2k}},
$$
where we denoted $f(x,y)= \frac{\max(x,2y)}{\sqrt{\max(x-y,y)}}$. Now using that $f$ is increasing
coordinate-wise and that $\Vert \stnd\Vert_1\leq \big(1+(\log\log\log n)^{-1}\big) d$ by \eqref{eq: bound-Y}, we deduce that 
$$
f(\maxnorm, \Vert \stnd\Vert_1)\leq \big(1+(\log\log\log n)^{-1}\big) \w\rho_n', 
$$
and thus 
$$
\lim_{n\to\infty} \Prob\big\{\Vert M_n\Vert \mathbf{1}_{\Event_{\maxnorm}\cap \Event_g\cap \Event_{mjr}}\geq (1+\varepsilon)\w\rho_n'\big\}\leq 
\lim_{n\to \infty} \frac{n\, \big(1+(\log\log\log n)^{-1}\big)^{2k}\, e^{C\frac{k}{\log k}}}{(1+\varepsilon)^{2k}}.
$$
In view of the conditions on $k$, the above probability tends to zero. 

To finish the proof, it remains to show that 
$$
\lim_{n\to\infty}\Prob\big\{(\Event_{\maxnorm}\cap \Event_g\cap \Event_{mjr})^c\big\}=0
$$
This follows by combining Proposition~\ref{prop: heavyn}, Lemma~\ref{lem: tangle-free-prob} and Theorem~\ref{th: tal}.
\end{proof}

\section{Non-centered matrices}

In this section, we prove a non-symmetric version of Theorem~\ref{th: centered upper}
from the previous section, which would provide
upper bounds for the second largest, by absolute value, eigenvalue of a random symmetric matrix
with non-centered entries.
In particular, this will allow us to identify the necessary
and sufficient conditions for the presence of non-trivial outliers in the spectrum of adjacency matrices of the Erd\H os--Renyi graphs.
The reduction of the main statement of this section (Theorem~\ref{th: non-centered upper} below) to Theorem~\ref{th: centered upper}
is done by means of a special coupling between sequences of symmetric centered and non-centered random matrices,
combined with Talagrand's inequality for product measures.

Given a non-centered (sparse) symmetric random matrix $A$ with i.i.d.\ entries and denoting by $\lambda_{|2|}(A)$
the second largest (by absolute value) eigenvalue of $A$, we clearly have
$$|\lambda_{|2|}(A)|\leq\|A-\E A\|,$$
where $A-\E A$ is a symmetric centered random matrix.
However, in our context this trivial symmetrization is of no use: unlike the original matrix $A$,
the matrix $A-\E A$ is non-sparse, and a direct application of Theorem~\ref{th: centered upper}
is not possible. On the other hand, results existing in the literature (such as \cite{BGBK})
do not give sufficiently strong estimates for $\|A-\E A\|$ since the entries of the matrix are very spiky.
The basic idea which we employ to obtain the required result is to ``replace'' the matrix $A-\E A$
with a sparse centered matrix of the form $A-B\bullet \E A$,
where ``$\bullet$'' denotes entry-wise matrix product and
$B$ is an appropriately rescaled symmetric matrix with i.i.d.\ Bernoulli ($0/1$) entries
such that $\E B={\bf 1}{\bf 1}^\top$ (in fact, our definition will be slightly different,
although quite close to this one). In this case, the problem lies
in finding a relation between the new centered matrix and $\lambda_{|2|}(A)$.

The main technical statement of the section is Proposition~\ref{prop: coupling} defining the coupling.
For better readability, we extract a part of its proof into the two following lemmas.
\begin{lemma}\label{l: entrywise coupling}
Let $\varepsilon,p\in(0,1)$ with $p/ \varepsilon\leq 1$.
Let $\xi$ be a random variable of unit second moment,
and let $a$ and
$b'$ be Bernoulli ($0/1$) random variables with probabilities of success $ \varepsilon$ and
$p/ \varepsilon$, respectively, such that $\xi,a,b'$ are jointly independent.
Define
$$\xi_\varepsilon:=a\,\xi-\frac{ \varepsilon(1-a)\E\xi}{1- \varepsilon}\quad \text{ and } \quad \xi_\varepsilon':=\frac{\xi_{\varepsilon}}
{\sqrt{\Var\big(\xi_{\varepsilon}\big)}},$$
and set $b:=a\,b'$.
Then
\begin{itemize}

\item[(a)] $\xi_\varepsilon'$ is of zero mean and unit variance;

\item[(b)] $b$ is Bernoulli with probability of success $p$;

\item[(c)] Let $(W,W')$ be a pair of symmetric random matrices 
such that the collection of the pairs of entries $\{(w_{ij}, w_{ij}'), j\geq i\}$ from $(W,W')$ are i.i.d and   equidistributed with $(b \xi,b'\xi_{\varepsilon}')$. Then
$$
\E\max\limits_{i\leq n}\|\Row_i(W')\|_2^2\leq 
\beta^{-2} \bigg(1+\frac{\sqrt{\varepsilon}\, \vert \E\, \xi\vert}{1-\varepsilon}\bigg)^2\,\E\max\limits_{i\leq n}\|\Row_i(W)\|_2^2,
$$
where $\beta:=\sqrt{\Var\big(\xi_{\varepsilon}\big)}
=\sqrt{\varepsilon+\frac{ \varepsilon^{\,2}(\E\xi)^2}{1- \varepsilon}}$.
\end{itemize}

\end{lemma}
\begin{proof}
The first two assertions can be easily verified. 
Next, we consider property (c) of the lemma.
Let $W=(w_{ij})$ and $W'=(w_{ij}')$ be as stated above 
and  denote by $b_{ij},b'_{ij},a_{ij},\xi_{ij},(\xi_{\varepsilon}')_{ij}$
the variables associated with the couple $(w_{ij},w_{ij}')$.
Define an auxiliary $n\times n$ random matrix
$$\widetilde W:=(w_{ij}'{\bf 1}_{\{a_{ij}=0\}}).$$
Then, using the convexity of $\Vert\cdot\Vert_2^2$,   we get
$$
\max\limits_{i\leq n}\|\Row_i(W')\|_2^2
\leq (1+x^{-1})\max\limits_{i\leq n}\|\Row_i(\widetilde W)\|_2^2+
 (1+x)\max\limits_{i\leq n}\|\Row_i(W'-\widetilde W)\|_2^2,
$$
where $x= \sqrt{\varepsilon}\,  \vert \E \xi\vert/(1-\varepsilon)$. 
Observe that
$$w'_{ij}{\bf 1}_{\{a_{ij}=1\}}=b'_{ij}(\xi_{\varepsilon}')_{ij}{\bf 1}_{\{a_{ij}=1\}}
=\frac{b_{ij}\xi_{ij}{\bf 1}_{\{a_{ij}=1\}}}{\beta}=\frac{w_{ij}{\bf 1}_{\{a_{ij}=1\}}}{\beta}
=\frac{w_{ij}}{\beta},\quad i,j=1,\dots,n,$$
whence
$$
\max\limits_{i\leq n}\|\Row_i(W'-\widetilde W)\|_2^2= \beta^{-2}\max\limits_{i\leq n}\|\Row_i(W)\|_2^2.
$$
Further, since $\w w_{ij}= - \varepsilon\big((1-\varepsilon)\beta\big)^{-1} b_{ij}' (\E\, \xi_{ij})\, \mathbf{1}_{a_{ij}=0}$, then we have
\begin{equation}\label{eq: aux 309250235}
\max\limits_{i\leq n}\|\Row_i(\widetilde W)\|_2^2
\leq \bigg(\frac{ \varepsilon\,\E\xi}{(1- \varepsilon)\beta}\bigg)^2
\max\limits_{i\leq n}\sum\limits_{j=1}^n b_{ij}'.
\end{equation}
On the other hand,
$$\E\big(b_{ij}^2\xi_{ij}^2\;|\;b_{ij}'=1\big)=
 \varepsilon,$$
whence by Jensen's inequality 
$$
\E\big(\max\limits_{i\leq n}\|\Row_i(W)\|_2^2\;|\;b_{ij}',\;i,j=1,\dots,n\big)\geq
 \varepsilon\,\max\limits_{i\leq n}\sum\limits_{j=1}^n b_{ij}'.
$$
Combining this with relation \eqref{eq: aux 309250235}, we obtain
$$
\E\,\max\limits_{i\leq n}\|\Row_i(W)\|_2^2
\geq  \varepsilon\,\bigg(\frac{ \varepsilon\,\E\xi}{(1- \varepsilon)\beta}\bigg)^{-2}
\E\,\max\limits_{i\leq n}\|\Row_i(\widetilde W)\|_2^2.
$$
Thus,
\begin{align*}
\E\max\limits_{i\leq n}\|\Row_i(W')\|_2^2&\leq
(1+x^{-1}) \E\max\limits_{i\leq n}\|\Row_i(\widetilde W)\|_2^2+
(1+x)\E\max\limits_{i\leq n}\|\Row_i(W'-\widetilde W)\|_2^2\\
&\leq 
\beta^{-2} \bigg(1+\frac{\sqrt{\varepsilon}\, \vert \E\, \xi\vert}{1-\varepsilon}\bigg)^2\, \E\max\limits_{i\leq n}\|\Row_i(W)\|_2^2.
\end{align*}
The result follows.
\end{proof}
\begin{lemma}\label{l: coupling 2}
Let the variable $\xi$ and the matrices $W,W'$ be as in the last lemma,
and assume additionally that $|\xi|$ is uniformly bounded and that $p\leq \varepsilon$. Further,
let $X$ be a unit random eigenvector of $W-\E\, W$
measurable with respect to $W$, corresponding to the largest (by absolute value) eigenvalue of $W-\E\, W$.
Then
\begin{align*}
\beta^{2}(1-p)^2\E\|W' X\|_2^2-\E\|W-\E\, W\|^2\geq
-Ch p \max(\log n,n p)^2-C h p^2 n \max(\log n,n p),
\end{align*}
where $h$ denotes the uniform upper bound for $\xi^2$ and $C>0$ is a universal constant.
\end{lemma}
\begin{proof}
Let variables $a,b',b,\xi_\varepsilon$ be as in the last lemma. Conditioned on $\xi$ and on $b=0$, we can calculate 
$$
\Prob\{b'=1\,\vert\,\xi; b=0\} =  \frac{\Prob\{b'=1 \text{ and } b=0\vert\,\xi\}}{\Prob\{b=0\vert\,\xi\}}= \frac{\Prob\{b'=1 \text{ and } a=0\vert\,\xi\}}{\Prob\{b=0\}}=
\frac{p/\varepsilon-p }{1-p},
$$
where we used that $\xi$ is independent from $a,b,b'$. We also deduce that 
$$
\Prob\{b'=0\,\vert\,\xi;b=0\}=\frac{1-p/ \varepsilon}{1-p}.
$$
Observe that on the event $\{b=0\}$ we have
$$
b'\xi_\varepsilon'=\beta^{-1}\, a b'\,\xi-\frac{\beta^{-1} \varepsilon(b'-ab')\E\xi}{1- \varepsilon}
=-\frac{\beta^{-1} \varepsilon\,b'\,\E\xi}{1- \varepsilon}.
$$
Hence, the conditional expectation of $b'\xi_\varepsilon'$ given $\xi$ and $\{b=0\}$, is
\begin{equation}\label{eq: aux 24928357023985}
\E\big(b'\xi_\varepsilon'\,|\,\xi,\,b=0\big)=\frac{1-p/ \varepsilon}{1-p}\cdot 0-\frac{p/\varepsilon-p}{1-p}
\frac{\beta^{-1}\varepsilon\,\E\xi}{1- \varepsilon}
=-\frac{\beta^{-1} p\,\E\xi}{1-p},
\end{equation}
and the conditional second moment of $b'\xi_\varepsilon'$ is
\begin{equation}\label{eq: aux 2096305250958}
\E\big((b'\xi_\varepsilon')^2\,|\,\xi,\,b=0\big)=
\frac{p/\varepsilon-p}{1-p}\frac{\beta^{-2} \varepsilon^2\,(\E\xi)^2}{(1- \varepsilon)^2}
=\frac{p}{1-p}\frac{\beta^{-2} \varepsilon\,(\E\xi)^2}{1- \varepsilon}.
\end{equation}
Let the entries of $W$ be represented in the form $(b_{ij}\xi_{ij})_{ij}$, where each pair $(b_{ij},\xi_{ij})$
is equidistributed with $(b,\xi)$.
Condition on any realization of $(b_{ij},\xi_{ij})_{1\leq i,j\leq n}$,
fix an index $k\leq n$, and let $J_k\subset[n]$ be the collection of all indices $j$
such that $b_{kj}=1$. Then, denoting the entries of $W'$ by $w_{ij}'$,
we get $w'_{kj}=\beta^{-1} \xi_{kj}$ for all $j\in J_k$.
The square of the scalar product of $X=(x_1,x_2,\dots,x_n)\in S^{n-1}$ and the $k$-th row of
$A:=W-\E\, W$
can be written as
$$\langle \Row_k(A),X\rangle^2=\Big(\sum\limits_{j\in J_k}(\xi_{kj}-p\E\xi)x_j+\sum\limits_{j\notin J_k}(-p\E\xi)x_j\Big)^2.$$
On the other hand, in view of \eqref{eq: aux 24928357023985}--\eqref{eq: aux 2096305250958} the conditional second moment of
the scalar product of $X$ with the $k$-th row of $W'$ can be computed as
\begin{align*}
\E_C\langle \Row_k(W'),X\rangle^2&=
\E_C\Big(\sum\limits_{j\in J_k}w'_{kj} x_j+\sum\limits_{j\notin J_k}w'_{kj}x_j\Big)^2\\
&\hspace{-2.5cm}=\E_C\Big(\sum\limits_{j\in J_k}\beta^{-1} \xi_{kj} x_j+\sum\limits_{j\notin J_k}w'_{kj}x_j\Big)^2\\
&\hspace{-2.5cm}=\Big(\sum\limits_{j\in J_k}\beta^{-1} \xi_{kj} x_j+\sum\limits_{j\notin J_k}\E_C w'_{kj}x_j\Big)^2
+\sum\limits_{j\notin J_k}\E_C(w'_{kj})^2\,x_j^2-\sum\limits_{j\notin J_k}(\E_C w'_{kj})^2\,x_j^2\\
&\hspace{-2.5cm}=\Big(\sum\limits_{j\in J_k}\beta^{-1} \xi_{kj} x_j-\sum\limits_{j\notin J_k}\frac{\beta^{-1} p\,\E\xi}{1-p}\,x_j\Big)^2
+\sum\limits_{j\notin J_k}\frac{p}{1-p}\frac{\beta^{-2} \varepsilon\,(\E\xi)^2}{1-\varepsilon}\,x_j^2-
\sum\limits_{j\notin J_k}\frac{\beta^{-2} p^2\,(\E\xi)^2}{(1-p)^2}\,x_j^2,
\end{align*}
where, for brevity, we write $\E_C$ for conditional expectation given a realization of $(b_{ij},\xi_{ij})_{1\leq i,j\leq n}$.
Thus,
\begin{align*}
\beta^{2}(1-p)^2&\E_C\langle \Row_k(W'),X\rangle^2-\langle \Row_k(A),X\rangle^2\\
&=\Big(\sum\limits_{j\in J_k}(1-p) \xi_{kj} x_j-\sum\limits_{j\notin J_k}p\,\E\xi\,x_j\Big)^2
-\Big(\sum\limits_{j\in J_k}(\xi_{kj}-p\E\xi)x_j-\sum\limits_{j\notin J_k}p\, \E\xi\, x_j\Big)^2\\
&\hspace{1cm}+\sum\limits_{j\notin J_k}p(1-p)\frac{ \varepsilon\,(\E\xi)^2}{1-\varepsilon}\,x_j^2-
\sum\limits_{j\notin J_k}p^2\,(\E\xi)^2\,x_j^2.
\end{align*}
Factorizing the first difference and using that $p\leq \varepsilon$ for the second one, we get  
\begin{align*}
\beta^{2}(1-p)^2&\E_C\langle \Row_k(W'),X\rangle^2-\langle \Row_k(A),X\rangle^2\\
&\geq\Big(\sum\limits_{j\in J_k}p(\E\xi- \xi_{kj}) x_j\Big)
\Big(\sum\limits_{j\in J_k}(2\xi_{kj}-p\xi_{kj}-p\E\xi)x_j-2\sum\limits_{j\notin J_k}p\,\E\xi\, x_j\Big)\\
&\geq -2\sqrt{h}\, p\sum\limits_{j\in J_k}|x_j|\,\Big(2\sqrt{h}\sum\limits_{j\in J_k}|x_j|+4\sqrt{h}\,  p\sum\limits_{j=1}^n |x_j|\Big)\\
&\geq -4h p|J_k|\sum\limits_{j\in J_k}x_j^2-8h p^2\sqrt{n}\sum\limits_{j\in J_k}|x_j|,
\end{align*}
where $h$ denotes the uniform upper bound for $\xi^2$.

Set $K:=\max\limits_{k\leq n}|J_k|$, and note that by the symmetry of $W$, every index $j\in [n]$ can belong to at most $K$ of the sets $J_k$'s.  
Using this and taking the sum over all $k\leq n$ in the above relation, we get 
\begin{align*}
\beta^{2}(1-p)^2\E_C\|W' X\|_2^2-\|A X\|_2^2
&\geq \sum\limits_{k=1}^n\Big(-4h p|J_k|\sum\limits_{j\in J_k}x_j^2-8h p^2\sqrt{n}\sum\limits_{j\in J_k}|x_j|\Big)\\
&\geq -4h p K^2-8h p^2\sqrt{n} K\sum\limits_{k=1}^n |x_k|\\
&\geq -4h p K^2-8h p^2 n K.
\end{align*}
Removing the conditioning on a realization of $W$, we obtain from the last relation
\begin{align*}
\beta^{2}(1-p)^2\E\|W' X\|_2^2-\E\|A\|^2&\geq
-4h p\, \E K^2-8h p^2 n\,  \E K\\
&\geq
-Ch p \max(\log n,n p)^2-C h p^2 n \max(\log n,n p)
\end{align*}
for a universal constant $C>0$, where we have applied Bernstein's inequality to estimate the
moments of $K$.
\end{proof}

\begin{prop}\label{prop: coupling}
Let $\xi$ be a uniformly bounded real random variable with $\E\xi\neq 0$ and  unit second moment,
and let $(p_n)_{n\geq 1}$
be a sequence of positive real numbers in $(0,1]$ such that $\lim\limits_{n\to\infty} np_n=\infty$
and
$$\lim\limits_{n\to\infty}\big(p_n \max(\log n,n p_n)^2\big)=0.$$
Further, for each $n\geq 1$, let $W_n$ be an $n\times n$ symmetric random matrix with i.i.d.\ entries
(up to the symmetry constraint) equidistributed with
$b_n\xi$, where $b_n$ is a Bernoulli ($0/1$) random variable with probability of success $p_n$, independent from $\xi$.
Then for any $\varepsilon>0$, there is $n_\varepsilon\geq 1$,
a uniformly bounded {\bf centered} random variable $\xi_\varepsilon'$ of unit variance and a
coupling $(W_n,W_n')_{n=n_\varepsilon}^\infty$ of sequences of random matrices with the following properties:
\begin{itemize}

\item For every $n\geq n_\varepsilon$, $W_n'$ is an $n\times n$ symmetric random matrix with independent entries;

\item All entries of $W_n'$ are equidistributed with $b_n'\xi_\varepsilon'$,
where $b_n'$ is a Bernoulli random variable independent from $\xi_\varepsilon'$, with $\lim\limits_{n\to\infty}(n\,\Prob\{b_n'=1\})=\infty$;

\item $
{\E\max\limits_{i\leq n}\|\Row_i(W_n')\|_2^2}\big/{\E\max\limits_{i\leq n}\|\Row_i(W_n)\|_2^2}\leq
\frac{(1+\varepsilon)\Prob\{b_n'=1\}}{p_n}$ for all $n\geq n_\varepsilon$;

\item We have
$$\lim\limits_{n\to\infty}\Prob\big\{\Vert W_n-\E W_n\Vert \leq (1+\varepsilon)\,\sqrt{p_n/\Prob\{b_n'=1\}}\,\|W_n'\|\big\}=1.$$

\end{itemize} 
\end{prop}
\begin{proof}
Fix any $\varepsilon\in(0,1/2]$, and set $\varepsilon'>0$ to be the largest number in $(0,1]$
satisfying the conditions
$$\frac{ \varepsilon'}{1- \varepsilon'}\leq \frac{\varepsilon}{(\E\xi)^2}\quad\mbox{ and }
\quad \bigg(1+\frac{\varepsilon'\,\vert\E\, \xi\vert}{1-\varepsilon'}\bigg)^2
\leq (1+\varepsilon)\bigg(1+\frac{\varepsilon'\,(\E\,\xi)^2}{1-\varepsilon'}\bigg).$$
We can choose $n_\varepsilon$ large enough so that $p_n\leq\varepsilon'$ for all $n\geq n_\varepsilon$.
Fix for a moment any $n\geq n_\varepsilon$.
As the first step of the proof, we define the random variables $\xi_{\varepsilon'}$, $a$, $b_n$, $b_n'$:
we assume that $\xi,a,b_n'$ are jointly independent, where $a$ and $b_n'$
are Bernoulli with probabilities of success $\varepsilon'$ and $p_n/\varepsilon'$,
respectively, and set
$$\xi_{\varepsilon'}':= \beta^{-1}\,a\,\xi-\frac{\beta^{-1}\, \varepsilon'}{1- \varepsilon'}(1-a)\E\xi;$$
where $\beta:=\sqrt{\Var\big(a\,\xi-\frac{ \varepsilon'}{1- \varepsilon'}(1-a)\E\xi\big)}
=\sqrt{\varepsilon'+\frac{ {\varepsilon'}^{2}(\E\xi)^2}{1-\varepsilon'}}$.

Now we can define the required coupling: for the given $n$, let $W_n'$ be an $n\times n$ symmetric random matrix
such that pairs of respective entries of $W_n$ and $W_n'$,
$\big\{\big((W_n)_{ij},(W_n')_{ij}\big),\;j\geq i\big\}$, are i.i.d. and equidistributed
with the pair $(b_n\xi,b_n'\xi_{\varepsilon'}')$.
Then Lemma~\ref{l: entrywise coupling}, the above definitions, and the choice of $\varepsilon'$, imply that matrices
$W_n'$, $n\geq n_\varepsilon$, satisfy the first three assertions of the proposition.
Thus, it remains to verify the fourth assertion.

For each $n$, let $X_n$ be a unit random eigenvector of $A_n:=W_n-\E\, W_n$
measurable with respect to $W_n$, corresponding to the largest (by absolute value) eigenvalue of $A_n$.
Applying Lemma~\ref{l: coupling 2}, we get
$$
\beta^{2}(1-p_n)^2\E\|W' X\|_2^2-\E\|A_n\|^2\geq
-Ch p \max(\log n,n p)^2-C h p^2 n \max(\log n,n p),
$$
for a universal constant $C>0$, where $h$ is the uniform upper bound for $\xi^2$.
Hence, in view of the condition $\lim\limits_{n\to\infty}\big(p_n \max(\log n,n p_n)^2\big)=0$,
and since $\lim\limits_{n\to\infty}\E\|A_n\|^2=\infty$, we get that for every $\delta>0$,
$$\lim\limits_{n\to\infty}\Big(\beta^2(1-p_n)^2(1+\delta)\E\|W_n'\|^2-\E\|A_n\|^2\Big)=\infty.$$

It remains to apply Talagrand's concentration inequality (Theorem~\ref{th: tal})
to $\|W_n'\|,\|A_n\|$:
\begin{align*}
&\lim\limits_{n\to\infty}\Prob\big\{\|W_n'\|^2\leq (1-\delta)\E\|W_n'\|^2\big\}
=0;\\
&\lim\limits_{n\to\infty}\Prob\big\{\|A_n\|^2\geq (1+\delta)\E\|A_n\|^2\big\}=0.
\end{align*}
This, together with the previous assertion, gives
$$\lim\limits_{n\to\infty}\Prob\bigg\{\frac{\beta^2(1-p_n)^2(1+\delta)}{(1-\delta)}\|W_n'\|^2\geq
(1+\delta)^{-1}\|A_n\|^2\bigg\}=1\quad \mbox{ for every $\delta>0$}.$$
Finally, choose $\delta>0$ so that $\frac{\beta^2(1-p_n)^2(1+\delta)^2}{(1-\delta)}\leq\varepsilon' (1+\varepsilon)^2$
for all large $n$
(such $\delta$ exists in view of the choice of $\varepsilon'$). The result follows.
\end{proof}

Using the coupling provided by the previous proposition, we can now prove the main statement of this section. 
\begin{theorem}\label{th: non-centered upper}
Let $\xi$ be a uniformly bounded real random variable with unit second moment.
For each $n$, let $W_n$ be an $n\times n$ symmetric random matrix with i.i.d.\ entries (up to the symmetry constraint),
with each entry equidistributed with the product $b_n\xi$, where $b_n$ is $0/1$ (Bernoulli) random variable independent of $\xi$,
with a probability of success equal to $p_n$.
Assume further that $n p_n\to\infty$ with $n$.
For each $n$, define the random quantities
$$\rho_n:=\theta_n+\frac{n p_n}{\theta_n},\quad \theta_n:=\sqrt{\max\big(\max\limits_{i\leq n}\|\Row_i(W_n)\|_2^2-np_n,n p_n\big)}.$$
Then for any $\varepsilon>0$ we have
$$\lim\limits_{n\to\infty}\Prob\big\{\Vert W_n-\E\, W_n\Vert \leq (1+\varepsilon)\rho_n\big\}=1.$$
\end{theorem}
\begin{proof}
Without loss of generality, we can assume that $\lim\limits_{n\to\infty}\big(p_n \max(\log n,n p_n)^2\big)=0$:
indeed, when $\lim\limits_{n\to\infty}\frac{n p_n}{\log n}=\infty$, standard concentration inequalities imply that the sequence
$$((np_n)^{-1}\max\limits_{i\leq n}\|\Row_i(W_n)\|_2^2)_{n=1}^\infty$$
converges to one
in probability, so the assertion of the theorem is equivalent to 
$$\lim\limits_{n\to\infty}\Prob\big\{\Vert W_n-\E\, W_n\Vert \leq (1+\varepsilon)2\sqrt{n p_n}\big\}=1\quad\mbox{ for any $\varepsilon>0$}.$$
This, in turn, is a known result; see \cite{Khorunzhy,BGBK,LVY}.

Below, we work under the assumption $\lim\limits_{n\to\infty}\big(p_n \max(\log n,n p_n)^2\big)=0$.
Take any $\varepsilon>0$, and let $b_n'$, $\xi_\varepsilon'$ be the random variables
and $(W_n,W_n')_n$ be the coupling of sequences of random matrices
from Proposition~\ref{prop: coupling}. For each $n$, denote $p_n':=\Prob\{b_n'=1\}$, so that $\lim\limits_{n\to\infty} p_n' n=\infty$.
Hence, applying Theorem~\ref{th: centered upper} to $W_n'$, we get
$$
\lim\limits_{n\to\infty}\Prob\big\{\|W_n'\|\leq (1+\varepsilon)\rho_n'\big\}=1,
$$
where 
$\rho_n':=\theta_n'+\frac{n p_n'}{\theta_n'}$
and
$$\theta_n':=\sqrt{\max\big(\max\limits_{i\leq n}\|\Row_i(W_n')\|_2^2-np_n',n p_n'\big)}.$$
Combined with Proposition~\ref{prop: coupling}, this gives
$$
\lim\limits_{n\to\infty}\Prob\big\{\Vert W_n-\E\, W_n\Vert\leq (1+\varepsilon)^2\,\sqrt{p_n/p_n'}\,\rho_n'\big\}=1.
$$
To prove the statement, it remains to compare the quantities $\rho_n$ and $\rho_n'$
for $n$ tending to infinity.
According to Proposition~\ref{prop: coupling},
$${\E\max\limits_{i\leq n}\|\Row_i(W_n')\|_2^2}\big/{\E\max\limits_{i\leq n}\|\Row_i(W_n)\|_2^2}\leq
\frac{(1+\varepsilon)p_n'}{p_n}$$
for all $n\geq n_\varepsilon$.
Together with Talagrand's concentration inequality (Theorem~\ref{th: tal}) and the assumption $\lim\limits_{n\to\infty} p_n n=\infty$,
this implies
$$\lim\limits_{n\to\infty}\Prob\bigg\{\max\limits_{i\leq n}\|\Row_i(W_n')\|_2^2
\leq \frac{(1+\varepsilon)^2 p_n'}{p_n}\,\max\limits_{i\leq n}\|\Row_i(W_n)\|_2^2\bigg\}=1.$$
Noting that the function $f(x,y)= \frac{\max(x,2y)}{\sqrt{\max (x-y,y)}}$ is coordinate-wise increasing, we get with probability tending to one with $n$ that 
\begin{align*}
\rho_n'= f\big( \max\limits_{i\leq n}\|\Row_i(W_n')\|_2^2, np_n'\big)&\leq   f\big( (1+\varepsilon)^2p_n' \max\limits_{i\leq n}\|\Row_i(W_n)\|_2^2/p_n, (1+\varepsilon)^2np_n'\big)\\
&=(1+\varepsilon) \sqrt{\frac{p_n'}{p_n}}\, f\big( \max\limits_{i\leq n}\|\Row_i(W_n)\|_2^2, np_n\big).
\end{align*}
This implies that 
$$
\lim\limits_{n\to\infty}\Prob\big\{\rho_n'\leq (1+\varepsilon)\sqrt{p_n'/p_n}\,\rho_n\big\}=1,
$$
and the result follows.
\end{proof}

\section{Lower bound for largest eigenvalues}

The main result of this section provides a lower bound on the largest eigenvalues
of a sparse random symmetric matrix, which asymptotically matches the upper bounds given in Theorems~\ref{th: centered upper}
and~\ref{th: non-centered upper}, thus completing the proofs of Theorems~A and~B
from the introduction. Our approach to finding lower bounds on the $k$--largest eigenvalue
is completely different from the combinatorial methods from the first part of the paper,
and is based on explicitly constructing an ``approximate eigenspace'' corresponding to the first few largest eigenvalues
(or singular values) of the matrix.

\begin{theorem}\label{th: lower top}
Let $\xi$ be a uniformly bounded random variable with unit second moment.
Further, let $(W_n)$ be a sequence of random matrices, where for each $n\geq 1$,
$W_n$ is $n\times n$ symmetric, with i.i.d.\ entries above the main diagonal (and zeros on the diagonal),
each entry equidistributed with $b_n\xi$, where $b_n$ is a Bernoulli ($0/1$) random variable
with probability of success $p_n$. Assume further that $\lim\limits_{n\to\infty}n p_n=\infty$,
and $n p_n\leq \log^{2}n$ for all large enough $n$. Set
$$\rho_n:=\theta_n+\frac{n p_n}{\theta_n},\quad \theta_n:=\sqrt{\max\big(\max\limits_{i\leq n}\|\Row_i(W_n)\|_2^2-np_n,n p_n\big)}.$$
Then for any $k\in\N$, denoting by $\lambda_{|k|}(W_n)$ the $k$-th largest (by absolute value)
eigenvalue of $W_n$, we have
$$\lim\limits_{n\to\infty}\Prob\Big\{|\lambda_{|k|}(W_n)|\geq (1-\varepsilon)\rho_n\Big\}=1\quad \mbox{for all }\varepsilon>0.$$  
\end{theorem}
\begin{Remark}
In view of Lemma~\ref{lem: upper-zero-diag}, the theorem remains true if the condition that matrices $W_n$ have zero diagonal
is replaced with the condition that the diagonal entries are independent and equidistributed with the off-diagonal.
\end{Remark}
\begin{Remark}
The case of relatively denser matrices, with $\lim\limits_{n\to\infty}\frac{n p_n}{\log n}=\infty$
and $p_n\to 0$, is known
and immediately follows from standard Bernstein--type concentration inequalities
for $\Row_i(W_n)$ and the property that the spectrum of appropriately normalized matrices $W_n$
converges to the semi-circle distribution.
Indeed, it is not difficult to check that in this setting $\rho_n/(2\sqrt{n p_n})$ converges to one in probability,
i.e.\ the theorem amounts to checking that $\frac{1}{\sqrt{n p_n}}\lambda_{|k|}(W_n)$ is asymptotically outside of the support
of the semi-circle distribution.
Although only the ``global'' limiting law is required here,
let us remark that a local semi-circle law for sparse symmetric matrices was recently established in 
\cite{EKYY,HKM}.
\end{Remark}

To prove the theorem, we consider several preparatory lemmas.
The next lemma encapsulates some simple structural properties of sparse Erd\H os--Renyi graphs,
and its proof is given mostly for completeness.
The only element which makes the lemma different from absolutely standard observations
is the assumption that some edges of the graph are ``frozen'' (non-random).
This assumption will be important later in a decoupling trick which we apply in the proof of the theorem.

\begin{lemma}\label{l: graph str}
For any positive integers $k,q$ and a real number $\delta>0$ there is $n_0=n_0(k,q,\delta)\in\N$
with the following property.
Let $n\geq n_0$, $\ell\geq k+1$, and let $B=(b_{ij})$ be an $n\times n$ symmetric $0/1$ random matrix such that
the entries
$$b_{ji}=b_{ij},\quad 1\leq i\leq k,\;\;\ell\leq j\leq n$$
are fixed (non-random), with $1\leq \sum_{j=\ell}^n b_{ij}\leq \log^{4} n$ for every $1\leq i\leq k$,
$\sum_{i=1}^k b_{ij}\leq 1$ for every $\ell\leq j\leq n$;
and the remaining off-diagonal
entries are i.i.d (up to the symmetry constraint) Bernoulli random variables with probability of success $p$ satisfying 
$np\in [(\log\log n)^2,  n^{\frac{1}{4q}}]$. Further, let $G$ be the simple random graph on $[n]$ with adjacency matrix $B$.
Then
$$\Prob\big\{\mbox{$i_1$ connected to $i_2$ by a path of length at most $2q$ in $G$ for some $1\leq i_1\neq i_2\leq k$}\big\}
\leq n^{-1/8},$$
and
\begin{align}
\Prob\big\{&\mbox{$q$--neighborhood of every vertex $i\in[k]$ in $G$ is a tree, where}\nonumber\\
&\mbox{every leaf has depth $q$, and}\label{lower-bound-event}\\
&\mbox{the degrees of all vertices of the trees, except for the roots}\nonumber\\
&\mbox{and leaves, are in the range $[(1-\delta )pn,(1+\delta)pn]$}\big\}\geq 1-\delta.\nonumber
\end{align}

\end{lemma}
\begin{proof}
The only point that needs some attention is that the edges connecting $[k]$ with $[n]\setminus[\ell]$
are ``frozen''. Hence, our argument will involve, as an additional step, separate treatment of those edges
and the rest of the graph.

To prove the first assertion of the lemma, let us fix $1\leq i_1\neq i_2\leq k$.
Clearly,
\begin{align*}
\Prob\big\{\mbox{$i_1, i_2$ connected by a path of length $1$}\big\}=p.
\end{align*}
Further, since for any $\ell\leq j\leq n$, we have $\sum_{i=1}^k b_{ij}=1$, then no two entries $b_{i_1 j}$, $b_{i_2 j}$ are simultaneously equal to one. Thus, we have
\begin{align*}
\Prob\big\{\mbox{$i_1, i_2$ connected by a path of length $2$}\big\}
= \Prob\big\{\mbox{$i_1\edg u\edg i_2$ for some $u< \ell$}\big\}\leq np^2.
\end{align*}
For any $3\leq v \leq 2q$,
\begin{align*}
\Prob\big\{&\mbox{$i_1, i_2$ connected by a path of length $v$ not passing through $[k]\setminus\{i_1,i_2\}$}\big\}\\
&= \Prob\big\{ \exists u_1\neq\ldots\neq u_{v-1}\in [n]\setminus[k]: \, i_1\edg u_1\edg u_2\ldots u_{v-1}\edg i_2\big\}\\ 
&\leq \Prob\big\{ \exists u_1\neq\ldots\neq u_{v-1}\in [n]\setminus[k]: u_1,u_{v-1}<\ell,\; i_1\edg u_1\edg u_2\ldots u_{v-1}\edg i_2\big\}\\
&\hspace{1cm}+\Prob\big\{ \exists u_1\neq\ldots\neq u_{v-1}\in [n]\setminus[k]: u_1<\ell,\,u_{v-1}\geq\ell,\; i_1\edg u_1\edg u_2\ldots u_{v-1}\edg i_2\big\}\\
&\hspace{1cm}+\Prob\big\{ \exists u_1\neq\ldots\neq u_{v-1}\in [n]\setminus[k]: u_1\geq\ell,\,u_{v-1}<\ell,\; i_1\edg u_1\edg u_2\ldots u_{v-1}\edg i_2\big\}\\
&\hspace{1cm}+\Prob\big\{ \exists u_1\neq\ldots\neq u_{v-1}\in [n]\setminus[k]: u_1,u_{v-1}\geq\ell,\; i_1\edg u_1\edg u_2\ldots u_{v-1}\edg i_2\big\}\\
&\leq n^{v-1}p^{v}
+ 2\big(\log^{4} n\big)n^{v-2}p^{v-1}
+\big(\log^{4} n\big)^2 n^{v-3}p^{v-2}.
\end{align*}
Summing these probabilities over all $v \in \{1,\ldots, 2q\}$ and using that $np\leq n^{\frac{1}{4q}}$, we get 
$$
\Prob\big\{\mbox{$i_1, i_2$ connected by a path of length at most $2q$ not passing through $[k]$}\big\}
\leq 4(2q)p \sqrt{n}.
$$
Now a union bound over all choices of $1\leq i_1\neq i_2\leq k$ implies that 
\begin{align*}
\Prob\big\{&\mbox{$i_1$ connected to $i_2$ by a path of length at most $2q$ in $G$ for some $1\leq i_1\neq i_2\leq k$}\big\}\\
&\leq k^2\cdot 4(2q)p \sqrt{n}\leq n^{-1/8},
\end{align*}
where we used that $p\leq n^{-3/4}$
and assumed that $n_0$ is sufficiently large.

\medskip

To prove the second part of the lemma, let us denote by $\Event$ the event appearing in \eqref{lower-bound-event}. 
For every $i\in [k]$, we denote by $\Event_i$ the event that $B(i,q)$ --- the $q$--neighborhood of $i$ in $G$
--- is a tree, where every leaf has depth $q$, and the degrees of all vertices of the tree, except for the root and leaves,
are in the range $[(1-\delta )pn,(1+\delta)pn]$. With these notations, we have $\Event=\bigcap_{i\in [k]} \Event_i$. 

Let us fix $i\in [k]$, and denote by $\Event_{i,1}$ the event that $B(i,q)$ is a tree, by
$\Event_{i,2}$ the event that every degree one vertex in $B(i,q)$, except for the root $i$, is at distance $q$ from $i$, and by
$\Event_{i,3}$ the event that the degrees of all vertices in $B(i,q)$ except for $i$ and those at distance $q$, are in the range
$[(1-\delta )pn,(1+\delta)pn]$. Clearly, we have $\Event_i=\Event_{i,1}\cap \Event_{i,2}\cap \Event_{i,3}$.
We will show that each of these events has probability close to one. 
Observe that 
\begin{align*}
\Prob&\big( \Event_{i,1}^c\big)\\
& \leq \Prob\big\{ \exists\; 3\leq v\leq 2q+1 \mbox{ and a subgraph
of $B(i,q)$ with $v$ edges and $v$ vertices containing $i$}
\big\}\\
&\leq\sum_{v=3}^{2q+1} \Prob\big\{\mbox{There exists a subgraph of $B(i,q)$ with $v$ edges and $v$ vertices containing $i$}\big\}\\
&\leq\sum_{v=3}^{2q+1}\sum\limits_{w=0}^{v-1}
\Prob\big\{\mbox{There exists a subgraph of $G$ with $v$ edges and $v$ vertices,}\\
&\hspace{2.5cm}\mbox{containing $i$ and with $w$ edges connecting $[k]$ to $[n]\setminus[\ell-1]$}\big\}.
\end{align*}
Fix any admissible parameters $v,w$.
Then the total number of choices of vertices of the subgraph can be roughly estimated from above by
$
k^w\big(\log^{4} n\big)^w\,n^{v-w-1},
$
and the probability that there exist $v$ edges among the chosen vertices
--- by $p^{v-w}(2q)^{2(v-w)}$.
Hence,
$$
\Prob\big(\Event_{i,1}^c\big)\leq \sum_{v=3}^{2q+1}\sum\limits_{w=0}^{v-1}
(2q)^{2(v-w)}k^w (\log n)^{4w} \,n^{v-w-1} p^{v-w}\leq
\frac{4q^2}{n} \big(4q^2k\max(np,\log^{4} n)\big)^{2q+1}\leq \frac{\delta}{3k},
$$
where we have used that $np\leq n^{\frac{1}{4q}}$ and $n$ is sufficiently large. 

To estimate $\Event_{i,2}^c,\Event_{i,3}^c$, let us denote by $\widetilde \Event$
the event that any two vertices of $G$ from $[k]$ are at distance at least $2q+1$ from each other.
Observe that, for $q=1$, $\Prob(\Event_{i,2}^c\;|\;\Event_{i,1})=0$ (hence $\Prob\big(\Event_{i,2}^c\big)
\leq \frac{\delta}{3k}$), while for $q>1$ we can write 
\begin{align*}
\Prob&\big(\widetilde \Event\cap\Event_{i,2}^c\big)\\
&= \Prob\big(\widetilde \Event\cap\big\{\mbox{There exists a degree $1$ vertex of $G$ $(\neq i)$ at distance at most $q-1$ from $i$}\big\}\big)\\
&\leq \sum\limits_{v=1}^{q-1}\Prob\big\{\exists\,
u_1\neq\ldots\neq u_{v}\in [n]\setminus[k]: u_1\geq\ell,\; i\edg u_1\edg u_2\ldots \edg u_{v},\;\degree(u_{v})=1\big\}\\
&\hspace{0.5cm}+\sum\limits_{v=1}^{q-1}\Prob\big\{\exists\,
u_1\neq \ldots\neq u_{v}\in [n]\setminus[k]: u_1<\ell,\; i\edg u_1\edg u_2\ldots \edg u_{v},\;\degree(u_{v})=1\big\}\\
&\leq \log^{4} n\,\sum\limits_{v=1}^{q-1}(np)^{v-1}(1-p)^{n-2}
+\sum\limits_{v=1}^{q-1}(np)^{v}(1-p)^{n-2}\leq \frac{\delta}{6k},
\end{align*}
where the last inequality follows assuming $n$ is sufficiently large, since $pn\geq (\log\log n)^2$.

Finally, we have for $q\geq 2$:
\begin{align*}
\Prob&\big(\widetilde \Event\cap\Event_{i,3}^c\big)\\
&\leq \sum_{v=1}^{q-1} \Prob\big(\widetilde \Event\cap\big\{\mbox{$\exists u\in[n]$ at distance $v$ from $i$
with $\degree(u)\not\in [(1-\delta )pn,(1+\delta)pn]$}\big\}\big)\\
&\leq \sum_{v=1}^{q-1} \Prob\big(\widetilde \Event\cap\big\{ \exists u_1\neq u_2\neq\ldots  \neq u_v\in[n]\setminus[k] \mbox{ such that }
 i\edg u_1\edg u_2\ldots u_{v-1}\edg u_v, \\
&\qquad\qquad\qquad  \mbox{and $\degree(u_v)\not\in [(1-\delta )pn,(1+\delta)pn]$}\big\}\big)\\
&\leq \sum_{v=1}^{q-1}\sum_{\substack{u_1\neq \ldots \neq u_v\in[n]\setminus[k]}}
\Prob\big\{i\edg u_1\edg u_2\ldots u_{v-1}\edg u_v\big\} \\
&\qquad\qquad\qquad\qquad \cdot   \Prob\big\{\mbox{ $\big\vert \nbr(u_v)
\setminus \{u_{v-1}\}\big\vert\not\in [(1-\delta )pn-1,(1+\delta)pn-1]$}\big\}\\
&\leq q\big(\log^{4} n\big)(np)^{q-2}e^{-c\delta^2 np}
+
q(np)^{q-1} e^{-c\delta^2 np},
\end{align*}
for some constant $c$, where we have used Lemma~\ref{l:bernstein} to get the last inequality,
and the two terms correspond to the cases $u_1\geq \ell$ and $u_1<\ell$ in the path representation.
It remains to note that this quantity can be bounded by $\frac{\delta}{6k}$ using that $np\geq (\log\log n)^2$
and $n$ is sufficiently large. 

Putting these estimates together, we deduce that for any $i\in [k]$
$$
\Prob\big(\Event_i^c\big)\leq \Prob\big( \Event_{i,1}^c\big)+\Prob\big(\widetilde \Event\cap \Event_{i,2}^c\big)+
\Prob\big(\widetilde \Event\cap\Event_{i,3}^c\big)+\Prob\big(\widetilde \Event^c\big)\leq \frac{\delta}{k},
$$
as long as $n$ is large enough so that $n^{-1/8}\leq \frac{\delta}{3k}$.
It remains to apply a union bound over all $i\in [k]$ to finish the proof. 
\end{proof}

The next lemma provides a simple decoupling argument for
the quantity $\max\limits_{1\leq i\leq n}\|\Row_i(W)\|_2$,
where $W$ is a symmetric random matrix with independent entries.
While the symmetry constraint induces dependencies between the matrix rows,
it can be shown that, under some assumptions on the distribution of the entries,
the maximum of the row-norms is close to the maximum taken in a rectangular submatrix of $W$
entirely contained above the main diagonal, so that its entries are jointly independent.

\begin{lemma}\label{l: rows restricted}
Let $\xi$ be a uniformly bounded random variable. 
Then for any $\delta, \alpha>0$ there are
$n_{\ref{l: rows restricted}}\in \N$ and
$r_{\ref{l: rows restricted}}\in(0,1)$ depending only on the distribution of $\xi$, on $\delta$ and $\alpha$, with the following property.
Let $n\geq n_{\ref{l: rows restricted}}$, and let $W=(w_{ij})$ be an $n\times n$ random symmetric matrix with independent
(up to the symmetry constraint) entries equidistributed with $b\xi$, where $b$ is Bernoulli ($0/1$)
random variable independent from $\xi$, with probability of success $p\geq \alpha \log n/n$
(the diagonal entries of $W$ may be either all zeros or be random variables jointly independent and equidistributed
with the off-diagonal entries of $W$).
Denote
$$\eta:=\max\limits_{1\leq i\leq r_{\ref{l: rows restricted}}n}\Big\|\sum\limits_{j=\lfloor r_{\ref{l: rows restricted}}  n\rfloor+1}^n w_{ij}e_j\Big\|_2.$$
Then
$$\Prob\big\{\max\limits_{1\leq i\leq n}\|\Row_i(W)\|_2\geq (1+\delta)\eta\big\}\leq \delta.$$
\end{lemma}
\begin{proof}
Without loss of generality, we suppose that $\E\xi^2=1$. We will also assume that the diagonal entries of the matrix
are equidistributed with the off-diagonal (the case of zero diagonal can be treated with the same method as below,
up to minor adjustments).
Let $h$ be a uniform upper bound for $\xi^2$.
Let $\delta\in (0,1)$ and 
$$
r_{\ref{l: rows restricted}}= \frac{\delta^2}{32} \exp\Big(-\frac{128}{\delta^2\alpha}\Big). 
$$
We will first show that $\eta$ and $\max\limits_{1\leq i\leq r_{\ref{l: rows restricted}}n}\Big\|\sum\limits_{j=1}^n w_{ij}e_j\Big\|_2$ are of the same order with a large probability.
To this aim, we will estimate
the contribution of $\widetilde \eta:= \max\limits_{1\leq i\leq r_{\ref{l: rows restricted}}n}\Big\|\sum\limits_{j\leq \lfloor r_{\ref{l: rows restricted}}  n\rfloor} w_{ij}e_j\Big\|_2$.
We can write
\begin{align}
\Prob\Big\{ \widetilde \eta\geq \frac{\delta}{2} \eta\Big\}
\leq \Prob\Big\{ \widetilde \eta\geq \frac{\delta}{4} \sqrt{np}\Big\}
+\Prob\Big\{ \eta\leq \sqrt{np}/2\Big\}
\leq \Prob\Big\{ \widetilde \eta\geq \frac{\delta}{4} \sqrt{np}\Big\}
+e^{-cnp},
\end{align}
for some constant $c>0$ depending only on the distribution of $\xi$,
where at the last step we applied Bernstein's inequality (Lemma~\ref{l:bernstein}) to
$\sum\limits_{j=\lfloor r_{\ref{l: rows restricted}}  n\rfloor+1}^n w_{1j}e_j$.
Now, by the union bound, 
\begin{align*}
\Prob\Big\{ \widetilde \eta \geq \frac{\delta}{4} \, \sqrt{np}\Big\} &\leq 
 \sum_{i=1}^{\lfloor r_{\ref{l: rows restricted}}  n\rfloor}\, \Prob\Big\{ \Big\|\sum\limits_{j\leq \lfloor r_{\ref{l: rows restricted}}  n\rfloor} w_{ij}e_j\Big\|_2\geq \frac{\delta}{4} \, \sqrt{np}\Big\}\\
  &\leq  \sum_{i=1}^{\lfloor r_{\ref{l: rows restricted}}  n\rfloor} \, \Prob\Big\{ \Big\|\sum\limits_{j\leq \lfloor r_{\ref{l: rows restricted}}  n\rfloor} w_{ij}e_j\Big\|_2^2- \E\,  \Big\|\sum\limits_{j\leq \lfloor r_{\ref{l: rows restricted}}  n\rfloor} w_{ij}e_j\Big\|_2^2 \geq \frac{\delta^2}{32} \, np\Big\}\\
&  =  \sum_{i=1}^{\lfloor r_{\ref{l: rows restricted}}  n\rfloor} \, \Prob\Big\{ \sum\limits_{j\leq \lfloor r_{\ref{l: rows restricted}}  n\rfloor} (w_{ij}^2-\E\, w_{ij}^2)  \geq \frac{\delta^2}{32} \, np\Big\},\\
\end{align*}
where in the second inequality we have used that $\E\,  \Big\|\sum\limits_{j\leq \lfloor r_{\ref{l: rows restricted}}  n\rfloor} w_{ij}e_j\Big\|_2^2= \lfloor r_{\ref{l: rows restricted}}  n\rfloor p \leq \frac{\delta^2}{32} \, np$ by the choice of $r_{\ref{l: rows restricted}}$. Note that the random variables $(w_{ij}^2-\E\, w_{ij}^2)_{j\leq \lfloor r_{\ref{l: rows restricted}}  n\rfloor}$ are independent, centered, of variance at most $ph$ and bounded by $h$. Then an application of Bennett's inequality (see, for example, \cite[Theorem~2.9]{BLM})
implies that for any $i\leq \lfloor r_{\ref{l: rows restricted}}  n\rfloor$ 
$$
\Prob\Big\{ \sum\limits_{j\leq \lfloor r_{\ref{l: rows restricted}}  n\rfloor} (w_{ij}^2-\E\, w_{ij}^2)  \geq \frac{\delta^2}{32} \, np\Big\}
\leq \exp\Big(-\frac{\delta^2 np}{64 } \log \Big(\frac{\delta^2}{32 r_{\ref{l: rows restricted}}}\Big)\Big)\leq \frac{1}{n^2},
$$
by the choice of $r_{\ref{l: rows restricted}}$. We deduce that
$$
\Prob\Big\{ \widetilde \eta\geq \frac{\delta}{2} \eta\Big\}\leq e^{-cnp}+\frac{1}{n}. 
$$
Now note that if $\w \eta\leq \frac{\delta}{2}\eta$ then 
$$
\max\limits_{1\leq i\leq \lfloor r_{\ref{l: rows restricted}}  n\rfloor}\|\Row_i(W)\|_2\leq \eta+\w\eta\leq \big( 1+\frac{\delta}{2}\big) \eta\leq \frac{\big( 1+\delta\big)}{\big( 1+\frac{\delta}{3}\big)} \eta. 
$$
Therefore, it follows that 
\begin{align*}
\Prob&\big\{\max\limits_{1\leq i\leq n}\|\Row_i(W)\|_2\geq (1+\delta)\eta\big\}\\
&\leq \Prob\Big\{\max\limits_{1\leq i\leq n}\|\Row_i(W)\|_2\geq \Big(1+\frac{\delta}{3}\Big)\max\limits_{1\leq i\leq \lfloor r_{\ref{l: rows restricted}}  n\rfloor}\|\Row_i(W)\|_2\Big\}+ e^{-cnp}+\frac{1}{n}.
\end{align*}
It remains to apply Lemma~\ref{l: maxnorm-restrict} and use that $n$ is large enough to finish the proof. 
\end{proof}

The next lemma encapsulates the main construction step of the proof.
In its essense, it gives a procedure for finding a random vector $Y$
which, for a special random matrix $M$ with a ``tree--structure''
of non-zero elements, provides a good approximation of the matrix norm: $\|M\|\approx \|MY\|_2/\|Y\|_2$
with a large probability. In turn, the ratio $\|MY\|_2/\|Y\|_2$ is then related to a quantity
depending on the Euclidean norms of rows of $M$, which will be ultimately connected with the value of $\rho$
in the main theorem.

\begin{lemma}\label{l: main lower}
For any $h>0$ and $\varepsilon>0$
there is odd integer $q_{\ref{l: main lower}}=q_{\ref{l: main lower}}(h,\varepsilon)$,
numbers $d_{\ref{l: main lower}}=d_{\ref{l: main lower}}(h,\varepsilon)>0$ and $\delta_{\ref{l: main lower}}=\delta_{\ref{l: main lower}}(h,\varepsilon)>0$
depending only on $h$ and $\varepsilon$
with the following property.
Let $\xi$ be a random variable such that $\E\xi^2=1$ and $\xi^2\leq h$ a.e.
Let $\tree=(V,E)$ be a rooted tree with a vertex set $V\subset[n]$ and a root $v$,
and assume that every leaf of $\tree$ has depth $q_{\ref{l: main lower}}$. Further, assume that
the degree of every vertex except for the root
and leaves,
is at least $d'$ and at most $\w d$,
where $d',\w d$ satisfy the relations $(1+\delta_{\ref{l: main lower}})d'\geq\w d\geq d_{\ref{l: main lower}}$. 
Let $M=(\mu_{ij})$ be an $n\times n$ random symmetric matrix where
$\mu_{ij}$, $i\edg j\in E$, $i,j\neq v$, are independent (up to the symmetry constraint)
copies of $\xi$, $\mu_{ij}=0$ whenever $i\edg j\notin E$, and $\mu_{vi}=\mu_{iv}$ are fixed
numbers in $\big[-\sqrt{h},\sqrt{h}\big]$ for $v\edg i\in E$
(so that the $v$-th row and column of $M$ are non-random).
Assume further that $\|\Row_v(M)\|_2^2\geq 2(1+\varepsilon)\w d$. 
Then
$$\Prob\Bigg\{\;\|M\|\geq \frac{(1-\varepsilon)\|\Row_v(M)\|_2^2}{\sqrt{\|\Row_v(M)\|_2^2-\widetilde d}}\;\Bigg\}
\geq 1-\varepsilon.$$
\end{lemma}
\begin{proof}
Without loss of generality, we can assume that $\varepsilon$ is bounded above by a small universal constant.
Fix the quantities $q,d$ and $\delta$
(we will discuss below how they should be chosen).
Let the tree $\tree$ and the matrix $M$, as well as the numbers $d',\w d$,
be as in the above statement.

Given any vertex $u\in V$ of depth $r\geq 1$, let $P_u:[0,r]\to V$ be the (unique) path on $\tree$
starting at the root $v$ of the tree and ending at $u$.
Further, for any integer $r\in[0, q]$, denote by $V_r$ the set of all vertices of $\tree$
having depth $r$ (so that, in particular, $V_0=\{v\}$). 
Since every leaf of $\tree$ has depth $q$, then we necessarily have
$$
\nbr(v)=\{P_u(1):\, u\in V_r\},
$$
for any $r\geq 1$. Therefore, we have 
$$
\sum\limits_{u\in V_r}\mu_{P_u(0),P_u(1)}^2 = 
\sum_{P: [0,r]\setminus\{1\}\to V} \sum_{w\in \nbr(v)} \mu_{v,w}^2
= \|\Row_v(M)\|_2^2\cdot \vert\{P: [0,r]\setminus\{1\}\to V\}\vert.
$$
The conditions on the degrees of the
vertices of $\tree$  imply
\begin{equation}\label{eq: aux 087205295}
\|\Row_v(M)\|_2^2\, (d'-1)^{r-1}\leq \sum\limits_{u\in V_r}\mu_{P_u(0),P_u(1)}^2
\leq \|\Row_v(M)\|_2^2\, (\widetilde d-1)^{r-1},\quad r\geq 1.
\end{equation}
Let $\delta_r$, $r\in \Z_2\cap[0,q]$ (where $\Z_2$ are all even integers), be non-negative parameters whose values will be chosen later.
We define a random vector $Y$ in $\R^n$ as
$$Y:=\sum\limits_{r\in \Z_2\cap[0,q]}\sum\limits_{u\in V_r}Y_u,$$
where for each $u\in V_r$ we set
$$Y_u:=\sum\limits_{\tiny\mbox{$z\in V$: $z$ is a child of $u$}}\delta_r \mu_{P_z(0),P_z(1)}
\bigg(\prod\limits_{t=1}^{r}\mu_{P_z(t),P_z(t+1)}\bigg)\; e_z,$$
where $e_z$ is the $z$-th element of the canonical basis. We note that the above product is empty when $r=0$, so that $Y_v=\delta_0\, \Row_v(M)$. In words, $Y_u$ is an $n$-dimensional vector whose coordinates $Y_{u,z}$ are zero if $z$ is not a child of $u$, and equal to 
(up to $\delta_r$) the product of the weights of edges of the unique path leading to $z$ when $z$ is a child of $u$. 
Note that with this definition, we have the vectors $Y_u$ ($u\in V_0\cup V_2\cup\dots\cup V_{q-1}$) have disjoint supports.
Further, by trivially bounding absolute values of all matrix entries by $\sqrt{h}$ and taking into account
the definition of $Y_v$ and \eqref{eq: aux 087205295}, we get a deterministic two-sided inequality
\begin{align*}
\delta_0^2\|\Row_v(M)\|_2^2\leq \|Y\|_2^2&\leq \delta_0^2\|\Row_v(M)\|_2^2+
\sum\limits_{r\in \Z_2\cap[2,q]}\sum\limits_{u\in V_r} \delta_r^2
\mu_{P_u(0),P_u(1)}^2 (\widetilde d-1)h^{r}\\
&\leq
\|\Row_v(M)\|_2^2\sum\limits_{r\in \Z_2\cap[0,q]}\delta_r^2 (\widetilde d-1)^r h^{r}.
\end{align*}
Our proof shall proceed by bounding from below the expectation of $\|MY\|_2^2$,
and bounding from above $\E\|Y\|_2^2$, followed by application of standard concentration inequalities.

For every $r\in \Z_2\cap[2,q]$ and $u\in V_r$, we have
\begin{align*}
&\langle \Row_u(M),Y\rangle^2
=\bigg(\mu_{u,P_u(r-1)}\langle Y,e_{P_u(r-1)}\rangle
+\sum\limits_{\tiny\mbox{$z\in V$: $z$ is a child of $u$}}\mu_{u,z}\langle Y,e_{z}\rangle
\bigg)^2\\
&\hspace{0.5cm}=
\bigg(\delta_{r-2} \mu_{u,P_u(r-1)}\prod\limits_{t=0}^{r-2}\mu_{P_u(t),P_u(t+1)}
+\sum\limits_{\tiny\mbox{$z\in V$: $z$ is a child of $u$}}\delta_r\mu_{u,z} \prod\limits_{t=0}^{r}\mu_{P_z(t),P_z(t+1)}
\bigg)^2\\
&\hspace{0.5cm}=
\bigg(
\Big(\delta_{r-2}+\delta_r\sum_{\tiny\mbox{$z$: $z$ child of $u$}}\mu_{u,z}^2\Big)\prod\limits_{t=0}^{r-1}\mu_{P_u(t),P_u(t+1)}
\bigg)^2,
\end{align*}
where the second line of the formula appears by noting that
$\langle Y,e_{P_u(r-1)}\rangle=\langle Y_{P_{u}(r-2)},e_{P_u(r-1)}\rangle$, and
$\langle Y,e_{z}\rangle=\langle Y_u,e_{z}\rangle$ for any vertex $z$ which is a child of $u$.
Taking the expectation, we get
\begin{align*}
\E \langle \Row_u(M),Y\rangle^2
&\geq \E\bigg(\delta_{r-2}+\delta_r\sum_{\tiny\mbox{$z$: $z$ child of $u$}}\mu_{u,z}^2\bigg)^2
\mu_{P_u(0),P_u(1)}^2\\
&\geq \big(\delta_{r-2}^2+2\delta_{r-2}\delta_r(\deg(u)-1)+\delta_r^2(\deg(u)-1)^2\big)\mu_{P_u(0),P_u(1)}^2\\
&\geq \big(\delta_{r-2}+\delta_r(d'-1)\big)^2\mu_{P_u(0),P_u(1)}^2.
\end{align*}
Further, for the tree root $v$ we have the deterministic identity
$$\langle \Row_v(M),Y\rangle^2=\delta_0^2\|\Row_v(M)\|_2^4.$$
Thus, using \eqref{eq: aux 087205295} and the above relations, we obtain
\begin{align*}
\E\|MY\|_2^2&=\E\sum\limits_{u\in V}\langle \Row_u(M),Y\rangle^2
=\E\sum\limits_{r\in \Z_2\cap[0,q]}\sum\limits_{u\in V_r}\langle \Row_u(M),Y\rangle^2\\
&\geq \delta_0^2\|\Row_v(M)\|_2^4+\sum\limits_{r\in \Z_2\cap[2,q]}\|\Row_v(M)\|_2^2(d'-1)^{r-1}
\big(\delta_{r-2}+\delta_r(d'-1)\big)^2.
\end{align*}
On the other hand, again in view of \eqref{eq: aux 087205295},
\begin{align*}
\E\|Y\|_2^2&=\sum\limits_{r\in \Z_2\cap[0,q]}\sum\limits_{u\in V_r}\E\|Y_u\|_2^2\\
&\leq \delta_0^2 \|\Row_v(M)\|_2^2+\sum\limits_{r\in \Z_2\cap[2,q]}\sum\limits_{u\in V_r}\delta_r^2 \mu_{P_u(0),P_u(1)}^2(\w d-1)\\
&\leq \delta_0^2 \|\Row_v(M)\|_2^2+\sum\limits_{r\in \Z_2\cap[2,q]} \delta_r^2 \|\Row_v(M)\|_2^2(\w d-1)^{r}.
\end{align*}
Now, we recursively define $\delta_r:=\frac{\delta_{r-2}\w d}{(\w d-1)(\|\Row_v(M)\|_2^2-\w d)}$, $r\in\Z_2\cap[2,q]$.
Note that this formula gives
$$\delta_r=\frac{\delta_{0}\w d^{r/2}}{(\w d-1)^{r/2}(\|\Row_v(M)\|_2^2-\w d)^{r/2}},\quad r\in\Z_2\cap[2,q].$$
Plugging in the values into the above estimate of $\E\|Y\|_2^2$, we get
$$\E\|Y\|_2^2\leq \delta_0^2\|\Row_v(M)\|_2^2\big(1+\rho^2+\rho^4+\dots+\rho^{q-1}\big)
=\delta_0^2\|\Row_v(M)\|_2^2\frac{1-\rho^{q+1}}{1-\rho^2},$$
where $\rho:=\frac{\w d}{\|\Row_v(M)\|_2^2-\w d}\leq \frac{1}{1+2\varepsilon}$.
On the other hand, using the same procedure for $\E\|MY\|_2^2$ and that $d'\leq \w d$, we get 
\begin{align*}
\E\|MY\|_2^2&\geq \delta_0^2\|\Row_v(M)\|_2^4+\sum\limits_{r\in \Z_2\cap[2,q]}\delta_0^2\|\Row_v(M)\|_2^2
\frac{(d'-1)^{r+1}}{(\w d-1)^{r}}
\big(\rho^{(r-2)/2}+\rho^{r/2}\big)^2\\
&\geq \delta_0^2\|\Row_v(M)\|_2^4
+\delta_0^2\|\Row_v(M)\|_2^2(1+\rho)^2\frac{(d'-1)^3}{(\w d-1)^2}\frac{1-\w \rho^{q-1}}{1-\w\rho^2},
\end{align*}
where $\w\rho:=\frac{d'-1}{\w d-1}\rho$. Taking the ratio and using that $\w\rho\leq \rho'$, we obtain
\begin{align*}
\frac{\E\|MY\|_2^2}{\E\|Y\|_2^2}&\geq \|\Row_v(M)\|_2^2(1-\rho^2)+
\frac{(d'-1)^3}{(\w d-1)^2}\frac{(1+\rho)^2(1-{\w\rho}^{q-1})}{(1-\rho^{q+1})(1-\w\rho^2)}(1-\rho^2)\\
&\geq \|\Row_v(M)\|_2^2(1-\rho^2)+(d'-1)(1+\rho)^2\Big[ \frac{(d'-1)^2}{(\w d-1)^2}\frac{(1-{\w\rho}^{q-1})(1-\rho)}{(1-\rho^{q+1})(1-\w\rho)}\Big].
\end{align*}
Choosing $d$ sufficiently large (depending on $\varepsilon$), we have $\w d-1\leq (1+2\delta)(d'-1)$ and $\rho\leq (1+2\delta)\w \rho$. Since $\rho\leq \frac{1}{1+2\varepsilon}$ then this implies that $(1-\rho)\geq (1-\delta/\varepsilon) (1-\w \rho)$. Replacing these in the above relation, we get 
$$
\frac{\E\|MY\|_2^2}{\E\|Y\|_2^2}\geq \|\Row_v(M)\|_2^2(1-\rho^2)+(d'-1)(1+\rho)^2\Big[ \frac{(1-{\w\rho}^{q-1})(1-\delta/\varepsilon)}{(1+2\delta)^2(1-\rho^{q+1})}\Big]
$$
Now, note that for every sufficiently large $q$ and sufficiently small $\delta$ (depending on $\varepsilon$), we get
$$\frac{\E\|MY\|_2^2}{\E\|Y\|_2^2}\geq \|\Row_v(M)\|_2^2(1-\rho^2)+(d'-1)(1+\rho)^2-\varepsilon d'/16.$$
Further, using the definition of $\rho$, it is not difficult to check that
\begin{align*}
\|\Row_v(M)\|_2^2(1-\rho^2)+(d'-1)(1+\rho)^2
&=\frac{\|\Row_v(M)\|_2^4}{\|\Row_v(M)\|_2^2-\widetilde d}+\big(d'-1-\widetilde d\big)\frac{\|\Row_v(M)\|_2^4}{(\|\Row_v(M)\|_2^2-\widetilde d)^2}\\
&\geq \frac{\|\Row_v(M)\|_2^4}{\|\Row_v(M)\|_2^2-\widetilde d} -\varepsilon d'/16,
\end{align*}
where we used the assumption on $\|\Row_v(M)\|_2^2$ and that $\delta$ is sufficiently small in terms of $\varepsilon$.  
Thus, setting $y:=\frac{Y}{(\E\|Y\|_2^2)^{1/2}}$, we get
$$\E(\|M\|^2\|y\|_2^2)\geq\E\|My\|_2^2\geq  \frac{\|\Row_v(M)\|_2^4}{\|\Row_v(M)\|_2^2-\widetilde d}
-\varepsilon d'/8.$$
Denoting by $\Event$ the event that $\big|\|M\|-(\E\|M\|^2)^{1/2}\big|\geq (\varepsilon/4) (\E\|M\|^2)^{1/2}$,
we obtain from the above
\begin{align*}
(1+\varepsilon/4)^2\E\|M\|^2&=\E\big((1+\varepsilon/4)^2(\E\|M\|^2)\|y\|_2^2\big)\\
&\geq \E\big((1+\varepsilon/4)^2(\E\|M\|^2)\|y\|_2^2{\bf 1}_{\Event^c}\big)\\
&\geq\E\big(\|M\|^2\|y\|_2^2{\bf 1}_{\Event^c}\big)\\
&\geq \frac{\|\Row_v(M)\|_2^4}{\|\Row_v(M)\|_2^2-\widetilde d}-\varepsilon d'/8
-\E\big(\|M\|^2\|y\|_2^2{\bf 1}_{\Event}\big).
\end{align*}
Applying the two-sided deterministic inequality for $\|Y\|_2^2$ from the beginning of the proof and using that $\w d\leq (\Vert \Row_v(M)\Vert_2^2-\w d)$, we get
that
$$\|y\|_2^2\leq \delta_0^{-2}\sum\limits_{r\in \Z_2\cap[0,q]}\delta_r^2 (\widetilde d-1)^r h^{r}
\leq q h^{q-1}$$
everywhere on the probability space. Hence,
$$\E\big(\|M\|^2\|y\|_2^2{\bf 1}_{\Event}\big)
\leq q h^{q-1} \E\big(\|M\|^2{\bf 1}_{\Event}\big).$$
Further, applying Talagrand's concentration inequality (Theorem~\ref{th: tal}) to $\|M\|$, we get that $\Vert M\Vert -\E\, \Vert M\Vert$ is a subgaussian variable. Using this, and applying 
Talagrand's inequality again to bound $\Prob(\Event)$, we can write 
$$
q h^{q-1} \E\big(\|M\|^2{\bf 1}_{\Event}\big)\leq q h^{q-1}  \sqrt{\E\, \Vert M\Vert^4}\cdot \sqrt{\Prob(\Event)} \leq \frac{\varepsilon}{8} \E\, \Vert M\Vert^2,
$$
as long as
$d=d(h,\varepsilon)$ is chosen sufficiently large. 
Since $d'\leq \max_{u\leq n} \E\, \Vert \Row_u(M)\Vert_2^2\leq \E\, \Vert M\Vert^2$, we get from the above that 
$$
(1+\varepsilon/4)^2\E\|M\|^2\geq \frac{\|\Row_v(M)\|_2^4}{\|\Row_v(M)\|_2^2-\widetilde d}
-\frac{\varepsilon}{4} \E\, \Vert M\Vert^2,
$$
implying that 
$$
\E\|M\|^2\geq \frac{(1-\varepsilon)\|\Row_v(M)\|_2^4}{\|\Row_v(M)\|_2^2-\widetilde d}.
$$
It remains to apply Talagrand's inequality to $\|M\|$ in a similar manner to above to get the result.
\end{proof}

\medskip
\begin{proof}[Proof of Theorem~\ref{th: lower top}]
In the case $\frac{n p_n}{\log n}\to 0$, Lemma~\ref{l: sublog convergence} implies that 
$$
\frac{\rho_n}{\max\limits_{i\leq n}\|\Row_i(W_n)\|_2}\underset{n\to\infty}{\overset{\Prob}{\longrightarrow}} 1.
$$
On the other hand, Lemma~\ref{l: maxnorm-restrict} together with a simple observation that,
under the sparsity assumption, with high probability supports of any two rows of the matrix intersect on
at most a constant number of coordinates, implies that for any fixed $k$ and any $\varepsilon>0$,
$\Prob\big\{\lambda_{|k|}(W_n)\geq (1-\varepsilon)\max\limits_{i\leq n}\|\Row_i(W_n)\|_2\big\}$
converges to $1$ with $n\to\infty$. Then the theorem follows.

In the remaining, we assume that
$n p_n\geq \alpha \log n$ for some $\alpha>0$ and all large $n$.
Let $k\in \N$ be fixed. For any $n\geq 1$, denote the entries of $W_n$ by $w_{ij}^n$.
We will use the representation $w_{ij}^n=b_{ij}^n\xi_{ij}^n$, where
$\xi_{ij}^n$ is a copy of $\xi$, and $b_{ij}^n$ is an independent Bernoulli random variable
with probability of success $p_n$.

It is not difficult to check that it is sufficient to prove the statement for $\xi$ having an absolutely continuous
distribution (for example, the reduction can be performed by convolving the original distribution of $\xi$
with a uniform distribution on $[-\kappa,\kappa]$ for a small $\kappa$; note that in this case the maximum Euclidean
norm of the rows of the original and perturbed matrix are close with high probability).

Define the event $$\Event_n:=\Big\{\sum\limits_{j=1}^n b_{ij}^n\leq \log^{4}n\mbox{ for all $i\leq n$}\Big\},$$
and note that in view of our assumptions on the $p_n$'s, a direct application of Bernstein inequality (Lemma~\ref{l:bernstein}) implies
$\Prob(\Event_n)\geq 1-n^{-\log n},$ for all sufficiently large $n$. 

Let $r=r_{\ref{l: rows restricted}}$ be taken from Lemma~\ref{l: rows restricted},
for $\delta:=\varepsilon/8$. Then, according to Lemma~\ref{l: rows restricted} and our condition on $p_n$,
setting
$$\eta_n:=\max\limits_{1\leq i\leq rn}\Big\|\sum\limits_{j=\lfloor rn\rfloor+1}^n w_{ij}^n e_j\Big\|_2,$$
for every sufficiently large $n$ we get
\begin{equation}\label{eq: aux --498523-592835}
\Prob\big\{\max\limits_{1\leq i\leq n}\|\Row_i(W_n)\|_2\geq (1+\varepsilon/8)\eta_n\big\}\leq \varepsilon/8.
\end{equation}
Fix for a moment any $n$. Let
$U_n$ be a median of $\eta_n$ i.e.\ a number such that
$\Prob\{\eta_n\geq U_n\}=\Prob\{\eta_n\leq U_n\}= 1/2$ and  
set 
$$
\tau_n:=\Prob\Big\{\big\|\sum_{j=\lfloor rn\rfloor+1}^n w_{1j}^n e_j\big\|_2\geq U_n\Big\},\quad \text{ and }\quad  m_n:=\min(\lfloor rn\rfloor,\lfloor k/\tau_n\rfloor).$$
Now, take any $k$--subset $I$ of $[m_n]$, and define the event 
$$
\Event_{n,I}:=\Big\{\Big\|\sum\nolimits_{j=\lfloor rn\rfloor+1}^n w_{ij}^n e_j\Big\|_2\geq U_n\mbox{ for all $i\in I$ and
$\sum\limits_{i\in I}b_{ij}^n\leq 1$ for all $j\geq \lfloor rn\rfloor+1$}\Big\}.$$
By independence, the probability of $\Event_{n,I}$ is bounded above by $\tau_n^k$.
Further, the above estimate of the probability of $\Event_n$, together with the observation
$\tau_n\geq 1/(2n)$, yields
\begin{align*}
\Prob&\Big(\Event_n\;\Big\vert\;\Big\{\Big\|\sum\nolimits_{j=\lfloor rn\rfloor+1}^n
w_{ij}^n e_j\Big\|_2\geq U_n\mbox{ for all $i\in I$}\Big\}\Big) \\
&\geq1-\frac{\Prob\big(\Event_n^c\big)}{\Prob\Big\{\Big\|\sum\nolimits_{j=\lfloor rn\rfloor+1}^n
w_{ij}^n e_j\Big\|_2\geq U_n\mbox{ for all $i\in I$}\Big\}}\\
&\geq 1-\frac{n^{-\log n}}{\tau_n^k}\geq 3/4.
\end{align*}
On the other hand, we can similarly write 
\begin{align*}
\Prob\Big\{&\mbox{$\sum\limits_{i\in I}b_{ij}^n\leq 1\;\;$ $\forall\;j\geq \lfloor rn\rfloor+1$}\;\Big\vert\;
\Big\|\sum\nolimits_{j=\lfloor rn\rfloor+1}^n w_{ij}^n e_j\Big\|_2\geq U_n\mbox{
and $\sum\nolimits_{j=1}^n b_{ij}^n\leq \log^{4}n\;\;$  $\forall\;i\in I$}\Big\}\\
&\geq 1- 
\frac{\Prob\Big\{\mbox{$\exists\;j\geq \lfloor rn\rfloor+1,\; \sum\limits_{i\in I}b_{ij}^n\geq  2$}\Big\}}{\Prob\Big\{\Big\|\sum\nolimits_{j=\lfloor rn\rfloor+1}^n w_{ij}^n e_j\Big\|_2\geq U_n\mbox{ and $\sum\nolimits_{j=1}^n b_{ij}^n\leq \log^{4}n\;\;$  $\forall\;i\in I$}\Big\}}\\
\, \\
&\geq 1-\frac{n.{k\choose 2} p_n^2}{\tau_n^k-n^{-\log n}} \geq 3/4. 
\end{align*}
Putting the above estimates together, we deduce that $\Prob(\Event_{n,I})\geq \tau_n^k/2$ assuming that $n$ is sufficiently large.

Let $G_n$ be the random graph whose adjacency matrix is given by the matrix $B_n=(b_{ij}^n)$.
We set $\widetilde \delta:=\min\big(\delta_{\ref{l: main lower}}(h,\varepsilon/(32ke^k))/4,
\varepsilon/(32e^k)\big)$ and $q:=q_{\ref{l: main lower}}(h,\varepsilon/(32ke^k))$.
Define two more events
$$\Event_{n,I}':=\big\{\mbox{for any $i_1\neq i_2\in I$, the distance between $i_1$ and $i_2$ in $G_n$
is at least $2q+1$}\big\},$$
and
\begin{align*}
\Event_{n,I}'':=\big\{&\mbox{$q$--neighborhood of every vertex $i\in I$ in $G_n$ is a tree, and}\\
&\mbox{the degrees of all vertices of the trees, except for the roots}\\
&\mbox{and leaves, are in the range $[(1-\widetilde\delta )n p_n,(1+\widetilde\delta)n p_n]$}\big\}.
\end{align*}
Then from the above definition of $\Event_{n,I}$ (and the lower bound for the probability) and from Lemma~\ref{l: graph str}
we obtain that for all sufficiently large $n$,
\begin{align}
\Prob\big(\Event_{n,I}'\cap
\Event_{n,I}\cap \Event_n\big)
&\geq
\Prob\Big(\Event_{n,I}'\cap
\Event_{n,I}\cap \Big\{\sum\nolimits_{j=\lfloor rn\rfloor+1}^n b_{ij}^n\leq \log^{4}n\mbox{ for all $i\in I$}\Big\}\Big)
-n^{-\log n}\nonumber\\
&\geq \big(1-n^{-1/8}\big)\Prob(\Event_{n,I}\cap \Event_n)-n^{-\log n}\nonumber\\
&\geq \big(1-n^{-1/16}\big)\Prob(\Event_{n,I}\cap \Event_n);\label{eq: aux 20872652=095=23}
\end{align}
similarly,
\begin{equation}\label{eq: aux 0982260-986}
\Prob\big(\Event_{n,I}''\cap\Event_{n,I}\cap \Event_n\big)\geq \big(1-\varepsilon/(32e^k)\big)\Prob(\Event_{n,I}\cap \Event_n).
\end{equation}
Observe that the intersection $\Event_{n,I}\cap\Event_{n,I}'\cap\Event_{n,I}''\cap \Event_n$
is measurable with respect to the algebra generated by entries
$\big\{w_{ij}^n:\;j\geq \lfloor rn\rfloor+1\mbox{ and }i\in I\big\}$ and by matrix $B_n$.

Condition for a moment on any realization of $B_n$ and of $\big\{a_{ij}:\;j\geq \lfloor rn\rfloor+1,\;i\in I\big\}$
which belongs to the intersection $\Event_{n,I}\cap\Event_{n,I}'\cap\Event_{n,I}''\cap \Event_n$.
For each $i\in I$, let $\tree_i$ denote the tree of depth $q$ (in $G_n$) rooted in $i$.
Further, for each $i$ denote by $B_n(i)$ the $n\times n$ adjacency matrix of $\tree_i$ (where we treat every number
$j$ which is not a node/leaf of $\tree_i$ as an isolated vertex),
and by $W_n(i)$ the entry-wise product of $B_n(i)$ with $(\xi_{ij}^n)$.
Since, by our assumption (see definition of $\Event_{n,I}'$),
the trees do not have common vertices, the sum $\sum_{i\in I}W_n(i)$
is a permutation of a block diagonal matrix with $|I|=k$ blocks having spectral norms $\|W_n(i)\|$, $i\in I$.
Hence, $|\lambda_{|k|}(\sum_{i\in I}W_n(i))|\geq \min_{i\in I}\|W_n(i)\|$, and, moreover,
since $\sum_{i\in I}W_n(i)$ is a compression of $W_n$, the Cauchy interlacing theorem implies that
$|\lambda_{|k|}(W_n)|\geq |\lambda_{|k|}(\sum_{i\in I}W_n(i))|$.
On the other hand, applying Lemma~\ref{l: main lower} and taking into account our choice of parameters,
we obtain that, as long as $n$ is sufficiently large,
\begin{align*}
\Prob\bigg\{\|W_n(i)\|\geq \frac{(1-\varepsilon/8)\|\Row_i(W_n)\|_2^2\,{\bf 1}_{\{\|\Row_i(W_n)\|_2^2
\geq 2(1+\varepsilon/32)\widetilde d\}}}{\sqrt{\|\Row_i(W_n)\|_2^2-\widetilde d}}\;\Big\vert\;\mathcal R\bigg\}
\geq 1-\varepsilon/(32ke^k),\quad i\in I,
\end{align*}
where $\widetilde d:=(1+\widetilde\delta)n p_n$ and ``$\mathcal R$'' is meant to emphasize that
we are conditioning on an appropriate realization of $B_n$ and $\big\{a_{ij}:\;j\geq \lfloor rn\rfloor+1,\;i\in I\big\}$.
By our definition of $\Event_{n,I}$, we have $\|\Row_i(W_n)\|_2^2\geq U_n^2$, $i\in I$, assuming the conditioning.
Hence, using also the choice of $\w d$, we get
\begin{align*}
\Prob\bigg\{\|W_n(i)\|&\geq \frac{(1-\varepsilon/8)U_n^2\,{\bf 1}_{\{U_n^2
\geq 2(1+\varepsilon/8)np_n\}}}{\sqrt{U_n^2-np_n}}\;\Big\vert\;\mathcal R\bigg\}\\
&=
\Prob\bigg\{\|W_n(i)\|\geq \frac{(1-\varepsilon/8)U_n^2\,{\bf 1}_{\{U_n^2
\geq 2(1+\varepsilon/32)\widetilde d\}}}{\sqrt{U_n^2-\widetilde d}}\;\Big\vert\;\mathcal R\bigg\}
\geq 1-\varepsilon/(32ke^k),\quad i\in I.
\end{align*}
Thus, setting $\rho_n':=\big(\sqrt{U_n^2-n p_n}
+n p_n/\sqrt{U_n^2-n p_n}\big){\bf 1}_{\{U_n^2\geq 2(1+\varepsilon/8)n p_n\}}$ we get
$$
\Prob\big\{|\lambda_{|k|}(W_n)|\geq (1-\varepsilon/8)\rho_n'\;\vert\;
\Event_{n,I}\cap\Event_{n,I}'\cap\Event_{n,I}''\cap \Event_n\big\}\geq 1-\varepsilon/(32e^k),
$$
whence, using \eqref{eq: aux 20872652=095=23}--\eqref{eq: aux 0982260-986},
$$
\Prob\big\{|\lambda_{|k|}(W_n)|\geq (1-\varepsilon/8)\rho_n'\;\vert\;
\Event_{n,I}\cap\Event_n\big\}\geq
\big(1-\varepsilon/(32e^k)\big)\frac{\Prob(\Event_{n,I}\cap\Event_{n,I}'\cap\Event_{n,I}''\cap \Event_n)}
{\Prob(\Event_{n,I}\cap\Event_n)}\geq 1-\frac{\varepsilon}{e^k}.
$$
Define the event $\widetilde\Event_n:=\{|\lambda_{|k|}(W_n)|\geq (1-\varepsilon/8)\rho_n'\}$. 
Then, by the above, $\Prob(\widetilde\Event_n^c\cap \Event_{n,I}\cap \Event_n)
\leq \varepsilon e^{-k}\Prob(\Event_{n,I}\cap \Event_n)$
for all $I\subset[m_n]$ with $|I|=k$.
However,
$$\sum\limits_{I\subset[m_n],|I|=k}\Prob(\Event_{n,I})
\leq\frac{m_n^k}{k!}\tau_n^k\leq \frac{k^k}{k!},$$
whence
\begin{equation}\label{eq: aux 2-9845y05634029}
\Prob\Big(\widetilde\Event_n^c\cap \Event_n\cap \bigcup\limits_{I\subset[m_n],|I|=k}\Event_{n,I}\Big)\leq \varepsilon.
\end{equation}
Further, we note that, essentially repeating the argument we used to estimate $\Prob(\Event_{n,I})$ from below,
we get that
\begin{align*}
\Prob&\Big(\bigcup\limits_{I\subset[m_n],|I|=k}\Event_{n,I}\Big)\\
&\geq \frac{1}{2}\Prob\Big\{\mbox{There is $I\subset[m_n]$ with $|I|=k$ such that }
\Big\|\sum\nolimits_{j=\lfloor rn\rfloor+1}^n a_{ij}^n e_j\Big\|_2\geq U_n\mbox{ for all $i\in I$}\Big\}\\
&= \frac{1}{2}\sum\limits_{\ell=k}^{m_n}{m_n \choose \ell}\,\tau_n^\ell(1-\tau_n)^{m_n -\ell},
\end{align*}
while the definition of $\tau_n$ implies
$(1-\tau_n)^{m_n}\geq 1/2$ and $m_n \tau_n\geq 1/2$.
These relations together give
$$
\Prob\Big(\bigcup\limits_{I\subset[m_n],|I|=k}\Event_{n,I}\Big)
\geq f(k)
$$
for some strictly positive function of $k$.
Applying this estimate with \eqref{eq: aux 2-9845y05634029}
and the lower bound for $\Prob(\Event_n)$, we obtain
$$\Prob(\widetilde\Event_n)
\geq \Prob\Big(\Event_n\cap \bigcup\limits_{I\subset[m_n],|I|=k}\Event_{n,I}\Big)
-
\Prob\Big(\widetilde\Event_n^c\cap \Event_n\cap \bigcup\limits_{I\subset[m_n],|I|=k}\Event_{n,I}\Big)
\geq f(k)-n^{-\log n}-\varepsilon.$$
Thus, for sufficiently small $\varepsilon$ and assuming that $n$ is large, we have
$\Prob(\widetilde\Event_n)\geq f(k)/2$. 
Note that $|\lambda_{|k|}(W_n)|$ is a $1$--Lipschitz function of the matrix $W_n$, so that Talagrand concentration inequality
(Theorem~\ref{th: tal}) gives for such $\varepsilon$ and all large enough $n$,
$$
\Prob\big\{|\lambda_{|k|}(W_n)|\geq (1-\varepsilon/7)\rho_n'\big\}\geq 1-\varepsilon/10.
$$
Indeed, in the above we applied the geometric form of the concentration of measure phenomenon as if $\lambda_{|k|}(W_n)|\leq (1-\varepsilon/7)\rho_n'$, then 
necessarily $W_n$ is at distance (with respect to the Frobenius norm) at least $\varepsilon \rho_n'/56$ from $\w \Event$. 

As a final step of the proof, we replace $\rho'$ in the last relation with
$$\rho=\sqrt{\max\big(\max\limits_{i\leq n}\|\Row_i(W_n)\|_2^2-np_n,n p_n\big)}
+\frac{n p_n}{\sqrt{\max\big(\max\limits_{i\leq n}\|\Row_i(W_n)\|_2^2-np_n,n p_n\big)}}.$$
Observe that, in view of \eqref{eq: aux --498523-592835} and Talagrand's concentration inequality applied to $\eta_n$, we have
$$
\Prob\big\{(1-\varepsilon/8)U_n\leq
\max\limits_{1\leq i\leq n}\|\Row_i(W_n)\|_2\leq (1+\varepsilon/7)U_n\big\}\geq 1-\varepsilon/7,
$$
and
$$
\Prob\big\{{\bf 1}_{\{U_n^2
\geq 2(1+\varepsilon/8)\widetilde d\}}
\geq {\bf 1}_{\{\max\limits_{1\leq i\leq n}\|\Row_i(W_n)\|_2^2
\geq 2(1+\varepsilon)n p_n\}}\big\}\geq 1-\varepsilon/7.
$$
Hence, we get
$$
\Prob\Big\{|\lambda_{|k|}(W_n)|\geq (1-\varepsilon)\rho_n {\bf 1}_{\{\max\limits_{1\leq i\leq n}\|\Row_i(W_n)\|_2^2
\geq 2(1+\varepsilon)n p_n\}}\Big\}\geq 1-\varepsilon/3.
$$
It remains to note that on event $\big\{\max\limits_{1\leq i\leq n}\|\Row_i(W_n)\|_2^2
< 2(1+\varepsilon)n p_n\big\}$ we have
$\rho_n\leq (2+\varepsilon)\sqrt{n p_n}$, whence for all large $n$
$$
\Prob\Big\{|\lambda_{|k|}(W_n)|\geq (1-\varepsilon)\rho_n {\bf 1}_{\{\max\limits_{1\leq i\leq n}\|\Row_i(W_n)\|_2^2
< 2(1+\varepsilon)n p_n\}}\Big\}\geq 1-\varepsilon/3,
$$
just using that the spectral distribution of $A_n$, rescaled by $1/\sqrt{n p_n}$, converges to the semi-circle
law.
Combining the last two relations, we get that for all $\varepsilon\in(0,f(k)/4)$ and all
sufficiently large $n$,
$$
\Prob\big\{|\lambda_{|k|}(W_n)|\geq (1-\varepsilon)\rho_n\big\}\geq 1-\varepsilon.
$$
The result follows.
\end{proof}

\section{Outliers in the spectrum of the Erd\H os--Renyi graphs}\label{s: outliers-erdos-renyi}

In this section, we prove Corollary~C which is simply an application of the main theorems of this paper to the adjacency matrices
of the Erd\H os--Renyi graphs. 
The next proposition follows from \cite[Theorem~3.1]{book-bollobas}, we record it for ease of future references. 

\begin{prop}\label{prop: maxdegree}
Let $\mathcal{G}(n,p_n)$ be the undirected Erd\H{o}s--Renyi graph with the edge probability equal to $p_n$. Suppose further that $p_n\to 0$ and $n p_n\to\infty$ with $n$. 
Then 
$$
\frac{\max_{i\leq n} \degree(i)}{enp_n \exp\Big[\mathcal{W}_0\big( \frac{\log n-np_n}{enp_n}\big)\Big]} \underset{n\to\infty}{\overset{\Prob}{\longrightarrow}} 1, 
$$
where $\mathcal{W}_0$ denotes the main branch of the Lambert function.
\end{prop}
\begin{proof}
Denote $\gamma_n:= enp_n \exp\Big[\mathcal{W}_0\big( \frac{\log n-np_n}{enp_n}\big)\Big]$. 
For any $\ell\leq n$, denote $\alpha_{\ell}(n):= n{n \choose \ell} p_n^{\ell}(1-p_n)^{n-\ell}$. 
In \cite[Theorem~3.1]{book-bollobas}, it is  stated that if $\ell=\ell(n)\leq n/2$ is a positive integer sequence with
$$
\lim_{n} \alpha_\ell(n)=0,
$$
then $\lim_n \Prob\big\{ \max_{i\leq n} \degree(i)\leq \ell\big\}=1$. On the other hand, if 
$$
\lim_{n} \alpha_\ell(n)=\infty,
$$
then $\lim_n \Prob\big\{ \max_{i\leq n} \degree(i)\geq \ell\big\}=1$. 
Let $\varepsilon\in (0,1)$. To get that $\lim_n \Prob\big\{ \max_{i\leq n} \degree(i)\leq \lceil(1+\varepsilon)\gamma_n\rceil\big\}=1$, it is sufficient to verify that 
$\lim_{n} \alpha_{\lceil (1+\varepsilon)\gamma_n\rceil} (n)=0$. 
Notice that by our assumptions $\gamma_n=o(n)$ and that $\gamma_n\to \infty$ with $n$. Therefore, using Stirling's approximation
and the estimating the rate of change of $\mathcal{W}_0(z)$ near $z=-1/e$, we get
$$
\alpha_{\lceil (1+\varepsilon)\gamma_n\rceil} (n)\sim \frac{ne^{-np_n}}{\sqrt{2\pi \lceil (1+\varepsilon)\gamma_n\rceil}} \Big(\frac{e np_n}{\lceil (1+\varepsilon)\gamma_n\rceil}\Big)^{\lceil (1+\varepsilon)\gamma_n\rceil}.
$$
Thus, using that $ \mathcal{W}_0(z)e^{\mathcal{W}_0(z)}=z$, we obtain
\begin{align*}
\lim_n \alpha_{\lceil (1+\varepsilon)\gamma_n\rceil} (n)&\leq C
\lim_n  \frac{ne^{-np_n}}{\sqrt{2\pi (1+\varepsilon)\gamma_n}} \exp\Big[-(1+\varepsilon) (\log n-np_n)\Big]\,  (1+\varepsilon)^{ -(1+\varepsilon) \gamma_n}\\
& \leq C\lim_n   \frac{1}{\sqrt{2\pi \lfloor (1+\varepsilon)\gamma_n\rfloor}}\exp\Big[-\varepsilon (\log n-np_n)\Big]\,  (1+\varepsilon)^{ -(1+\varepsilon) \gamma_n}.
\end{align*}
To evaluate the last limit, note that since $\mathcal{W}_0(\cdot)\geq -1$, then $\gamma_n\geq np_n$ and we get 
$$
\lim_n \alpha_{\lceil (1+\varepsilon)\gamma_n\rceil} (n)\leq 
C \lim_n   \frac{1}{\sqrt{2\pi \lfloor (1+\varepsilon)\gamma_n\rfloor}}\exp\Big[-\varepsilon \log n-np_n\big( (1+\varepsilon)\log(1+\varepsilon)-\varepsilon\big) \Big].
$$
Since $(1+\varepsilon)\log(1+\varepsilon)-\varepsilon>0$ 
we deduce that the above limit is zero.

To prove the lower bound, we first note that if $\frac{np_n}{\log n}\to \infty$, then by standard concentration inequalities, we have 
$$
\lim_n \Prob\big\{ \max_{i\leq n} \degree(i)\geq (1-\varepsilon)np_n\big\}=1\quad\mbox{for all $\varepsilon>0$}.
$$
Moreover, in this case, we also have $\frac{\gamma_n}{np_n}\to 1$ with $n$. 

Therefore, we may suppose that $\frac{np_n}{\log n}\leq C\in (0,\infty)$ for all $n$.
It follows that there exists $\varepsilon_0\in (0,1/2)$ such that for all sufficiently large $n$, we have $\gamma_n\geq (1+\varepsilon_0)np_n$. Let $\varepsilon\leq \varepsilon_0^2$.  To show that $\lim_n \Prob\big\{ \max_{i\leq n} \degree(i)\geq \lfloor (1-\varepsilon)\gamma_n\rfloor\big\}=1$,
we will verify that 
$\lim_{n} \alpha_{\lfloor (1-\varepsilon)\gamma_n\rfloor} (n)=\infty$. 
As before, applying Stirling's approximation formula, we get 
\begin{align*}
\lim_n \alpha_{\lfloor (1-\varepsilon)\gamma_n\rfloor} (n)& \geq c\lim_n  \frac{1}{\sqrt{2\pi \lceil (1-\varepsilon)\gamma_n\rceil}}\exp\Big[\varepsilon (\log n-np_n)\Big]\,  (1-\varepsilon)^{ -(1-\varepsilon) \gamma_n}\\
&\geq  c\lim_n \exp\Big[-\gamma_n \Big((1-\varepsilon) \log (1-\varepsilon) +\frac{\varepsilon}{1+\varepsilon_0} \Big) \Big],
\end{align*}
where in the last inequality we used that $\gamma_n\geq (1+\varepsilon_0)np_n$,
and that $n^\varepsilon\gg \gamma_n$. Now, in view of the inequality $\log (1-\varepsilon)\leq -\varepsilon$, we get 
\begin{align*}
\lim_n \alpha_{\lfloor (1-\varepsilon)\gamma_n\rfloor} (n)&\geq  c\lim_n \exp\Big[\varepsilon \gamma_n \Big((1-\varepsilon) -\frac{1}{1+\varepsilon_0} \Big) \Big].
\end{align*}
By the choice of $\varepsilon$ and $\varepsilon_0$, we have $\frac{1}{1+\varepsilon_0}< 1-\varepsilon$ which implies that the above limit is infinite and finishes the proof.  
\end{proof}

\medskip

\begin{proof}[Proof of Corollary~C]
The first part of the corollary follows immediately by combining the above Proposition~\ref{prop: maxdegree}
with Theorem~B. Hence, it only remains to explicitly compute the point of the phase transition.

Suppose that $\liminf\limits_{n\to\infty} \frac{n p_n}{\log n}\geq\frac{1}{\log (4/e)}$. For any $\varepsilon\in (0, 1/2)$,  we have $n p_n\geq \frac{1-\varepsilon}{\log (4/e)}\log n$ for all sufficiently large $n$. 
If follows from Proposition~\ref{prop: maxdegree} that with probability going to one with $n$, we have 
$$
\frac{\max_{i\leq n} \degree(i)}{np_n} \leq (1+\varepsilon) e\, \exp\Big[\mathcal{W}_0\big( \frac{\log n-np_n}{enp_n}\big)\Big] 
\leq  (1+\varepsilon) e\, \exp\Big[\mathcal{W}_0\big( \frac{2\log (2/e) +\varepsilon}{e(1-\varepsilon)}\big)\Big]. 
$$
Noting that  $\mathcal{W}_0\big( \frac{2\log (2/e)}{e}\big)= \log(2/e)$ and using an approximation of $\mathcal{W}_0$, we deduce that with probability going to one with $n$, we have
$$
\frac{\max_{i\leq n} \degree(i)}{np_n} \leq 2(1+c\varepsilon), 
$$
for some universal constant $c$. Applying Theorem~B, we get the first part of Corollary~C.

Assume now that $\limsup\limits_{n\to\infty} \frac{n p_n}{\log n}<\frac{1}{\log (4/e)}$.
This implies that there is $\varepsilon_0\in (0,1/2)$ such that for all sufficiently large $n$ 
we have $\frac{n p_n}{\log n}<\frac{1-\varepsilon_0}{\log (4/e)}$. It follows from Proposition~\ref{prop: maxdegree} that for any $\varepsilon\in (0,1)$,  we have 
$$
\frac{\max_{i\leq n} \degree(i)}{np_n} \geq (1-\varepsilon) e\, \exp\Big[\mathcal{W}_0\big( \frac{\log n-np_n}{enp_n}\big)\Big]
\geq (1-\varepsilon) e\, \exp\Big[\mathcal{W}_0\big( \frac{2\log (2/e) +\varepsilon_0}{e(1-\varepsilon_0)}\big)\Big],
$$
with probability going to one with $n$. 
As before, this implies that  for any $\varepsilon\in (0,1)$, with probability going to one with $n$ we have 
$$
\frac{\max_{i\leq n} \degree(i)}{np_n} 
\geq 2(1-\varepsilon)(1+c'\varepsilon_0),
$$
where $c'>0$ is a universal constant. Therefore, for any $\varepsilon \leq c'\varepsilon_0/6$, with probability going to one with $n$ we have
$$
\frac{\max_{i\leq n} \degree(i)}{np_n} \geq 2(1+\varepsilon).
$$
An application of Theorem~\ref{th: lower top} gives the result.
\end{proof}

\medskip

\begin{Remark} 
We note that the phase transition point can be computed without using Proposition~\ref{prop: maxdegree}
and relying on completely standard estimates. We provide an argument below.

We let $A_n$ be the adjacency matrix of $\mathcal{G}(n,p_n)$ and assume that $p_n\to 0$ and $n p_n\to\infty$ with $n$. 
First, we consider the case when $\liminf\limits_{n\to\infty} \frac{n p_n}{\log n}\geq\frac{1}{\log (4/e)}$.
Take any $\varepsilon\in (0, 1/4)$ and note that by our assumption,  we have $n p_n\geq \frac{1-\varepsilon}{\log (4/e)}\log n$ for all sufficiently large $n$. 
Applying Bennett's inequality (see, for example, \cite[Theorem~2.9]{BLM}), we get 
$$
\Prob\big\{\|\Row_1(A_n)\|_2^2\geq (2+\varepsilon)n p_n\big\} \leq \exp\Big(-np_n H\big( 1+\varepsilon\big)\Big),
$$
where $H$ is defined by $H(x)=(1+x)\log (1+x)- x$. Now, it is easy to check that 
$$
(1-\varepsilon)\, H(1+\varepsilon) \geq (1+\varepsilon/4)\log (4/e),
$$
for any $0\leq \varepsilon\leq 1/4$.  This, together with the condition on $np_n$, implies that 
$$
\Prob\big\{\|\Row_1(A_n)\|_2^2\geq (2+\varepsilon)n p_n\big\} \leq \frac{1}{n^{1+\varepsilon/4}}.
$$
Thus, for all $\varepsilon\in (0,1/4)$ we have
$$
\lim\limits_{n\to\infty }\big(n\Prob\big\{\|\Row_1(A_n)\|_2^2\geq (2+\varepsilon)n p_n\big\}\big)=0,
$$
implying that
$$\frac{\rho_n}{2\sqrt{n p_n}}=\frac{\sqrt{\max(\max\nolimits_{i}\|\Row_i(A_n)\|_2^2-np_n,n p_n)}}{2\sqrt{n p_n}}
+\frac{\sqrt{n p_n}}{2\sqrt{\max(\max\nolimits_{i}\|\Row_i(A_n)\|_2^2-np_n,n p_n)}}$$
converges in probability to $1$ when $n$ tends to infinity.
Applying Theorem~B, we get the case of ``no non-trivial outliers''. 

\medskip

Next, assume that $\limsup\limits_{n\to\infty} \frac{n p_n}{\log n}<\frac{1}{\log (4/e)}$.
This implies that there is $\varepsilon_0\in (0,1/2)$ such that for all sufficiently large $n$ 
we have $\frac{n p_n}{\log n}<\frac{1-\varepsilon_0}{\log (4/e)}$. Let $\varepsilon= \varepsilon_0/12$. 
Denoting the entries of $A_n$ by $a_{ij}^n$,
standard estimates on the tail of the binomial distribution (see for instance \cite[Lemma~4.7.2]{ash-book}) imply 
$$
\Prob\Big\{\displaystyle\sum_{j=\lfloor \varepsilon n\rfloor+1}^n a_{1j}^n\geq (2+\varepsilon)n p_n\Big\}
\geq \frac{1}{\sqrt{8n\gamma (1-\gamma)}}\, \exp\Big(-n\Big( \gamma \log \frac{\gamma}{p_n}+ (1-\gamma) \log \frac{1-\gamma}{1-p_n}\Big)\Big),$$
where $\gamma= \frac{(2+\varepsilon)p_n}{(1-\varepsilon)}$. Now it is easy to check that 
$$
\gamma \log \frac{\gamma}{p_n}\leq p_n\Big(\log 4 + 3\varepsilon\Big) \quad \text{and}\quad  (1-\gamma) \log \frac{1-\gamma}{1-p_n}\leq -p_n (1+\varepsilon),
$$
where we have used that $p_n$ is small enough for large $n$. Using this together with the condition on $np_n$, we deduce that 
$$
\Prob\Big\{\displaystyle\sum_{j=\lfloor \varepsilon n\rfloor+1}^n a_{1j}^n\geq (2+\varepsilon)n p_n\Big\}
\geq \frac{c}{\sqrt{\log n}}\, \exp\Big(-(1-\varepsilon_0)(1+6\varepsilon) \log n\Big)
\geq \frac{c}{\sqrt{\log n}}\, n^{c\varepsilon-1} ,$$
for an appropriate universal constant $c$. 

Therefore, we can write 
\begin{align*}
\Prob\big\{\max\limits_{i\leq n}\|\Row_i(A_n)\|_2^2<(2+\varepsilon)n p_n\big\}
&\leq \Prob\Big\{\sum\nolimits_{j=\lfloor \varepsilon n\rfloor+1}^n a_{1j}^n<(2+\varepsilon)n p_n\Big\}^{\lfloor \varepsilon n\rfloor}\\
&\leq \bigg(1-\frac{c}{\sqrt{\log n}}n^{c\varepsilon-1}\bigg)^{\lfloor \varepsilon n\rfloor}\stackrel{n\to\infty}{\longrightarrow}0,
\end{align*}
implying that
$$
\lim\limits_{n\to\infty}\Prob\bigg\{\frac{\rho_n}{\sqrt{n p_n}}\geq\frac{1}{\sqrt{1+\varepsilon}}+\sqrt{1+\varepsilon}\bigg\}=1.
$$
Together with Theorem~\ref{th: lower top}, this gives the result.
\end{Remark}

\bigskip
\bigskip

\noindent {\small Konstantin Tikhomirov,}\\
{\small School of Mathematics, GeorgiaTech,}\\
{\small E-mail: ktikhomirov6@gatech.edu}

\bigskip

\noindent {\small Pierre Youssef,}\\
{\small Laboratoire de Probabilit\'es, Statistique et Mod\'elisation, 
Universit\'e Paris Diderot,}\\
{\small E-mail: youssef@lpsm.paris}

\end{document}